\documentclass[12pt,psamsfonts,leqno,oneside,letterpaper]{amsart}
\usepackage{ifpdf}
\ifpdf
\usepackage[text={6.5truein,9truein},left=1truein,top=1truein]{geometry}
\usepackage[pdftex]{graphicx}
\else
\usepackage[dvips,text={6.5truein,9truein},left=1truein,top=1truein]{geometry}
\usepackage{graphicx}
\fi
\usepackage{amssymb,amsmath,amscd,enumerate}
\usepackage{url}

\usepackage[colorlinks,linkcolor=blue,citecolor=blue,pdfstartview=FitH]{hyperref}
\input xy
\xyoption{all}
\SelectTips{cm}{12}

\usepackage{color}

\newcommand\delete[1]{{\color{green}#1}}

\theoremstyle{definition}
\newtheorem{para}{}[section]
\newtheorem{subpara}{}[para]
\newtheorem{remark}[para]{Remark}
\newtheorem{remarks}[para]{Remarks}
\newtheorem{notation}[para]{Notation}
\newtheorem{convention}[para]{Convention}
\newtheorem{definition}[para]{Definition}
\newtheorem{definitions}[para]{Definitions}

\newcommand\Alternatives{\begin{enumerate}[(i)]}
\newcommand\EndAlternatives{\end{enumerate}}
\newcommand\Conditions{\begin{enumerate}[(1)]}
\newcommand\EndConditions{\end{enumerate}}

\theoremstyle{plain}
\newtheorem{theorem}[para]{Theorem}
\newtheorem{lemma}[para]{Lemma}
\newtheorem{proposition}[para]{Proposition}
\newtheorem{corollary}[para]{Corollary}
\newtheorem{conjecture}[para]{Conjecture}
\newtheorem{claim}[subpara]{}

\numberwithin{equation}{para}
\numberwithin{figure}{section}

\newcommand\Number{\begin{para}}
\newcommand\EndNumber{\end{para}}
\newcommand\Definition{\begin{definition}}
\newcommand\EndDefinition{\end{definition}}
\newcommand\Definitions{\begin{definitions}}
\newcommand\EndDefinitions{\end{definitions}}
\newcommand\Theorem{\begin{theorem}}
\newcommand\EndTheorem{\end{theorem}}
\newcommand\Conjecture{\begin{conjecture}}
\newcommand\EndConjecture{\end{conjecture}}
\newcommand\Remark{\begin{remark}}
\newcommand\EndRemark{\end{remark}}
\newcommand\Remarks{\begin{remarks}}
\newcommand\EndRemarks{\end{remarks}}
\newcommand\Convention{\begin{convention}}
\newcommand\EndConvention{\end{convention}}
\newcommand\Notation{\begin{notation}}
\newcommand\EndNotation{\end{notation}}
\newcommand\Lemma{\begin{lemma}}
\newcommand\EndLemma{\end{lemma}}
\newcommand\Proposition{\begin{proposition}}
\newcommand\EndProposition{\end{proposition}}
\newcommand\Corollary{\begin{corollary}}
\newcommand\EndCorollary{\end{corollary}}
\newcommand\Claim{\begin{claim}}
\newcommand\EndClaim{\end{claim}}
\newcommand\Proof{\begin{proof}}
\newcommand\EndProof{\end{proof}}
\newcommand\Equation{\begin{equation}}
\newcommand\EndEquation{\end{equation}}
\newcommand\NoProof{{\hfill$\square$}}
\newcommand\Bullets{\begin{itemize}}
\newcommand\EndBullets{\end{itemize}}

\newcommand\threefreevolume{9.3}

\newcommand\bigradcor{6.2}
\newcommand\threefreevol{9.3}

\newcommand\ifnotwhynot{7.1}
\newcommand\firstAST{9.4}
\newcommand\nonfibroidsection{6}
\newcommand\geomsix{1.2}

\newcommand\Vbor{V_{\text{B\"or}}}
\newcommand\semithick{semithick}

\newcommand\nameit{\mathfrak N}
\newcommand\Nameit{\mathfrak N}

\newcommand\Vnearnought{V_{\rm near}^*}
\newcommand\mnearnought{m_{\rm near}^*}
\newcommand\Vfar{V_{\rm far}}

\newcommand\Vnear{V_{\rm near}}
\newcommand\mnear{m_{\rm near}}

\newcommand\lc{labeling-compatible{}}

\newcommand\hyph{-}
\newcommand\inter{\mathop{\rm int}}
\newcommand\cosech{\mathop{\rm cosech}}
\newcommand\st{\mathop{\rm st}}
\newcommand\lk{\mathop{\rm lk}}

\newcommand\Mthick{M_{\rm thick}}
\newcommand\Mthin{M_{\rm thin}}

\newcommand\tc{\widetilde c}

\newcommand\tY{\widetilde Y}
\newcommand\tP{\widetilde P}

\newcommand\calc{{\mathcal C}}
\newcommand\cale{{\mathcal E}}
\newcommand\cald{{\mathcal D}}
\newcommand\calL{{\mathcal L}}
\newcommand\calk{{\mathcal K}}
\newcommand\calt{{\mathcal T}}
\newcommand\calr{{\mathcal R}}
\newcommand\calb{{\mathcal B}}
\newcommand\calg{{\mathcal G}}
\newcommand\cali{{\mathcal I}}

\newcommand\calu{{\mathcal U}}
\newcommand\calx{{\mathcal X}}
\newcommand\caly{{\mathcal Y}}
\newcommand\calw{{\mathcal W}}
\newcommand\calv{{\mathcal V}}

\newcommand\Isom{{\rm Isom}}
\newcommand\arccosh{\mathop{\rm arccosh}}

\newcommand\arcsech{{\rm arcsech}}
\newcommand\arcsinh{{\rm arcsinh}}
\newcommand\sech{{\rm sech}}
\newcommand\arcsec{{\rm arcsec}}
\newcommand\ZZ{{\mathbb Z}}
\newcommand\Z{{\mathbb Z}}

\newcommand\RR{{\mathbb R}}
\newcommand\EE{{\mathbb E}}

\newcommand\HH{{\mathbb H}}
\newcommand\UU{{\mathbb U}}
\newcommand\XX{{\mathfrak G}}
\newcommand\ff{{\mathfrak s}}

\newcommand\cals{{\mathcal S}}
\newcommand\calz{{\mathcal Z}}
\newcommand\dist{\mathop{\rm dist}}

\newcommand\vol{\mathop{\rm vol}}

\newcommand\length{\mathop{{\rm length}}}
\newcommand\rank{\mathop{{\rm rank}}}

\newcommand\isomplus{\mathop{{\rm Isom}_+}}

\newcommand{\D}{w}
\newcommand{\radone}{u}
\newcommand{\radtwo}{v}

\newcommand\whatvol{3.44}
\newcommand\whatslightlybiggervol{3.468}
\newcommand\whatMargulis{1.119}
\newcommand\whatfunnyradius{0.444}
\newcommand\whatfunnyvolume{0.381}
\newcommand\whatelse{1.392}
\newcommand\twicewhatelse{2.784}
\newcommand\fricewhatelse{5.568}

\begin{document}

\author{Marc Culler}
\address{Department of Mathematics, Statistics, and Computer Science (M/C 249)\\
University of Illinois at Chicago\\
851 S. Morgan St.\\
Chicago, IL 60607-7045}
\email{culler@math.uic.edu}

\author{Peter B. Shalen}
\address{Department of Mathematics, Statistics, and Computer Science (M/C 249)\\
University of Illinois at Chicago\\
851 S. Morgan St.\\
Chicago, IL 60607-7045}
\email{shalen@math.uic.edu}
\thanks{Partially supported by NSF {grants DMS-0204142 and DMS-0504975}}

\title{Four-free groups and hyperbolic geometry}

\begin{abstract}
  We give new information about the geometry of closed, orientable
  hyperbolic $3$-manifolds with $4$-free fundamental group. As an
  application we show that such a manifold has volume greater than
  $\whatvol$. This is in turn used to show that if $M$ is a closed
  orientable hyperbolic $3$-manifold such that $\vol M \le\whatvol$, then
  $\dim_{\Z_2}H_1(M;\Z_2) \le 7$.
\end{abstract}

\maketitle

\section{Introduction}

The theme of this paper is the interaction between the geometric
structure of hyperbolic $3$-manifolds and their topologically defined
properties. While the basic techniques of the paper are suggested by
those of \cite{log5}
and Section 9 of
\cite{accs}, the present paper involves deeper topological results and
much more subtle combinatorial ideas than those used in \cite{log5} or
\cite{accs}.

The following definitions give a context for some of our main results:

\Definitions\label{k-free def}
The {\it rank} of a finitely generated group is defined to be the
minimal cardinality of a generating set of the group.  A group
$\Gamma$ is said to be {\it $k$-free}, where $k$ is a given positive
integer, if every finitely generated subgroup of $\Gamma$ having rank
at most $k$ is free.
\EndDefinitions

The property of having $k$-free fundamental group can often be deduced
from natural conditions on the topology of a closed, orientable
hyperbolic $3$-manifold $M$. For example, according to a theorem
proved by Jaco and Shalen \cite[Theorem VI.4.1]{JS}, $\pi_1(M)$ is
always $2$-free unless $M$ has a finite cover $\widetilde M$ such that
$\pi_1(\widetilde M)$ has rank $2$. As a second example, it follows
from \cite[Proposition \ifnotwhynot]{lastplusone} that if $k\ge3$, and
if $H_1(M;\ZZ_2)$ has $\ZZ_2$-dimension at least $\max(3k-4,6)$, then
either $\pi_1(M)$ is $k$-free, or $M$ contains a closed incompressible
surface of genus at most $k-1$.

On the other hand, via the ``$\log(2k-1)$-Theorem'' proved in
\cite{paradoxical} and \cite{accs}, together with the Marden
conjecture proved in \cite{agol} and in \cite{cg}, the property of
being $k$-free interfaces directly with the geometry of the
manifold. For example, according to \cite[Corollary 4.2]{rankfour},
which is deduced from results in \cite{paradoxical}, \cite{agol} and
\cite{cg}, if $\pi_1(M)$ is $2$-free then $\log3$ is a ``Margulis
number'' for $M$ in a sense that will be reviewed below in Subsection
\ref{braes heights}. In particular this implies that there exists a
point of $M$ where the injectivity radius is at least $(\log 3)/2$.
Likewise, it follows from Corollary \threefreevolume\ of \cite{last},
that a closed orientable hyperbolic 3\hyph manifold with $3$-free
fundamental group contains a point where the injectivity radius is
$(\log5)/2$.

The methods used to prove these results in the $2$-free and $3$-free
cases depend not only on the $\log(2k-1)$ theorem but on subtle
topological and combinatorial arguments involving certain coverings of
$\HH^3$ by cylinders.  These arguments do not extend directly to give
stronger information in the $4$-free case.

In this paper we introduce much more subtle topological and
combinatorial methods which lead to a fundamental new fact
about the geometry of a $3$-manifold with $4$-free fundamental group,
which is stated as Theorem \ref{Real Key theorem} below. This result
will in turn be applied to prove a volume estimate for such manifolds,
which is stated as Theorem \ref{the big one} below. Furthermore,
combining this with the results of \cite{lastplusone} we obtain a new
connection between volume and homology, stated in Theorem \ref{sweet
  land of flub-a-dub} below.

The statement of Theorem \ref{Real Key theorem} requires a definition.

\Definition\label{bernard mergendeiler}
We shall say that a point $P$ of a hyperbolic $3$-manifold $ M$ is
{\it $\lambda$-thin}, where $\lambda$ is a given positive number, if
there is a homotopically non-trivial loop of length $<\lambda$ based
at $P$.  We will say that a point is {\it $\lambda$-thick} if
it is not $\lambda$-thin.
\EndDefinition

Equivalently, a point of $M$ is $\lambda$-thick if and only if it is
the center of a hyperbolic ball of radius $\lambda/2$ in $M$. Thus
Corollary \threefreevolume\ of \cite{last} is equivalent to the
assertion that if $M$ is closed and orientable and $\pi_1(M)$ is
$3$-free, then $M$ contains a $(\log5)$-thick point.

\Definition\label{cockeyed olives} We shall say that a point $P\in M$
is {\it $\lambda$-doubly thin} if there are two loops of length
$<\lambda$ based at $P$ which represent non-commuting elements of
$\pi_1(M,P)$.  A point $P\in M$ will be said to be {\it
  $\lambda$-\semithick} if it is not $\lambda$-doubly thin.
\EndDefinition

Section \ref{short'nin' bread} contains a general discussion of
properties of doubly-thin and \semithick\ points.

We can now state the main result of this paper, which is proved 
In Section \ref{Key  section}.

\Theorem\label{Real Key theorem} Suppose that $M$ is a closed, orientable
hyperbolic $3$-manifold such that $\pi_1(M)$ is $4$-free. Then  $M$
has a $\log7$-semithick point.
\EndTheorem

The proof begins from the same basic point of view as that of
\cite[Corollary \threefreevolume]{last}, but requires much deeper
arguments. In the following sketch we will use terms that are
well-known and are defined precisely in the body of the paper.

\Number\label{braes heights}
Consider a closed, orientable hyperbolic $3$-manifold written as
$M=\HH^3/\Gamma$, where $\Gamma\le\Isom_+(\HH^3)$ is discrete and
cocompact.  Fix a positive number $\lambda$.  For each maximal cyclic
subgroup $C$ of $\Gamma$ generated by a (loxodromic) element of
translation length $<\lambda$, we consider the hyperbolic cylinder
consisting of points of $\HH^3$ that are displaced through a distance
$<\lambda$ by some non-trivial element of $C$. This gives a family
$\calz = \{Z_C\}$ of cylinders indexed by certain maximal cyclic
subgroups.  This family of cylinders covers $\HH^3$ if and only if
the manifold $M$ contains no $\lambda$-thick point.

If $\Gamma$ is $2$-free then any pair of distinct maximal cyclic subgroups
generates a free group of rank $2$.  Thus, if we take $\lambda=\log 3$,
then it follows from the ``$\log3$-theorem''---the case
$k=2$ of the $\log(2k-1)$-Theorem---that
$Z_C\cap Z_{C'}=\emptyset$ whenever $C \not= C'$. This is
expressed by saying that $\log3$ is a Margulis number for $M$; see
Subsection \ref{rev} below.

Since hyperbolic space cannot be covered by pairwise disjoint open
cylinders, this shows that $M$ has a $\log 3$-thick point.

When $\Gamma$ is $k$-free, for $k>2$, this argument can be refined by
assuming that $\HH^3$ is covered by the cylinders in the family
$\calz$ and considering the nerve of this covering.  The nerve is an
abstract simplical complex $K$ whose associated space is denoted
$|K|$.  By definition each vertex of $K$ corresponds to an index for
the family $\calz$, i.e. a certain maximal cyclic subgroup of
$\Gamma$.  To each open simplex $\sigma$ of $K$ we assign the subgroup
$\Theta(\sigma) \le \Gamma$ which is generated by the cyclic subgroups
corresponding to the vertices of $\sigma$.  If we take $\lambda =
\log(2k-1)$ then the $\log(2k-1)$-Theorem, together with the
  $k$-free assumption and the assumption that $\calz$ covers, would
  give a contradiction if there existed a simplex $\sigma$ of
  dimension $k-1$ with $\Theta(\sigma) = k$.  Thus, for every simplex
  $\sigma$ of dimension $k-1$, $\Theta(\sigma)$ must be a free group
  of rank at most $k-1$.
For $k=3$, the group-theoretic arguments given in \cite{accs} use this
structure to derive a contradiction from the assumption that $\calz$
covers $\HH^3$, implying that $M$ has a $\log 5$-thick point.

In the case being considered in this paper, where $\Gamma$ is
$4$-free, we argue by contradiction and assume that $M$ has no
  $\log7$-semithick point. This implies that the family $\calz$
``doubly covers'' $\HH^3$ in the sense that every point lies in at
least two different cylinders in the family.  The ``double covering''
property, together with the contractibility of $\HH^3$, implies that
the complement of the $0$-skeleton $|K^{0}|$ relative to the
$3$-skeleton $|K^{3}|$ is connected and simply connected.

We know that $\Theta(\sigma)$ must be a free group of rank $2$ or $3$
for any open simplex $\sigma$ contained in $|K^3|-|K^0|$.
Set-theoretically we may therefore regard $|K^3|-|K^0|$ as the
disjoint union of sets $X_2$ and $X_3$, where $X_k$ is the union of
all open simplices $\sigma$ of $K_3$ for which $\Theta(\sigma)$ has
rank $k$.

As in the corresponding argument in \cite{accs}, we have a natural
simplicial action of $\Gamma$ on $K$. For $k=2,3$, this induces an
action on the set of connected components of $X_k$. A key step is
showing that the stabilizer of any component of $X_k$ under this
action is a free group. This is approached by the same basic ideas
involving the lattice of free subgroups in a $k$-free group that was
used in \cite{accs}, but it is much more difficult. In particular it
depends in a crucial way on the result recently proved by
R. Kent, \cite{kent} and independently by L. Louder and
D. B. McReynolds, \cite{l-mcr} that if two rank-$2$ subgroups of a
free group have a rank-$2$ intersection then they have a
rank-$2$ join. The relevant group theory, incorporating our
application of the Kent-Louder-McReynolds result, is done in Section
\ref{kentish section}.

The actions of $\Gamma$ on the sets of components of the $X_k$ give
rise to an action of $\Gamma$ on an abstract bipartite graph $T$. The
vertices of $T$ are the components of $X_2$ and of $X_3$. Two vertices
of $T$ are joined by an edge if they correspond to subsets of
$|K^3|-|K^0|$, one of which is a component of $X_2$ and one a
component of $X_3$, and if some simplex contained in one of these sets
is a face of a simplex contained in one of the other sets.

The $1$-connectedness of $|K^3|-|K^0|$ implies that $T$ is a tree. The
group $\Gamma$ acts on $T$ without inversions.  The vertex stabilizers
under the action of $\Gamma$ on $T$ are stabilizers of components of
the $X_k$ under the action of $\Gamma$, and are therefore free. As
$\Gamma$ is isomorphic to the fundamental group of a closed hyperbolic
$3$-manifold, basic facts about actions of $3$-manifold groups on
trees, which we quote from \cite{splittings}, then lead to a
contradiction, and Theorem \ref{Real Key theorem} is proved.
\EndNumber

Having established that interesting geometric properties of $M$ follow
from the assumption that $\pi_1(M)$ is $4$-free, it is natural to ask
whether these geometric conditions in fact imply stronger
volume estimates than the ones
which follow from the $2$-free and $3$-free hypotheses.  Sections
\ref{capsection} -- \ref{short section} address this question.
  
It is pointed out in \cite[Corollary 9.3]{last} that the existence of
a $(\log 5)$-thick point in a hyperbolic $3$-manifold $M$ with
$3$-free fundamental group leads to a lower bound of $3.08$ for the
volume of $M$.  As we explain below, in the case where $\pi_1(M)$ is
$4$-free, much more intricate methods are needed in order to pass 
from the existence of a $(\log 7)$-semithick point to an estimate for
the volume of $M$. The volume estimate which we eventually obtain with
these methods is summarized by the following theorem, which we prove
in Section \ref{short section}.

\Theorem\label{the big one} Let $M$ be a closed, orientable hyperbolic
$3$-manifold. If $\pi_1(M)$ is $4$-free then $\vol M>\whatvol$.
\EndTheorem

To compare the relative strengths of the bounds $3.08$ and $\whatvol$,
we  note that it is a consequence of \cite[Corollary
  6.6.3]{thurstonnotes} that the set $\calv$ of all volumes of closed,
  orientable hyperbolic $3$-manifolds is a well-ordered subset of
  $\RR$. There are $34$ known volumes of cusped orientable hyperbolic
  $3$-manifolds between $3.08$ and $\whatvol$; this can be shown to
  imply that the ordinal type of the set $\calv\cap(3.08,\whatvol)$ is
  at least $34\omega$. By contrast, there are only eight known volumes
  of cusped orientable hyperbolic $3$-manifolds less than $3.08$, and
  the best available lower bound for the ordinal type of the set
  $\calv\cap(0,3.08)$ is $8\omega+6$. Thus from the point of view of
  ordinal numbers, the lower bound of $3.44$ may be regarded as being
  more than four times stronger than the lower bound of $3.08$.

Section \ref{short section} also contains a proof of the following
result, which establishes a new connection between volume and
homology:

\Theorem\label{sweet land of flub-a-dub}
Let $M$ be a closed orientable hyperbolic $3$-manifold such that
$\vol M \le\whatvol$. Then $\dim_{\Z_2}H_1(M;\Z_2) \le 7$.
\EndTheorem

This result is analogous to Theorem \geomsix\ of \cite{lastplusone},
which asserts that if $\vol M \le3.08$ then $\dim_{\Z_2}H_1(M;\Z_2)
\le 5$.

\Number\label{jumabalaya}
We now sketch the proof of Theorem \ref{the
  big one}.  In the case where $M$ contains a $\log7$-thick point,
the sphere-packing results established in \cite{boroczky} give a lower
bound of more than $5.7$ for $\vol M$. If
$M$ contains a $\log7$-\semithick\ point but contains no $\log7$-thick
point, a continuity argument produces a point $P\in M$ which is 
$\log7$-\semithick, but is not $\lambda$-\semithick\ for any
$\lambda>\log7$. All loops based at the point $P$ which have length
$<\log7$ represent elements of a single maximal cyclic subgroup $C$ of
$\pi_1(M,P)$. (The relevant properties of the point $P$ are
summarized in Corollary \ref{Key theorem}, in somewhat different and
more technically convenient language.)

If $P$ is such a point, let $N\subset M$ denote the metric
neighborhood of $P$ with radius $(\log7)/2$, and let $D$ denote the
minimal length of a loop based at $P$ which represents a generator of
the maximal cyclic group $C$. Lemma \ref{nested cubes} gives a lower
bound for $\vol N$ in terms of $D$, subject to a suitable lower bound
$\delta$ on the length of the shortest geodesic in $M$.

To establish this lower bound one begins by interpreting 
$\vol N$ as the volume of a set $X$ obtained by removing finitely many caps
from a ball $B\subset\HH^3$ of radius $(\log7)/2$ (see Proposition
\ref{c.d. rivington}).  Each of these caps is associated with a
non-trivial element $\gamma$ of $C$. If we write $M=\HH^3/\Gamma$ as
above, take $B$ to be centered at a point $\tP$ in the pre-image of $P$
under the quotient map,
and identify $C$ with a subgroup of the deck transformation group,
then the cap corresponding to $\gamma$ is one ``half'' of the
intersection of $B$ with $\gamma\cdot B$. Thus an element $\gamma$ of $C$
corresponds to a cap in the construction of $X$ only if
$B\cap\gamma\cdot B\not=\emptyset$.  If $P$ happens to lie
close enough to the closed geodesic corresponding to the cyclic group
$C$ then the set $X$ is obtained by removing two disjoint caps
from the ball; in general, however, there may be more than two caps,
and some of them may overlap. 

Imposing a lower bound on the shortest geodesic in $M$ gives a lower
bound on the translation length of a generator of $C$. When this
translation length is sufficiently large, only relatively small powers
of the generator of $C$ can correspond to caps in the construction of
$X$. Furthermore, the lower bound for the translation length of the
generator gives lower bounds for the displacement of $\tP$ under these
powers of the generator, and this controls the volumes of the caps and
of their intersections. The details involve quite a bit of hyperbolic
trigonometry.

In order to apply Lemma \ref{nested cubes} for numerical estimates we
need to be able to calculate volumes of caps and of intersections of
caps. This  is the subject of the Appendix. 

The lower bound for $\vol N$ provided by Lemma \ref{nested cubes} is
in particular a lower bound for $\vol M$, but this lower bound
decreases as $D$ decreases, and by itself it turns out to be
insufficient to give the conclusion of Theorem \ref{the big one}
except for rather large values of $D$. To compensate for this we
obtain a lower bound for $\vol(M-N)$ which increases as $D$ decreases,
and use $\vol N+\vol(M-N)$ as a lower bound for $\vol M$.  To obtain a
lower bound for $\vol(M-N)$, we exploit the fact that there are
non-commuting elements of $\pi_1(M,P)$ represented by loops of length
$\log7$ and $D$, and again use the hypothesis that $\pi_1(M)$ is
$4$-free. These pieces of information are used, via results proved in
\cite{cusp} and adapted to the context of this paper in Section
\ref{distant point section}, to show that there is a point $Y\in M$
whose distance from $P$ is $\rho$, where the quantity $\rho$ is
explicitly defined as a monotonically decreasing function of $D$. The
results from \cite{cusp} also guarantee that we may take $Y$ to lie in
the $\mu$-thick part of $M$ if $\mu$ is any Margulis number for $M$.
(One could, for example, take $\mu$ to be $\log3$, which by
\cite[Corollary 4.2]{rankfour} is a Margulis number for any closed
orientable hyperbolic $3$-manifold with $2$-free fundamental group.)

When $\rho$ is sufficiently large it is easy to use the existence of
the point $Y$ to give a non-trivial lower bound for $\vol(M-N)$. For example, if $\rho>(\log7)/2$, it follows from the triangle
inequality that the $(\rho-(\log7)/2)$-neighborhood of $Y$ is contained in
$M-N$. If $(\rho-(\log7)/2)<\mu/2$, this neighborhood is a hyperbolic ball
and its volume, which can be calculated explicitly in terms of $\rho$,
is a lower bound for $\vol(M-N)$. If $(\rho-(\log7)/2)>\mu/2$, the
volume of $M-N$ is bounded below by the volume of a hyperbolic ball of
radius $\mu/2$.

In order to obtain the lower bound asserted in Theorem \ref{the big
  one}, this straightforward method for bounding $\vol(M-N)$ from
below must be refined in several ways. One of these involves the
choice of the Margulis constant $\mu$.  We mentioned above that
$\log3$ is a Margulis number for $M$ provided that $\pi_1(M)$ is
$2$-free. In Section \ref{Margulis section} we use the methods of
\cite{cusp} to give a stronger result when $\pi_1(M)$ is $k$-free for
a given $k>2$ and the diameter $\Delta$ of $M$ is known.  Corollary
\ref{biddy biddy bum bum to san fernando} asserts that a certain
quantity defined as a function of $k$ and $\Delta$, which is
monotonically increasing in $k$ and monotonically decreasing in
$\Delta$, is a Margulis number.  By combining this with
\cite[Corollary \threefreevolume]{last}, we show that if $\pi_1(M)$
were $4$-free and if $M$ had volume at most $\whatvol$ then $\mu=1.119$
would be a Margulis number for $M$. We may therefore use this improved choice
of $\mu$ in the proof of Theorem \ref{the big one}.

A second refinement of the straightforward lower bound
for $\vol (M-N)$ is based on the sphere-packing arguments given in
\cite{boroczky}.  If $(\rho-(\log7)/2)>h$, where $h$ denotes the
distance from the barycenter to a vertex of a regular hyberbolic
tetrahedron with sides of length $\mu$, then $M-N$ contains the metric
neighborhood $N'$ of radius $h$ about $Y$, and arguments in
\cite{boroczky} give an explicit lower bound for $\vol N'$ which is
significantly greater than the volume of a ball of radius $\mu/2$.
While this lower bound applies to any $\mu/2$-thick point in a
hyperbolic manifold, we have more information in the present
situation: not only is $Y$ a $\mu/2$-thick point, but $M$ contains the
point $P$ whose distance from $Y$ is $\rho$, a number which is often
considerably larger than $h$.  In Section \ref{Boroczky balls} we show
that the lower bounds given in \cite{boroczky} can be improved using
the existence of such a ``distant point'' $P$.

A third refinement of the lower bound for $\vol(M-N)$ is based on an
observation that was already used in \cite{cusp}.  When
$\rho-(\log7)/2>h$, there is a point $Y'$ of $M$ whose minimum
distances from both the $(\log7)/2$-neighborhood of $P$ and the
$h$-neighborhood of $Y$ are at least $\rho-((\log7)/2+h)/2$. It is
often possible to take such a point $Y'$ to be $\mu$-thick, and thus
to obtain an additional contribution to $\vol(M-N)$ from a suitable
metric neighborhood of $Y'$ in the same way as from a metric
neighborhood of $Y$.

We mentioned that our method for bounding $\vol N$ from below
requires a lower bound $\delta$ for the length of the shortest geodesic in
$M$. It turns out that if we take $\delta=0.58$, the methods
sketched above give the lower bound $\whatvol$ for 
$\vol M$. We have lower bounds for both $\vol N$ and $\vol(M-N)$ in  terms
of the parameter $D$, and their sum is a function of $D$ which is a
lower bound for $\vol M$, provided that no geodesic in $M$
has length $<0.58$. In Section \ref{Numbers} we prove by a rigorous
sampling argument that this function is bounded below by $\whatvol$ on the
relevant range.

In the case where $M$ does contain a short geodesic $c$ of length $l<0.58$,
the argument does not use Theorem \ref{Real Key theorem}, but it does
use many of the other ingredients described above. We fix a point $P$
on the closed geodesic $c$ and define $N\subset M$ to be the
$\lambda/2$-neighborhood of $P$, where $\lambda$ is the length of the
shortest loop based at $P$ that does not represent an element of the
cyclic subgroup of $\pi_1(M,P)$ determined by $c$. Since $P$ lies on
$c$, $\vol N$ is equal to the volume of a set obtained by
removing two disjoint caps from a ball of radius $\lambda/2$ in
$\HH^3$; the volumes of the caps are determined by $\lambda$ and the
length $l$ of $c$. Lemma 
\ref{nested
  cubes} again applies to give a lower bound for $\vol(M-N)$;
the quantity $\lambda$ plays the role that $\log7$ played in the
earlier application, and $l$ plays the role of $D$. We now have a
lower bound for $\vol M$ involving the two parameters
$\lambda$ and $l$. In Section \ref{short section} we prove by a rigorous
sampling argument that this function is bounded below by $\whatvol$ on the
relevant range of values of $l$ and $\lambda$, except possibly when
$l<.003$. In the latter case a different method, based on
estimates for tube radius established in \cite{accs} and tube-packing
estimates proved by Przeworski in \cite{prez}, gives the desired
lower bound of $\whatvol$ for $\vol M$.
\EndNumber

We are grateful to Rosemary Guzman for pointing out a number of errors
and clarifying some passages in earlier drafts of the paper.  We
are also grateful to Richard Kent, and to Lars Louder and Ben
  McReynolds, for proving their theorem in
\cite{kent} and \cite{l-mcr}. This theorem was an old conjecture, and we had
called attention to it after we realized its relevance to hyperbolic
geometry.

\section{General Conventions}

\Number\label{kurosh stuff}
Let $S$ be a subset of a group
$\Gamma$. We shall denote by $\langle S\rangle$ the subgroup of
$\Gamma$ generated by $S$. We shall say that $S$ is {\it independent}
(or that the elements of $S$ are independent) if $\langle S\rangle$ is
free on the generating set $S$.

It is a basic fact in the theory of free groups \cite[vol. 2, p. 59]{kurosh}
that a finite set $S\subset\Gamma$ is independent if and only if
$\langle S\rangle$ is free of rank $|S|$. This fact will often be used
without explicit mention.
\EndNumber

\Number\label{neither would you} For any positive real number $R$, let
us define the {\it cylinder of radius $R$} about a line $L$ in
$\HH^3$ to be the set of all points whose distance from $L$ is less
than $r$.
\EndNumber

\Number\label{monotonicity} A real-valued function $f(x_1,\ldots,x_n)$
of $n$ variables will be said to be {\it monotone increasing
(resp. decreasing) in the variable $x_i$} if
$f(x_1, \ldots, x_i, \ldots, x_n) \le
f(x_1, \ldots, x'_i, \ldots, x_n)$
whenever $x_i < x'_i$ (resp. $x_i > x'_i$) and 
both $(x_1, \ldots, x_i, \ldots, x_n)$ and
$(x_1, \ldots, x'_i, \ldots, x_n)$ lie in
the domain of $f$.
\EndNumber
 
\Number\label{omegastuff}
Let $\gamma$ be a loxodromic isometry of $\HH^3$. Let $l$ denote the
translation length of $\gamma$ and let $\theta$ denote its twist
angle. If $z\in\HH^3$ is a point, and if we set
$D=\dist(z,\gamma\cdot z)$, then the discussion in
  \cite[Subsection 1.3]{betti2} shows that the distance from $z$ to the axis of
$\gamma$ is equal to $\omega(l,\theta,D)$, where $\omega$ is the
function defined for $0\le l\le D$ and $\theta\in\RR$ by
\Equation\label{omega} \omega(l,\theta, D )\doteq \arcsinh\bigg(\bigg(
\frac{\cosh D -\cosh l}{\cosh l-\cos\theta}\bigg)^{1/2}\bigg).
\EndEquation The formula (\ref{omega}) shows that for any $\theta$ and
any $l>0$, the function $\omega(l,\theta,\cdot)$ is a continuous,
monotonically increasing function on $(l,\infty)$.

If $\lambda$ is any positive real number, we shall denote by
$Z_\lambda(\gamma)$ the set of points $ z\in\HH^3$ such that
$\dist(z,\gamma\cdot z)<\lambda$. Then $Z_\lambda(\gamma)$ is empty if the
translation length $l$ of $\gamma$ is at least $\lambda$. If
$l<\lambda$, then in view of the monotonicity of $\omega$ in the third
variable, $Z_\lambda(\gamma)$ is a cylinder  of radius
$\omega(l,\theta,\lambda)$ about the axis of $\gamma$.
Furthermore, the continuity of $\omega$ in the third variable implies
that if $l<\lambda$ then
\Equation\label{weak}
\overline{Z_\lambda(\gamma)}=\{z\in\HH^3:\dist(z,\gamma\cdot z)\le\lambda\}.
\EndEquation
\EndNumber

\Number\label{mystic caravan}
Let $C$ be a cyclic subgroup of $\isomplus(\HH^3)$ generated by a
loxodromic element $\gamma_0$.  The non-trivial elements of $C$
have a common axis which we shall denote by $A_C$. Now if $\lambda$ is
any positive real number, it follows from the discussion in
\ref{omegastuff} that the set
$Z_\lambda(C)\doteq\bigcup_{1\ne\gamma\in C}Z_\lambda(\gamma)$ is
empty if $\gamma_0$ has translation length $\ge\lambda$, and is a
cylinder about $A_C$ if $\gamma_0$ has translation length
$<\lambda$.  If $\gamma_0$ has length $l$ and twist angle $\theta$
then the radius of the cylinder $Z_\lambda(C)$ is
\Equation\label{Omega}
\Omega_C(\lambda)\doteq\max_{1\le n\le[\lambda/l] }\omega(nl, n\theta, \lambda).
\EndEquation
\EndNumber

\Number\label{same to me}
If $\Gamma$ is a discrete subgroup of $\isomplus(\HH^3)$, we shall
denote the quotient projection  $\HH^3\to\HH^3/\Gamma$ by
$q_\Gamma$. If $P$ is a point of $M=\HH^3/\Gamma$, each point $\tP$ of
$q_\Gamma^{-1}(P)$ determines an isomorphic identification of $\pi_1(M,P)$
with $\Gamma$.
\EndNumber

\Number\label{extrinsic}
 If $(X,d)$ is a metric space (for example a hyperbolic manifold
with the usual distance function) and $A$ is a bounded subset of $X$,
we define the {\it extrinsic diameter} of $A$ in $X$ to be the
quantity $\sup_{x,y\in A}d(x,y)$.

If $r$ is a non-negative real number and $P$ is a point in $X$, we
shall let
$$N_X(P,r) \dot= \{x \in X \;|\; d(x,P) < r\}$$
denote the metric ball with radius $r$ and center $P$.  We will
abbreviate this as $N(P,r)$ when it is clear which metric space
$(X,d)$ is meant.   Note that we are using the term ``ball'' to
mean an open ball.
\EndNumber

\Number\label{ti creves quand meme}
Let $M$ be a closed, orientable hyperbolic $3$-manifold, and let
$\lambda$ be a positive number. The notion of a $\lambda$-thin point
of $M$ was defined in \ref{bernard mergendeiler}. We denote the set of
$\lambda$-thin points of $M$ by $\Mthin(\lambda)$, and we set
$\Mthick(\lambda)=M-\Mthin(\lambda)$.

Suppose we write $M=\HH^3/\Gamma$, where $\Gamma$ is a discrete
subgroup of $\isomplus(\HH^3)$. Let $P$ be a point of $M$ and let
$\tP$ be a point of $q_\Gamma^{-1}(P)$. Then $P$ is a $\lambda$-thin
point if and only if we have $\dist_{\HH^3}(\gamma\cdot\tP,\tP)<\lambda$
for some $\gamma\in\Gamma-\{1\}$. Equivalently, we have
$P\in\Mthin(\lambda)$ if and only if $P$ lies in $Z_\lambda(C)$ for some
maximal cyclic subgroup $C$ of $\Gamma$.
\EndNumber

\section{Doubly-thin and \semithick\ points}\label{short'nin' bread}

The notion of a $\lambda$-doubly thin or $\lambda$-\semithick\
point of a hyperbolic manifold, where $\lambda$ is a given positive
number, was defined in \ref{cockeyed olives}.

\Proposition\label{the buck stops here} Let $M$ be a closed,
orientable hyperbolic $3$-manifold. For every point $P\in M$ there is
a unique number $\nameit>0$ such that $P$ is $\lambda$-doubly thin
for every $\lambda>\nameit$ and is $\lambda$-semithick for every
$\lambda$ with $0<\lambda\le\nameit$.
\EndProposition

\Proof Let $\cals\subset\pi_1(M,P)\times\pi_1(M,P)$ denote the set of
all non-commuting pairs of elements of $\pi_1(M,P)$. Since $M$ is
closed, $\pi_1(M,P)$ is non-abelian and hence $\cals\ne\emptyset$. For
each $g\in\pi_1(M,P)$, there is a unique loop $\alpha_g$ of minimal
length among all representatives of the based homotopy class $g$. The
set $\calL=\{\length\alpha_g:g\in\pi_1(M,P)\}$ is discrete, and hence
the set
$\calL'=\{\max(\length\alpha_g,\length\alpha_h):(g,h)\in\cals\}\subset\calL$
has a least element $\nameit$. We have $\nameit>0$ because
non-commuting elements of $\pi_1(M,P)$ must be non-trivial. It is
immediate from the definitions that $P$ is $\lambda$-doubly thin for
every $\lambda>\nameit$ and is $\lambda$-semithick for every
$\lambda$ with $0<\lambda\le\nameit$. Uniqueness is obvious.
\EndProof

\Notation\label{you asked for it} If $M$ is a closed, orientable
hyperbolic $3$-manifold, then for every point $P\in M$ we shall denote
by $\Nameit_M(P)$ the number $\nameit$ given by Proposition \ref{the
  buck stops here}. Thus $\Nameit_M$ is a positive-valued function
defined on $M$.
\EndNotation

\Proposition\label{you can't prove it ain't so} If $M$ is a closed,
orientable hyperbolic $3$-manifold, the function $\Nameit_M$ is
$2$-Lipschitz; that is,
$$|\Nameit_M(P)-\Nameit_M(Q)|\le2\dist_M(P,Q)$$
for all $P,Q\in M$. In particular, $\Nameit_M$ is continuous.
\EndProposition

\Proof By symmetry it is enough to show that
$\Nameit_M(Q)-\Nameit_M(P)\le2\dist_M(P,Q)$ for all points $P,Q\in M$.
Set $d=\dist_M(P,Q)$ and $\lambda=\Nameit(P)$. Then $P$ is
$\lambda$-doubly thin, so that there are loops $\alpha$ and
$\beta$ based at $P$, both of length $<\lambda$, and representing
non-commuting elements of $\pi_1(M,P)$. Let $\zeta$ be a path of
length $d$ from $Q$ to $P$. Then $\zeta\star\alpha\star\bar\zeta$ and
$\zeta\star\beta\star\bar\zeta$ have length $<\lambda+2d$ and
representing non-commuting elements of $\pi_1(M,Q)$. It follows that
$Q$ is $(\lambda+2d)$-doubly thin, i.e. that
$\Nameit_M(Q)\le\lambda+2d=\Nameit_M(P)+2\dist_M(P,Q)$. \EndProof

\Number\label{a baloney shell}
If $M$ is a closed, orientable hyperbolic $3$-manifold, we may write
$M=\HH^3/\Gamma$ where $\Gamma\le\Isom_+(\HH^3)$ is discrete and
torsion-free. Since $M$ is closed, $\Gamma$ is purely
loxodromic. Hence each non-trivial element $\gamma$ of $\Gamma$ lies
in a unique maximal cyclic subgroup, which is the centralizer of
$\gamma$. In particular, non-trivial elements which lie in distinct
maximal cyclic subgroups do not commute.

We shall denote by $\calc(\Gamma)=\calc(M)$ the set of all maximal
cyclic subgroups of $\Gamma$.  If $\lambda$ is a positive number,
we denote by $\calc_\lambda(\Gamma)=\calc_\lambda(M)$ the subset
consisting of all $C\in\calc(\Gamma)$ such that a generator of $C$ has
translation length $<\lambda$. It follows from \ref{mystic caravan}
that $Z_\lambda(C)$ is a cylinder of radius $\Omega_C(\lambda)$ if
$C\in\calc_\lambda(\Gamma)$ and is empty if
$C\in\calc(\Gamma)-\calc_\lambda(\Gamma)$.  

The discreteness of the group $\Gamma$ implies that the family $(
Z_\lambda(\gamma))_{1\ne\gamma\in\Gamma}$ is locally finite. Since for
each $C\in\calc_\lambda(\Gamma)$ we have $Z_\lambda(C)=Z_\lambda(\gamma)$ for
some $\gamma\in C-\{1\}$, the family
$(Z_\lambda(C))_{C\in\calc_\lambda(\Gamma)}$ is also locally finite.
\EndNumber

\Proposition\label{hootchy-kootcher} 
Let $M=\HH^3/\Gamma$ be a closed, orientable hyperbolic $3$-manifold,
let $P$ be a point of $M$, and let $\tP$ be a point of
$q_\Gamma^{-1}(P)$, and let $\lambda$ be a positive number. Then $P$
is a $\lambda$-doubly thin point of $M$ if and only if $\tP$ belongs
to the set
$$\cald\ \doteq\bigcup_{\begin{matrix}C,C'\in\calc_\lambda(\Gamma)\cr C\ne
    C'\end{matrix}}  Z_\lambda(C)\cap Z_\lambda(C').$$
\EndProposition

\Proof
We use the point $\tP$ to identify $\pi_1(M,P)$ with $\Gamma$ as in
\ref{same to me}.  If $\tP\in\cald$, there exist
$C,C'\in\calc_\lambda(\Gamma)$, with $C\ne C'$, such that $\tP$ lies
in both $Z_\lambda(C)$ and $Z_\lambda(C')$. Hence for some $x\in
C-\{1\}$ and $x'\in C'-\{1\}$ we have $\dist(x\cdot\tP,\tP)<\lambda$
and $\dist(x'\cdot\tP,\tP)<\lambda$. Thus $x,x'\in\pi_1(M,P)$ are
represented by loops of length $<\lambda$. Since $C$ and $C'$ are
distinct maximal cyclic subgroups of $\Gamma$, the elements $x$ and
$x'$ do not commute. By definition it follows that $P$ is
$\lambda$-doubly thin.

Conversely, suppose that $P$ is $\lambda$-doubly thin, so that
there are non-commuting elements $x,x'\in\pi_1(M,P)$ represented by
loops of length $<\lambda$. If we regard $x$ and $x'$ as elements of
$\Gamma$, they lie in distinct maximal cyclic subgroups
$C,C'\in\calc_\lambda(\Gamma)$. We have $\dist(x\cdot\tP,\tP)<\lambda$
and $\dist(x'\cdot\tP,\tP)<\lambda$, so that $\tP$ lies in both
$Z_\lambda(C)$ and $Z_\lambda(C')$ and hence $\tP\in\cald$.
\EndProof

\Number\label{rev} 
As in \cite{cusp}, we define a {\it Margulis number} for a closed,
orientable hyperbolic $3$-manifold $M$ to be a positive number $\mu$
such that for any two distinct subgroups $C,C'\in\calc_\mu(\Gamma)$ we
have $Z_\mu(C)\cap Z_\mu(C')=\emptyset$.  If $\mu$ is a Margulis
number for $M$, the components of $\Mthin(\mu)$ are tubes. In
particular, $\Mthick(\mu)$ is connected and non-empty.
\EndNumber

\Proposition\label{ti bouffes ti bouffes pas}
Let $M$ be a closed, orientable hyperbolic $3$-manifold, and let $\mu$
be a positive number. Then the following conditions are equivalent:
\Conditions
\item \label{griffoni non e buono} $\mu>0$ is a Margulis number for $M$;
\item \label{means griffoni ain't too good}$M$ has no $\mu$-doubly thin points;
\item \label{it's clearly understood} $\mu$ is a
lower bound for the function $\Nameit_M$.
\EndConditions
\EndProposition

\Proof 
The equivalence of (\ref{griffoni non e buono}) and (\ref{means
  griffoni ain't too good}) follows from 
Proposition \ref{hootchy-kootcher}.
The equivalence of
(\ref{means griffoni ain't too good}) and (\ref{it's clearly
  understood}) follows from Proposition \ref{the buck
  stops here} and the definition of $\Nameit_M$. 
\EndProof

\Number\label{eggs and fries}
Let $M$ be a closed, orientable hyperbolic $3$-manifold. For every
point $P$ of $M$, we shall denote by $\ell_P$ the smallest length of
any homotopically non-trivial loop in $M$ based at $P$.  The subgroup
of $\pi_1(M,P)$ generated by all homotopy classes that contain loops
of length $\ell_P$ will be denoted $\calL_P$.

We shall denote by $\XX_M$ the set of all points $P\in M$ such that
$\calL_P$ is cyclic.  If $P\in \XX_M$, then there is a unique maximal
cyclic subgroup containing $\calL_P$; we shall denote this subgroup by
$C_P$.

We shall define real-valued functions $D_M$ and $\ff_M$ with domain
$\XX_M$ as follows. For any $P\in\XX_M$ we define $D_M(P)$ to be the
minimal length of a loop based at $P$ which represents a generator of
$C_P$.  (Note that $D_M(P)$ is not necessarily equal to $\ell_P$,
since the shortest homotopically non-trivial loop based at $P$ may
represent a proper power of a generator of $C_P$.)

For any $P\in\XX_M$ we define $\ff_M(P)$ to be the smallest length of
any loop in $M$ based at $P$ which does {\it not } represent an
element of the cyclic group $C_P$.  In particular we have
$\ff_M(P)>\ell_P$.\EndNumber

\Proposition\label{you can always telephone girl}
Let $P$ be a point of a closed hyperbolic $3$-manifold $M$, and set
$\nameit=\Nameit_M(P)$. Then either
\Alternatives
\item $P$ is an $\nameit$-thick point of $M$, or
\item $P\in\XX_M$ and $\ff_M(P)=\nameit$.
\EndAlternatives
\EndProposition

\Proof
Suppose that $P$ is not an $\nameit$-thick point of $M$. Then by
definition there is a homotopically non-trivial loop of length
$<\nameit$ based at $P$. In particular, if $\alpha$ is the shortest
homotopically non-trivial loop based at $P$, we have
$\ell_P=\length\alpha<\nameit$. If $P$ did not lie in $\XX_M$, there
would be two loops based at $P$ having length $\ell_P$ and
representing non-commuting elements of $\pi_1(M,P)$. Hence $P$ would
be a $\lambda$-doubly thin point of $M$ for every $\lambda > \ell_p$.
This is impossible because
$\ell_P=\length\alpha<\nameit=\Nameit_M(P)$. Hence $P\in\XX_M$.

Now by the definition of $\ff_M(P)$, there is a loop $\beta$ of length
$\ff_M(P)$ based at $P$ such that $[\beta]\notin C_P$. It follows that
$[\alpha]$ and $[\beta]$ do not commute. Since
$\length\alpha=\ell_P<\ff_M(P)=\length\beta$, the point $P$ is
$\lambda$-doubly thin for every $\lambda > \ff_M(P)$.  Hence 
$\ff_M(P) \ge \nameit$.

Now assume that $\ff_M(P) > \nameit$.  Then since
$\nameit=\Nameit_M(P)$, the point $P$ is $\ff_M(P)$-doubly thin. Hence
there are loops $\beta$ and $\beta'$ based at $P$ having length
$< \ff_M(P) $ and representing non-commuting elements of
$\pi_1(M,P)$. In particular $[\beta]$ and $[\beta']$ cannot both lie
in $C_P$. After re-labeling if necessary we may assume that
$[\beta]\notin C_P$, and hence that $\length\beta\ge\ff_M(P)$.
This is a contradiction, which shows that $\nameit = \ff_M(P)$.
\EndProof

\section{Structure of 4-free groups}\label{kentish section}

The results of this section are analogous to those of \cite[Section
4]{log5}, but are deeper because they require the recent
group-theoretical results established in \cite{kent} and \cite{l-mcr}.

\Number\label{spanier stuff} We will follow the conventions of
\cite{spanier} regarding simplicial complexes and their associated
topological spaces.  By a {\it simplicial complex} we shall mean a set
$V$, whose elements are called {\it vertices}, together with a
collection $S$ of non-empty finite subsets of $V$, called {\it
  simplices}, such that every singleton subset of $V$ is a simplex and
every non-empty subset of a simplex is a simplex.  For $q\ge 0$, a
simplex consisting of $q+1$ vertices will be called a $q$-simplex and
will be said to have dimension $q$. The {\it dimension} of a
simplicial complex is defined to be the supremum of the dimensions of
its simplices.

We emphasize that we do not assume simplicial complexes to be
locally finite, and indeed  in the main
application, in Section \ref{Key section}, the complexes that arise
are locally infinite.

The $q$-dimensional skeleton of a simplicial complex $K$ is the
subcomplex $K^q$ consisting of all $p$-simplices of $K$ for $p\le q$.
(In particular, $K^{-1}$ is the empty complex.) 

The {\it space} $|K|$ of a simplicial complex $K$ is given the
coherent (i.e. weak) topology (see \cite[3.1.14]{spanier}).  If $L$ is
a subcomplex of $K$ then $|L|$ is naturally identified with a subspace
of $|K|$.  For each $q$-simplex $s$ of $K$, the space $|s|$, as a
subspace of $|K|$, is called a {\it closed simplex} and is
homeomorphic to the standard simplex in $\RR^{q+1}$.  In the coherent
topology a subset of $|K|$ is closed (or open) if and only if its
intersection with $|s|$ is closed (or open) for every simplex $s$ of
$K$.

If $s$ is a simplex in a simplicial complex $K$ then the set of all
proper faces of $s$ is a subcomplex denoted by $\dot s$.  For each
simplex $s$ of $K$, the subspace $|s| - |\dot s|$ of $|K|$ will be
called an {\it open simplex} in $|K|$.  If $t$ is a face of $s$ then
the open simplex $|t|-|\dot t|$ will be said to be a face of the open
simplex $|s| - |\dot s|$.

We shall say that a subset $X$ of $|K|$ is {\it saturated} if $X$ is a
union of open simplices. A saturated subset will always be understood
to be endowed with the subspace topology as a subset of $|K|$. Note
that if $X$ is saturated in $|K|$ then the connected components of $X$
are also saturated.

If $v$ is a vertex of a simplicial complex $K$, we will use the
notation $\st_K(v)$ to denote the open star of $v$ in $|K|$ (see
\cite[3.1.22]{spanier}).  The link of $v$ is a subcomplex $L_v$ of $K$,
and we will use the notation $\lk_K(v)$ to denote the closed set
$|L_v|\subset |K|$.

The first barycentric subdivision of a simplicial complex $K$ will
be denoted $K'$, and we shall identify the spaces $|K|$ and $|K'|$.

We define a {\it graph} to be a simplicial complex of dimension at
most $1$. (Thus a graph has no loops or multiple edges.)
\EndNumber

\Number\label{labeled def}
We recall some definitions from \cite{log5}.
 
A group $\Gamma$ will be said to have {\it local rank $\leq k$},
where $k$ is a positive integer, if every finitely generated
subgroup of $\Gamma$ is contained in a subgroup of rank $\leq k$.
The local rank is the smallest integer $k$ with this property, and
is defined to be $\infty$ if no such integer exists.
Note that for a finitely generated group, the local rank is equal
to the rank.
 
Let $\Gamma$ be a group. By a {\it $\Gamma$-labeled complex} we shall
mean an ordered pair $(K,(C_v)_v)$, where $K$ is a simplicial complex
and $(C_v)_v$ is a family of infinite cyclic subgroups of $\Gamma$
indexed by the vertices of $K$. If $(K,(C_v)_v)$ is a $\Gamma$-labeled
complex then for any saturated subset $W$ of $|K|$ we shall
denote by $\Theta(W)$ the subgroup of $\Gamma$ generated by all the
groups $C_v$, where $v$ ranges over the vertices of open simplices
contained in $W$.
\EndNumber

The following simple group-theoretic lemma is needed for the next
two results.  The notion of a $k$-free group was defined in
\ref{k-free def}.

\Lemma\label{free product} Suppose that $\Gamma$ is a $k$-free group
for a given integer $k>0$, that $R\le\Gamma$ is a subgroup of rank
$<k$, that $C$ is a cyclic subgroup of $\Gamma$, and that the rank of
$\langle R\cup C\rangle$ is strictly greater than that of $R$. Then $C$ is
infinite cyclic, and $\langle R\cup C\rangle$ is the free product of the
subgroups $R$ and $C$.  \EndLemma

\Proof Let $r<k$ denote the rank of $R$. Since $\Gamma$ is $k$-free,
$R$ is a free group and has a basis $\{x_1,\ldots,x_r\}$. If $t$
denotes a generator of $C$ then $\langle R\cup C\rangle$ is generated by
$x_1,\ldots,x_r,t$. In particular its rank is at most $r+1$. Since by
hypothesis $\langle R\cup C\rangle$ has rank $>r$, its rank is exactly
$r+1$. Since $r+1\le k$, the group $\langle R\cup C\rangle$ is free, and
hence the set $\{x_1,\ldots,x_r,t\}$ is a basis (see \ref{kurosh
  stuff}). The conclusions follow.
\EndProof

\Proposition\label{fish or fowl} Let $k$ and $r$ be integers with $k>
r\ge2$ and $k\ge4$, Let $\Gamma$ be a $k$-free group, and let
$(K,(C_v)_v)$ be a $\Gamma$-labeled complex. Let $W$ be a connected,
saturated subset of $|K|$ such that $\Theta(\sigma)$ has rank exactly
$r$ for every open simplex $\sigma\subset W$.  Suppose in addition
that {\it either}
\begin{enumerate}[(i)]
\item\label{fish}
there is an integer $n\ge2$ such that
every open simplex contained in $W$ has dimension $n$ or $n-1$,  {\it or}
\item\label{fowl}
$r=2$.
\end{enumerate}
Then $\Theta(W)$ has local rank 
at most $r$.
\EndProposition

\Proof We must show that if $E$ is a finitely generated subgroup 
of $\Theta(W)$, then $E$ is contained in some finitely generated
subgroup of $\Theta(W)$ whose rank is at most $r$. It follows from the
definition of $\Theta(W)$ that $E\le\Theta(W_0)$ for some saturated
subset $W_0\subset W$ which is a union of finitely many open
simplices. Since $W$ is connected, there is a connected (non-empty)
saturated set $V\supset W_0$ which is also a union of finitely many
open simplices. We shall prove the proposition by showing that under
either of the hypotheses (\ref{fish}) or (\ref{fowl}), the finitely
generated group $\Theta(V)\ge E$ has rank at most $r$.

By connectedness we may list the open simplices contained in $V$ as
$\sigma_0,\ldots,\sigma_m$, where $m\ge0$, and for each index $ j $ with
$0< j \le m$ there is an index $l$ with $0\le l< j $ such that either
$\sigma_l$ is a proper face of $\sigma_ j $, or $\sigma_ j $ is a proper
face of $\sigma_l$. We set  $V_ j =\sigma_0\cup\cdots\cup\sigma_ j $ 
for each $ j $ with $0\le  j \le m$. We shall show by induction on $ j $,
for $ j =0,\ldots,m$, that $\Theta(V_ j )$ has rank at most $r$.

By hypothesis, $\Theta(\sigma_ j )$ has rank $r$ for $0\le  j \le m$. In
particular, $\Theta(V_0)=\Theta(\sigma_0)$ has rank $r$; this is the
base case of the induction. Now suppose that $0< j \le m$ and that
$\Theta(V_{ j -1})$ has rank at most $r$. According to the rule for ordering the
open simplices contained in $V$, we may fix an index $l$ with $0\le
l< j $ such that either $\sigma_l$ is a proper face of $\sigma_ j $, or
$\sigma_ j $ is a proper face of $\sigma_l$. If $\sigma_ j $ is a face of
$\sigma_l$, then $\Theta(V_ j )=\Theta(V_{ j -1})$, and the induction step
is trivial. 

Now suppose that $\sigma_l$ is a proper face of $\sigma_ j $. Set
$P=\Theta(V_{ j -1})$, $Q=\Theta(\sigma_ j )$, and
$R=\Theta(\sigma_l)$.  Then $R\le P\cap Q$ and $\Theta(V_ j )=\langle
P\cup Q\rangle$. Then $P$ has rank at most $r$, while $Q$ and $R$ have
rank exactly $r$. We are required to show that $\langle P\cup
Q\rangle$ has rank at most $r$.

First consider the case in which hypothesis (\ref{fish}) holds. Then
since $\sigma_l$ is a proper face of $\sigma_ j $, we must have
$\dim\sigma_ j =n$ and $\dim\sigma_l=n-1$. Let $v$ denote the vertex of
$\sigma_ j $ which is not a vertex of $\sigma_l$, and set $C=C_v$. Then
$C$ is an infinite cyclic group, and $Q=\langle R\cup C\rangle$. Hence
$\langle P\cup Q\rangle=\langle P\cup C\rangle$.

We need to show that $\langle P\cup C\rangle$ has rank at most
$r$. Assume to the contrary that its rank is $>r$. Since $\Gamma$ is
$k$-free, with $k>r$, and $P$ has rank at most $r$, it then follows
from Lemma \ref{free product} that $\langle P\cup C\rangle$ is the
free product of the subgroups $P$ and $C$. In particular, since $R\le
P$, it follows that $Q=\langle R\cup C\rangle$ is the free product of the
subgroups $R$ and $C$. But this is impossible since $Q$ and $R$ have
rank exactly $r$ (and are free since $r<k$), whereas $C$ is infinite
cyclic. This completes the induction step in this case.

We now turn to the case in which hypothesis (\ref{fowl}) holds. In
this case, $Q$ and $R$ have rank $2$, while $P$ has rank at most $2$
and contains $R$. Since a group of rank $\le1$ cannot have a subgroup
of rank $2$, the rank of $P$ must be equal to $2$ as well. Since $P$
and $Q$ are of rank $2$, the group $\langle P\cup Q\rangle$ certainly
has rank at most $4$, and since $4\le k$ it follows that $\langle
P\cup Q\rangle$ is free. Since $P$ and $Q$ are rank-$2$ subgroups of a
free group, it follows from the main theorem of \cite{burns} that
$P\cap Q$ has rank at most $2$. But $P\cap Q$ contains the rank-$2$
subgroup $R$, and hence cannot have rank $\le1$. Thus $P\cap Q$ has
rank exactly $2$.

We now appeal to the main result of \cite{kent} and \cite{l-mcr}, which
asserts that if $P$ and $Q$ are rank-$2$ subgroups of a free group and
$P\cap Q$ has rank $2$, then $\langle P\cup Q\rangle$ also has rank
$2$. This completes the induction in this case.
\EndProof
 
\Proposition\label{normalizer is good}
Let $k> r\ge1$ be integers. Let $\Delta$ be a $k$-free group, and
suppose that $\Delta$ has a normal subgroup of local rank $r$. Then
$\Delta$ has local rank at most $r$.
\EndProposition

\Proof We must show that if $E$ is a finitely generated subgroup of
$\Delta$, then $E$ is contained in some finitely generated subgroup of
$\Delta$ whose rank is at most $r$. Let $\{x_1,\ldots,x_m\}$ be a
finite generating set for $E$. Let $N$ be a normal subgroup of
$\Delta$ whose local rank is $r$. Since $r\ge1$, we may select a
non-trivial element $t$ of $N$. The finitely generated subgroup
$\langle t,x_1tx_1^{-1},\ldots,x_mtx_m^{-1}\rangle$ of $N$  is contained
in some subgroup $H_0$ of $N$ whose rank is at most $r$. For
$j=1,\ldots,m$ we denote by $H_j$ the subgroup $\langle
H_0\cup\{x_1,\ldots,x_j\}\rangle$  of $\Delta$. We shall show by
induction on $j=0,\ldots,m$ that $H_j$ has rank at most $r$. Since
$E\le H_m$, this implies the conclusion.

The base case is clear because $H_0$ was chosen to have rank $\le r$.
Assume for a given $j$, with $0\le j<m$, that $H_j$ has rank $\le
r$. If the rank of $H_{j+1}=\langle H_j\cup \{x_{j+1}\}\rangle$ is
$>r$, then it follows from Lemma \ref{free product} that $H_{j+1}$ is
the free product of the subgroups $H_j$ and $\langle x_{j+1}\rangle$,
and that $x_{j+1}$ has infinite order. But $H_0\le H_j$ contains the
elements $t$ and $t'=x_{j+1}tx_{j+1}^{-1}$.   It follows that $t'$
is an element of a factor in the free product $H_j\star\langle
x_{j+1}\rangle$ which can also be written in a normal form of length
$3$ in the free product. This contradicts the uniqueness of normal
form.  Hence $H_{j+1}$ must have
rank at most $r$.  \EndProof

\section{Existence of \texorpdfstring{$\log7$}{log 7}-\semithick\ points}\label{Key section}

The goal of this section is to prove Theorem \ref{Real Key theorem},
which was stated in the introduction.

\Number\label{nerve}
In this section we will use methods involving nerves of coverings that
were first introduced in \cite{log5}. Recall that if
$\calu=(U_i)_{i\in I}$ is an indexed covering of a topological space
by non-empty open subsets, the {\it nerve} of $\calu$ is an abstract
simplicial complex whose vertices are in bijective correspondence with
the elements of the index set $I$. If we denote by $v_i$ the vertex
corresponding to the index $i$, then a $k$-simplex of the nerve is a
collection $\{v_{i_0},\dots,v_{i_k}\}$ of $k+1$ vertices such that
$U_{i_0}\cap\dots\cap U_{i_k}\neq\emptyset$.
\EndNumber

\Number\label{john, you have one big nerve}
  Let $\Gamma$ be a group and let $(K,(C_v)_v)$ be a $\Gamma$\hyph
labeled complex.  By a {\it\lc\ action of $\Gamma$ on} $(K,(C_v)_v)$ we
shall mean a simplicial action of $\Gamma$ on $K$ such that for each
vertex $v$ of $K$ we have $C_{\gamma\cdot v}=\gamma C_v\gamma^{-1}$.
(This concept was first introduced explicitly in \cite{accs}, where the
term ``natural action'' was used. We are afraid this term may be
confusing.)
\EndNumber

\Number\label{oy, it's so natural!}
It follows from the definition of a \lc\ action of a group $\Gamma$
on a $\Gamma$-labeled complex $K$ that for every
saturated subset $A$ of $|K|$, and every element $\gamma$ of $\Gamma$, we
have $\Theta(\gamma\cdot A)=\gamma \Theta(A)\gamma^{-1}$.

In particular, if a saturated subset $A$ of $|K|$ is invariant under an
element $\gamma$ of $\Gamma$, then  $\Theta(A)$ is normalized by
$\gamma$. In other words, the stabilizer in $\Gamma$ of any saturated
subset $A$ of $|K|$ is contained in the normalizer of $\Theta(A)$.
\EndNumber

\Number\label{no sich a person} Let $\Gamma$ be a discrete, purely
loxodromic subgroup of $\Isom_+(\HH^3)$, and let $\lambda>0$ be a
number such that the indexed family of cylinders
$\calz=(Z_{\lambda}(C))_{C\in\calc_\lambda(\Gamma)}$ covers $\HH^3$.
Let $K$ denote the nerve of the covering $\calz$.  According to the
definition of the nerve, the vertices of $K$ are in bijective
correspondence with elements of $\calc_{\lambda}(\Gamma)$.  If $C_v$
denotes the element of $\calc_{\lambda}(\Gamma)$ corresponding to a
vertex $v$, then $C_v$ is by definition a maximal, and hence infinite,
cyclic subgroup of $\Gamma$. Thus $(K,(C_v)_{v})$ is a
$\Gamma$-labeled complex.  \EndNumber

\Proposition\label{natural}Suppose that $\Gamma$ is a discrete,
torsion-free subgroup of $\Isom_+(\HH^3)$, and that $\lambda>0$ is a
number such that the indexed family of sets
$\calz=(Z_{\lambda}(C))_{C\in\calc_{\lambda}}(\Gamma)$ covers $\HH^3$. Let 
$(K,(C_v)_{v})$ be the $\Gamma$-labeled complex defined as in \ref{no
sich a person}. Then $(K,(C_v)_{v})$ admits a \lc\ $\Gamma$-action.
\EndProposition

\Proof For any vertex $v$ of $K$ we have $C_v\in\calc_{\lambda}(\Gamma)$, and
hence $\gamma C_v\gamma^{-1}\in\calc_{\lambda}(\Gamma)$ for every
$\gamma\in\Gamma$. We may therefore define an action of $\Gamma$ on
the set of vertices of $K$ by $C_{\gamma\cdot v}=\gamma
C_v\gamma^{-1}$. If $v_0,\dots,v_m$ are the vertices of an $m$-simplex
of $K$, then for every $\gamma\in\Gamma$ we have
$$\bigcap_{0\leq i\leq m}Z_{\lambda}(\gamma C_{v_i}\gamma^{-1})=
\bigcap_{0\leq i\leq m}\gamma\cdot Z_{\lambda}(C_{v_i})=
\gamma\cdot\bigcap_{0\leq i\leq m} Z_{\lambda}(C_{v_i})
\neq\emptyset,$$
so that $\gamma\cdot v_0,\dots,\gamma\cdot v_m$ are
the vertices of an $m$-simplex of $K$. Thus the action of $\Gamma$ on
the vertex set extends to a simplicial action on $K$. It is immediate
from the definitions that this is a labeling-compatible action on
$(K,(C_v)_v)$. 
\EndProof

\Lemma\label{down to 3}For any simplical
complex  $K$, the inclusion $|K^3|-|K^0| \to |K|-|K^0|$ induces
isomorphisms on $\pi_0$ and $\pi_1$.
\EndLemma

\Proof The space $|K|-|K^0|$ is the direct limit of the system of
subspaces $|L|-|K^0|$, where $L$ ranges over all subcomplexes of $K$
such that (1) $|K^3|\subset |L|$ and (2) $|L|-|K^3|$ contains only
finitely many open simplices. Hence it suffices to show that for any
subcomplex satisfying (1) and (2), the inclusion $|L|-|K^0|\to
|K|-|K^0|$ induces isomorphisms on $\pi_0$ and $\pi_1$. By induction
on the number of simplices in $|L|-|K^3|$ , this reduces to showing
that if $L$ and $L_1 $ are subcomplexes of $K$ satisfying (1) and (2),
and $|L_1 |=|L|\cup\sigma$ for some open simplex $\sigma$, then the
inclusion $|L|-|K^0|\to |L_1 |-|K^0|$ induces isomorphisms on $\pi_0$
and $\pi_1$.  If $d>3$ denotes the dimension
of $\sigma$, then $B=\overline{\sigma}$ is a closed topological
$d$-ball, and the vertex set $V$ of $\sigma$ is a finite subset of the
$(d-1)$-sphere $\partial B$. Since $d-1>2$, the set $(\partial B)-V$
is connected and simply connected. The set $ B-V$ is contractible.
Since $|L_1 |=|L|\cup(B-V)$ and $|L|\cap(B-V)=(\partial B)-V$, the
inclusion $|L|-|K^0|\to |L_1 |-|K^0|$ induces isomorphisms on $\pi_0$
and $\pi_1$ as required.  \EndProof

\Lemma\label{space to complex}
Let $H$ be a topological space which has the homotopy type of a
CW-complex.  Let $\calu=(U_i)_{i\in I}$ be a covering
of $H$ by contractible open sets.  Suppose that
\begin{enumerate}[(I)]
\item\label{lox} for every finite non-empty subset
  $\{{i_0},\dots,{i_k}\}$  of $I$, the set $U_{i_0}\cap\dots\cap
  U_{i_k}\subset H$ is either empty or contractible; and
\item\label{sable} for every point $P\in H$ there exist distinct
  indices $i,j\in I$ such that $P\in U_i\cap U_j$.  
\end{enumerate}

Let $K$ denote the nerve of $\calu$.  Then the space $|K|-|K^0|$ is
homotopy-equivalent to $H$.
\EndLemma

\Proof According to Borsuk's Nerve Theorem (see for example
\cite[Theorem 6 and Remark 7]{bjorner}), Property (\ref{lox}) of the
covering implies that the space $|K|$ is homotopy-equivalent to
$H$. We shall complete the proof by showing that the inclusion
$|K|-|K^0|\to |K|$ is a homotopy equivalence.

If $k\in I$ is any index, let us denote by $J_k$ the set of all
indices $i\in I$ such that $i\ne k$ but $U_i\cap
U_k\ne\emptyset$. Condition (\ref{sable}) implies that the indexed family
$\calv_k=(U_i\cap U_k)_{k\in J_k}$ is an open covering of the set
$U_k$, to which we assign the subspace topology.  

It follows from the definitions that the nerve of $\calv_k$ is
simplicially isomorphic to the link of the vertex $v_k$ in the nerve
$K$ of $\calu$. (Although this observation is purely formal, it
depends on defining the nerve via the index set $J_k$ as in
\ref{nerve}.  Different indices in $J_k$ may define the
same set in $\calv_k$ even if they define different sets in $\calu$.)

Since $\calu$ satisfies condition (\ref{lox}), it is clear that
(\ref{lox}) remains true when we replace $\calu$ by $\calv_k$. Hence
by Borsuk's Nerve Theorem, the nerve of the covering $\calv_k$ is
homotopy-equivalent to $U_k$, i.e. it is contractible.  This shows
that the link in $K$ of every vertex of $K$ is contractible.

To show that the inclusion $|K|-|K^0|\to |K|$ is a homotopy
equivalence, we first note that  since $|K|$ has the weak
topology,  we may regard $|K|$ as the topological direct limit of
the subspaces $X_F=(|K|-|K^0|)\cup F$, where $F$ ranges over the
finite subsets of $|K^0|$. Hence it suffices to prove that
$|K|-|K^0|\to X_F$ is a homotopy equivalence for every finite
$F\subset |K^0|$. By induction on the cardinality of $F$, this reduces
to showing that if $F\subset |K^0|$ is finite, if $v\in |K^0|-F$, and
if we set $F'=F\cup\{v\}$, then the inclusion $X_F\to X_{F'}$ is a
homotopy equivalence. 

We recall  from \ref{spanier stuff} 
that $K'$ denotes the first barycentric subdivision of $K$. For
general reasons, $X_{F'}-\st_{K'}(v)$ is a deformation retract of
$X_F$. On the other hand, $X_{F'}-\st_{K'}(v)$ is a deformation retract
of $X_{F'}$, because $\lk_{K'}(v)$ is homeomorphic to $\lk_K(v)$ and is
therefore contractible.  \EndProof

\Definitions\label{calgary} Let $K$ be a simplicial
complex, and let $A$ and $B$ be disjoint saturated subsets of
$|K|$. We shall say that $A$ and $B$ are {\it adjacent} if there are
open simplices $\sigma\subset A$ and $\tau\subset B$ such that either
$\sigma$ is a face of $\tau$ or $\tau$ is a face of $\sigma$.

If $K$ is a simplicial complex and if $X_0$ and $X_1$ are
disjoint, saturated subsets of $|K|$, we define a simplicial
complex $\calg(X_0,X_1)$ of dimension at most $1$ as
follows.  The vertices of $\calg(X_0,X_1)$ are the elements
of $\calw_0\cup\calw_1$, where $\calw_i$ denotes the set of connected
components of $X_i$. A $1$-simplex is a pair $\{W_0,W_1\}$, where
$W_i\in\calw_i$ for $i=0,1$, and $W_0$ and $W_1$ are adjacent.
For $W\in\calw_0\cup\calw_1$, we denote by $v_W$ the
vertex of $\calg(X_0,X_1)$ corresponding to $W$.
\EndDefinitions

\Definition\label{ILGWU}
Let $\calg$ be a graph. A simplicial action of a group
$\Gamma$ on $\calg$ will be said to have {\it no inversions} if for
each edge $e$ of $\calg$ and for each element $\gamma$ of $\Gamma$
such that $\gamma\cdot e=e$, we have $\gamma\cdot v=v$ for each
endpoint $v$ of $e$.
\EndDefinition

\Number\label{complex to graph nonsense}Suppose that a group $\Gamma$
acts simplicially on a simplicial complex $K$, and
that $X_0$ and $X_1$ are disjoint, saturated subsets of $|K|$, each of
which is invariant under the action of $\Gamma$. Then for $i=0,1$, an
arbitrary element of $\Gamma$ maps each connected component of $X_i$
onto a connected component of $X_i$, and so the action of $\Gamma$
defines an action on the vertex set of $\calg(X_0,X_1)$. Furthermore,
if $W_i$ is a component of $X_i$ for $i=0,1$, and $W_0$ and $W_i$ are
adjacent, then $\gamma\cdot W_0$ and $\gamma\cdot W_1$ are adjacent
for every $\gamma\in\Gamma$, since the action of $\Gamma$ on $K$ is
simplicial. Thus the action of $\Gamma$ on the vertex set of
$\calg(X_0,X_1)$ extends to an action on $\calg(X_0,X_1)$. Note that
this action has the property that the stabilizer in $\Gamma$ of any
vertex of $\calg(X_0,X_1)$ is the stabilizer of some component of
$X_0$ or $X_1$ under the action of $\Gamma$ on $K$. Note also that for
$i=0,1$, the set of vertices of $\calg(X_0,X_1)$ corresponding to
components of $X_i$ is $\Gamma$-invariant. In particular, $\Gamma$
acts without inversions on $\calg(X_0,X_1)$.
\EndNumber

\Definition
Let $X$ and $Y$ be topological spaces.
We define a {\it homotopy-retraction} from $X$ to $Y$ to be a map
$r:X\to Y$ which admits a right homotopy inverse. We shall say that
$Y$ is a {\it homotopy-retract} of $X$ if there exists a homotopy
retraction from $X$ to $Y$.
\EndDefinition

\Lemma\label{complex to tree}
Suppose that $K$ is a simplicial complex and that $X_0$ and
$X_1$ are saturated subsets of $|K|$. Then $|\calg(X_0,X_1)|$ is a
homotopy-retract of the saturated subset $X_0\cup X_1$ of $|K|$.
\EndLemma

\Proof We set $\calg=\calg(X_0,X_1)$, and we use the notation of
Definition \ref{calgary}. For $m=0,1$ we denote by $V_m$ the set of
vertices of $\calg$ corresponding to elements of $\calw_m$.

Recall  from \ref{spanier stuff}
  that the first barycentric subdivision of $K$ is denoted by
$K'$, and that we identify $|K'|$ with $|K|$.  For each open
simplex $\sigma$ in $K$ we shall denote by $b_\sigma$ the barycenter
of $\sigma$, which is a vertex of $ K'$.

For each $m\in\{0,1\}$, we let $B_m$ denote the set of all vertices of
$ K'$ that have the form $b_\sigma$ for some open simplex $\sigma$ of
$K$ such that $\sigma\subset X_m$.    If
$b=b_\sigma\in B_m$, we denote by $W_b$ the component of $X_m$
containing $\sigma$, and we set $v_b=v_{W_b}$. We claim:

\Claim\label{dinah}If $m\in\{0,1\}$ and $b,c\in B_m$ are
given, and if $b$ and $c$ are joined by an edge of $ K'$, 
then $v_b=v_c$.
\EndClaim

To prove this we write $b=b_\sigma$ and $c=b_\tau$, where the open
simplices $\sigma$ and $\tau$ of $K$ are contained in $X_m$. The
hypothesis that $b$ and $c$ are joined by an edge of $ K'$
means that one of the simplices $\sigma$, $\tau$ is a face of the
other. By symmetry we may assume that $\tau$ is a face of $\sigma$.
Then $\tau\subset\overline{\sigma}$, and hence $\sigma$ and $\tau$ are
contained in the same component of $X_m$. In the notation introduced
above this means that $W_b=W_c$, and \ref{dinah} follows.

Next we claim:

\Claim\label{hedgehog}If $b_0\in B_0$ and $b_1\in B_1$ are
 joined by an edge of $ K'$, then 
 $v_{b_0}$ and $v_{b_1}$ are joined by an edge of $\calg$.
\EndClaim

To prove this we write $b_m=b_{\sigma_m}$ for $m=0,1$, where the open
simplex $\sigma_m$ of $K$ is contained in $X_m$.  Here again we may
assume by symmetry that $\sigma_0$ is a face of $\sigma_1$. Then
$\sigma_0\subset\overline{\sigma_1}$, and hence $W_{b_0}\subset
\overline W_{b_1}$. By the definition of the graph
$\calg=\calg(X_0,X_1)$, the conclusion of \ref{hedgehog} follows.

Now let $S$ denote the set of all simplices $s$ of $ K'$ such that
$s\subset X_0\cup X_1$. For each $s\in S$, we let $\calb_s$ denote the
set of all vertices of $s$ that lie in $B_0\cup B_1$. Note that if
$\sigma$ denotes the open simplex of $K$ containing $s$, then
$\sigma\subset X_0\cup X_1$ and hence $b_\sigma\in\calb_s$; in
particular $\calb_s\ne\emptyset$ for any $s\in S$.   For  each 
$s\in S$ we
set 
$$\calv_s=\{v_b:b\in\calb_s\}.$$
It follows from \ref{dinah} that $\calv_s$ contains at most one vertex
from each of the sets $B_0$ and $B_1$, and it follows from
and \ref{hedgehog} that if $\calv$ contains a vertex in $B_0$ 
and a vertex in $B_1$, then these vertices are joined by an
edge. Hence $\calv_s$ is the vertex set of a simplex of dimension $0$
or $1$ in $\calg$. We denote the corresponding open simplex of
$|\calg|$ by $\alpha_s$.  Thus $\alpha_s$ is either a vertex or an
open edge of $|\calg|$.

It is immediate from the definition of $\calv_s$ that if $t,s\in S$
and if $t$ is a face of $s$, then $\calv_t\subset\calv_s$. Hence:

\Claim\label{flamingo} If $t,s\in S$ and if $t$ is a face of $s$, then
$\alpha_t$ is a face of $\alpha_s$ in the simplicial complex
$\calg$. (In other words, either $\alpha_t=\alpha_s$, or $\alpha_s$ is
an edge and $\alpha_t$ is one of its endpoints.)
\EndClaim

Using \ref{flamingo} we may construct, recursively for $k\ge-1$, a
continuous map $r_k:(K')^k\cap(X_0\cup X_1)\to|\calg|$  such that
$r(s)$ is contained in (the open simplex) $\alpha_s$ for every simplex
$s\in S$ of dimension at most $k$. We take $r_{-1}$ to be the empty
map, and construct $r_{k+1}$ as an extension of $r_k$. If $r_k$ has
been defined, and $s\in S$ is $(k+1)$-dimensional, then every face $t$
of $s$ which belongs to $S$ is mapped by $r_k$ into $\alpha_t$,
which by \ref{flamingo} is a face of $\alpha_s$. This allows us to
extend $r_k|(\partial s)\cap(X_0\cup X_1)$ to a continuous map of
$\overline{s}\cap(X_0\cup X_1)$ into $\overline{\alpha_s}$ which maps
$s$ into $\alpha_s$, and thus to complete the recursive
definition. Since $r_{k+1}$ is an extension of $r_k$, we may define a
map $r:X_0\cup X_1\to|\calg|$ by setting $r|(K')^k\cap(X_0\cup
X_1)=r_k$  for each $k$. The construction gives:
\Claim\label{dormouse} Each simplex $s\in S$ is mapped by $r$ into
$\alpha_s$.
\EndClaim

We now claim:
\Claim\label{more mouse}
For each $m\in\{0,1\}$ and each component $W$ of $X_m$, the set $r(W)$
is contained in $\st_{\calg}(v_W)$.
\EndClaim

To prove this, we consider any open simplex $s$ of $K'$ which is
contained in $W$. Since $W$ is saturated in $K$, the open simplex
$\sigma$ of $K$ containing $s$ is also contained in $W$. We have
$v_W=b_\sigma\in\calb_s$, so that $\alpha_s$ is contained in
$\st_\calg(v_W)$. By (\ref{dormouse}) it follows that $r(s)\subset
\st_{\calg}(v_W)$, and \ref{more mouse} follows.

Now we construct a map $i:|\calg|\to X_0\cup X_1$ as follows. For each
vertex $v=v_W$ of $\calg$ we take $i(|v|)$ to be a point of $W$ (chosen
arbitrarily). If $e$ is an edge of $\calg$ with endpoints $v_0$ and
$v_1$, where $v_m=v_{W_m}\in V_m$ for $m=0,1$, then the components
$W_0$ and $W_1$ of $X_0$ and $X_1$ are adjacent. In particular, one of
them meets the closure of the other, and hence $W_0\cup W_1$ is
connected. Hence $i(v_0)$ and $i(v_1)$ are joined by a path in
$W_0\cup W_1$, and we may use this path to extend $i$ to $|e|$. The
construction of the map $i$ gives:
\Claim\label{white knight}
For each vertex $v=v_W$ of $\calg$ we have $i(|v|)\subset W$, and for
each edge $e$ of $\calg$ with endpoints $v_{W_0}$ and $v_{W_1}$ we
have $i(|e|)\subset W_0\cup W_1$.
\EndClaim

From \ref{more mouse} and \ref{white knight} we immediately deduce the
following property of the composition $r\circ i:|\calg|\to|\calg|$.

\Claim\label{humpty-dumpty}
Each vertex $v$ of $|\calg|$ is mapped by $r\circ i$ into
$\st_{\calg}(v)$ and each edge of $|\calg|$ is mapped into the union
of the open stars of its endpoints.
\EndClaim

From \ref{humpty-dumpty} we can deduce that $r\circ i$ is homotopic to
the identity. First we construct a homotopy $H^{(0)}$ from the vertex
set $|\calg^{(0)}|$ into $|\calg|$ such $H^{(0)}_0=r\circ i$,
$H^{(0)}_1$ is the inclusion, and $H^{(0)}(\{|v|\}\times[0,1])$ is
contained in  $\st_{\calg}(v)$  for every vertex $v$; this is
possible by \ref{humpty-dumpty} and the connectedness of the open
stars. Then we consider an arbitrary edge $|e|$ of $|\calg|$, and let $Y$
denote the union of the open stars of its endpoints. According to
\ref{humpty-dumpty} and the construction of $H^{(0)}$, we have
$H^{(0)}((|\dot e|)\times[0,1])\cup r\circ i(|e|)\subset
W$. Since $W$ is simply connected, we may extend $H^{(0)}|(\partial
e)\times[0,1]$ to a homotopy $H^e:|e|\times[0,1]\to W$ such that
$H^e_1$ is the inclusion map from $|e|$ to $W$.

This shows that $r$ is a homotopy retraction from $X_0\cup X_1$ to $\calg$.
\EndProof

\Lemma\label{no such action} Let $M$ be a closed, orientable, aspherical
$3$-manifold. Then $\pi_1(M)$ does not admit a simplicial action
without inversions on a tree $T$ with the property that the stabilizer
in $\pi_1(M)$ of every vertex of $T$ is a locally free subgroup of
$\pi_1(M)$.  \EndLemma

\Proof We will prove this from the point of view used in
\cite{splittings}. Assume an action with the stated properties
exists. It follows from the Bass-Serre theory (see \cite[Theorem
7]{Tret}) that $\pi_1(M)$ is then isomorphic to a graph of groups in
which each vertex group or edge group is respectively isomorphic to
the stabilizer of a vertex or edge of $T$. If the graph of groups is
non-trivial in the sense that the full group $\pi_1(M)$ is not a
vertex group, it follows from \cite[Proposition 2.3.1]{splittings}
that there is an incompressible surface $F\subset M$ such that the
image (defined up to conjugacy) of the inclusion $\pi_1(F)\to\pi_1(M)$ is
contained in an edge group. Since $M$ is simple and closed, this image
is the fundamental group of a surface of genus $>1$. This is
impossible because any edge group is contained in a vertex group, and
is therefore locally free by our assumption. If the graph of groups is
trivial, then $\pi_1(M)$ is itself locally free, and since $M$ is
compact, $\pi_1(M)$ is in fact free. But since $M$ is closed,
orientable and aspherical, we have $H_3(\pi_1(M);\ZZ)\cong\ZZ$,
whereas a free group has trivial homology in dimensions $>1$. Thus
this case is impossible as well.  \EndProof

\Proof[Proof of Theorem \ref{Real Key theorem}]
Assume that the conclusion is false, so that every point of $M$ is
($\log 7$)-doubly thin. Then it follows from Proposition \ref{hootchy-kootcher}
that
$$\HH^3=\bigcup_{\begin{matrix}
    C,C'\in\calc_{\log7}(\Gamma)\cr C\ne C'\end{matrix}}
    Z_{\log7}(C)\cap Z_{\log7}(C').$$
Thus the indexed family of sets
$\calz=(Z_{\log7}(C))_{C\in\calc_{\log7}(\Gamma)}$ is an open covering of
$\HH^3$, and Condition (\ref{sable}) of Lemma \ref{space to complex}
holds with $H=\HH^3$ and $\calu=\calz$. Condition (\ref{lox}) of Lemma
\ref{space to complex} follows from the convexity of the cylinders
$Z_{\log7}(C)$. Hence, according to Lemma \ref{space to complex}, if
we let $K$ denote the nerve of $(U_i)$, the
space $|K|-|K^0|$ is homotopy-equivalent to $\HH^3$ and is therefore
contractible. It then follows by Lemma \ref{down to 3} that the set
$|K^3|-|K^0|$ is connected and simply connected.

According to the observation in \ref{no sich a person}, with
$\lambda=\log7$, we may define a $\Gamma$-labeled complex
$(K,(C_v)_{v})$ by taking $C_v$ to be the element of 
$\calc_{\lambda}(\Gamma)$ corresponding to the vertex $v$.  

Let $\sigma$ be any open simplex contained in $|K^3|-|K^0|$, and
consider the group $\Theta(\sigma)$, in the notation of \ref{labeled
  def}. Since $\sigma$ has dimension at most $3$, the rank of
$\Theta(\sigma)$ is at most $4$. 

Suppose that $\Theta(\sigma)$ has rank $4$. Then $\sigma$ must be
$3$-dimensional. Let $v_0,\ldots,v_3$ denote the vertices of $\sigma$,
and for $i=0,\ldots,3$ set $C_i=C_{v_i}$ and
$Z_i=Z_{\log7}(C_i)$. Choose a generator $x_i$ of each
$C_i$.  The definition of the nerve
implies that the intersection $Z_0\cap Z_1\cap Z_2\cap Z_3$ is
non-empty. We choose a point $z$ in this intersection. 

For $i=0,\ldots,3$,
since $z\in Z_i=Z_{\log7}(C_i)$,  there exists an integer $m_i\ne0$ such that 
\Equation\label{log 7 contradiction}
\dist(z,x_i^{m_i}\cdot z)<\log7 \delete{\qquad{\rm for\ }i=0,\ldots,{3}}.
\EndEquation

Since
$\Theta(\sigma)$ is assumed to have rank $4$, and since
$\Gamma\cong\pi_1(M)$ is $4$-free by hypothesis, $\Theta(\sigma)$ is
freely generated by $x_0,\ldots,x_3$. Hence  $x_0^{m_0},\ldots,x_3^{m_3}$ also freely
generate a rank-$4$ free group. But it follows from 
 (\cite[Theorem 6.1]{accs},
 together with the main result of \cite{agol} or \cite{cg}), that if
$\xi_1,\ldots,\xi_k\in\isomplus(\HH^3)$ freely generate a discrete
 subgroup of $\isomplus(\HH^3)$, then for any $z\in\HH^3$ we have 
$$\max_{1\le i \le k}\dist(z,\xi_i\cdot z)\ge\log(2k-1).$$
Taking $k=4$ and $\xi_i=x_{i+1}^{m_{i+1}}$, we obtain a contradiction. This
shows that $\Theta(\sigma)$ has rank at most $3$ for any
open simplex $\sigma\subset |K^3|-|K^0|$.

On the other hand, for any
open simplex $\sigma\subset |K^3|-|K^0|$, there are at least two
distinct vertices, say $v_C$ and $v_{C'}$, in $\overline{\sigma}$. Here $C$ and
$C'$ are by definition distinct maximal cyclic subgroups of $\Gamma$,
and hence $\Theta(\sigma)$, which contains $C$ and $C'$, is
non-abelian. This shows that $\Theta(\sigma)$ must have rank at least $2$.

Thus for any open simplex $\sigma\subset |K^3|-|K^0|$, the group
$\Theta(\sigma)$ has rank $2$ or $3$. We may therefore write 
$|K^3|-|K^0|$ as a disjoint union
$$|K^3|-|K^0|=X_2\dot\cup X_3,$$
where $X_k$ is the union of all open simplices $\sigma \subset |K^3|$
for which $\Theta(\sigma)$ has rank $k$.

We claim: 
\Claim\label{prove it ain't so}For any $m\in\{2,3\}$ and for any
component $W$ of $X_m$, the local rank of $\Theta(W)$ is at most $m$.
\EndClaim

This is an application of Proposition \ref{fish or fowl}. In the case
$m=2$, hypothesis (\ref{fowl}) of Proposition \ref{fish or fowl}
clearly holds, and hence $\Theta(W)$ has local rank at most
$2$. If $m=3$ then $\Theta(\sigma)$ has rank $3$ for every open
simplex $\sigma\subset W$. If $d=\dim\sigma$ then $\Theta(\sigma)$ is
generated by $d+1$ cyclic groups, and hence $d\ge2$. But $d\le3$ since
$\sigma\subset X_3\subset |K^3|-|K^0|$. Hence in this case,
hypothesis (\ref{fish}) of Proposition \ref{fish or fowl} holds with
$n=3$, and hence $\Theta(W)$ has local rank at most $3$.

We next claim:

\Claim\label{sez you}
If $W$ is a component of $X_2$ or $X_3$, the normalizer of $\Theta(W)$
in $\Gamma$ has local rank at most $3$.
\EndClaim

To prove \ref{sez you}, we let $r$ denote the local rank of
$\Theta(W)$. By \ref{prove it ain't so} we have $r\le3$. If $r\le1$
then $\Theta(W)$ is locally cyclic, and is therefore cyclic since
$\Gamma$ is a discrete subgroup of $\Isom(\HH^3)$. But this is
impossible because $\Theta(\sigma)\le\Theta(W)$ has rank $2$ or $3$
for every open simplex $\sigma\subset W$. It follows that $r$ is $2$
or $3$.

Now if $\Delta$ denotes the normalizer of $\Theta(W)$ in $\Gamma$,
then $\Delta$ is $4$-free and has the normal subgroup $\Theta(W)$ of
local rank $r$. Since $r$ is $2$ or $3$, it follows from Proposition
\ref{normalizer is good} that $\Delta$ has local rank at most $r$.
This establishes \ref{sez you}.

Let us now set $T=\calg(X_2,X_3)$ (with the notation of
\ref{calgary}). According to Lemma \ref{complex to tree}, $T$ is a
homotopy-retract of $X_2\cup X_3=|K^3|-|K^0|$. As we have seen
that $|K^3|-|K^0|$ is connected and simply connected, $T$ is a
tree.

According to Proposition \ref{natural}, $\Gamma$ admits a \lc\ action
on $K$. It follows from \ref{oy, it's so
natural!} that for any $\gamma\in\Gamma$ and for any simplex $\sigma$
of $K$, the groups $\Theta(\sigma)$ and $\Theta(\gamma\cdot\sigma)$
are conjugate in $\Gamma$ and hence have the same rank. It follows
that each of the saturated sets $X_2$ and $X_3$ is invariant under the
action of $\Gamma$.

Hence, according to \ref{complex to graph nonsense}, there is an
induced action, without inversions, of $\Gamma$ on $T$; and under this
induced action, if $v$ is any vertex of $T$, the stabilizer $\Gamma_v$
of $v$ in $\Gamma$ is the stabilizer of some component $W$ of $X_2$ or
$X_3$ under the action of $\Gamma$ on $|K|$.  Hence by \ref{oy, it's so
natural!}, $\Gamma_v$ is contained in the normalizer of
$\Theta(W)$. By \ref{prove it ain't so}, this normalizer has local
rank at most $3$. In particular, $\Gamma_v$ has local rank at most
$3$, and since $\Gamma$ is in particular $3$-free it follows that
$\Gamma_v$ is locally free.

Thus we have constructed a simplicial action, without inversions, of
$\Gamma\cong\pi_1(M)$ on the tree $T$ with the property that the
stabilizer in $\Gamma$ of every vertex of $T$ is a locally free
subgroup of $\Gamma$. This contradicts Lemma \ref{no such action}.
\EndProof

Theorem \ref{Real Key theorem} will be used in the sequel via the
following corollary.

\Corollary\label{Key theorem}
Suppose that $M$ is a closed, orientable hyperbolic $3$-manifold such
that $\pi_1(M)$ is $4$-free. Then either \Alternatives
\item\label{whenever richard corey went downtown}$M$ contains an
  embedded hyperbolic ball of radius $(\log7)/2$, or
\item\label{the people on the pavement looked at him}there is a point $P\in\XX_M$ with $\ff_M(P)=\log7$ (see
  \ref{eggs and fries}).
\EndAlternatives
\EndCorollary

\Proof
According to Theorem \ref{Real Key theorem} there is a
$\log7$-semithick point $P_0\in M$. By Proposition \ref{the buck stops
  here} and the definition of $\Nameit_M$ in \ref{you asked for it},
we have $\Nameit_M(P_0)\ge\log7$. If it happens that the function
$\Nameit_M(P)$ takes values $\ge\log7$ everywhere in $ M$, then by
Proposition \ref{ti bouffes ti bouffes pas}, the number $\log7$ is a
Margulis number for $M$. According to the discussion in \ref{rev}, $M$
then contains a $\log7$-thick point, and hence Conclusion
(\ref{whenever richard corey went downtown}) of the corollary holds. Now
suppose that $\Nameit_M(P)$ takes a value $<\log7$ somewhere in $M$.
Since $\Nameit_M$ is continuous according to Proposition \ref{you
  can't prove it ain't so}, it follows that there is a point $P\in M$
such that $\Nameit_M(P)=\log7$. According to Proposition \ref{you can
  always telephone girl}, either $P$ is a $\log 7$-thick point of
$M$, in which case Conclusion (\ref{whenever richard corey went
  downtown}) of the corollary holds; or $P\in\XX_M$ and
$\ff_M(P)=\log 7$, which gives Conclusion (\ref{the people on the
  pavement looked at him}) of the corollary.
\EndProof

\section{Caps}\label{capsection}

The term ``cap'' refers to the intersection of a closed ball with a
half-space whose interior does not contain the center of the ball.
We begin by introducing some notation for describing caps.

\Number\label{basic caps} 
Let $z_0$ be a point of $\HH^3$. For each positive real number $R$ we
shall denote by $S(R,z_0)$ the sphere of radius $R$ centered at
$z_0$. Thus $S(R,z_0)$ is the boundary of the closed ball
$\overline{N}(z_0,R) \dot= \overline{N(z_0,R)}$ of
radius $R$ about $z_0$.  For each point $\zeta\in S(R,z_0)$ we shall
denote by $\eta_\zeta$ the ray originating at $z_0$ and passing through
$\zeta$.  We shall endow $S(R,z_0)$ with the spherical metric in which the
distance between two points $\zeta,\zeta'\in S(R,z_0)$ is the angle between
$\eta_\zeta$ and $\eta_{\zeta'}$.

For each point $\zeta\in S(R,z_0)$ and each number $\D\ge0$ 
we shall denote by $\Pi(z_0,\zeta,\D)$ the plane which meets
$\eta_\zeta$ perpendicularly at a distance $\D$ from $z_0$, and by
$H(z_0,\zeta,\D)$ the closed half-space which is bounded by
$\Pi(z_0,\zeta,\D)$ and has unbounded intersection with
  $\eta_\zeta$.  We set
$$K(R,z_0,\zeta,\D)=\overline{N}(z_0,R)\cap H(z_0,\zeta,\D).$$
Thus $K(R,z_0,\zeta,\D)$ is a cap in the closed ball of radius $R$ cut out
by a plane at distance $\D$ from the center. 
Note that $K(R,z_0,\zeta,\D)=\emptyset$ when $\D\ge R$.
\EndNumber

Let $B(r)$ denote the volume of a ball of radius $r$ in $\HH^3$.
The following result was discussed in the Introduction.

\Proposition\label{c.d. rivington} Let $M$ be a closed, orientable
hyperbolic $3$-manifold, write $M=\HH^3/\Gamma$ where
$\Gamma\le\isomplus(\HH^3)$ is discrete and torsion-free, and set
$q=q_\Gamma$. Let $\lambda$ be a positive number, and suppose that $P$
is a point of $\XX_M$ with $\ff_M(P)\ge\lambda$. Let $\tP$ be a point
of $q^{-1}(P)$, let $j:\pi_1(M,P)\to\Gamma$ denote the isomorphism
determined by the base point $\tP\in\HH^3$ (see \ref{same to me}), and
let $x$ denote a generator of $j(C_P)$. For each integer $n\ne0$, set
$d_n=\dist(\tP,x^n\cdot\tP)/2$, and let $\zeta_n$ denote the point of
intersection of $S(\lambda/2,\tP)$ with the ray originating at $\tP$
and passing through $x^n\cdot\tP$.  Then
$$\vol(N(P,\lambda/2)) =B(\lambda/2)-\vol(\bigcup_{0\ne
n\in\ZZ}K(\lambda/2,\tP,\zeta_n,d_n)).$$ \EndProposition

\Proof 
For each $\gamma\in\Gamma-\{1\}$ let $\zeta_\gamma\in
S(\lambda/2,\tP)$ denote the intersection of $S(\lambda/2,\tP)$ with
the ray starting at $\tP$ which contains $\gamma\tP$, and let
$d_\gamma = \dist(\tP,\gamma\tP)/2$.  Let $\mathcal D$ denote the
Dirichlet domain for $\Gamma$ centered at the point $\tP$.  Then, by
definition, we have
$$\mathcal D = \bigcap_{1\not=\gamma\in\Gamma}H(\tP,\zeta_\gamma,d_\gamma) .$$
Since 
$\ff_M(P)\ge\lambda$ we have $d_\gamma \ge \lambda/2$ unless
$\gamma = x^n$ for some $n\in \ZZ$.   Thus
$$\begin{aligned}(\inter \mathcal D) \cap \overline{N}(\tP,\lambda/2) 
&= \overline{N}(\tP,\lambda/2) -
 \bigcup_{1\not=\gamma\in\Gamma}K(\lambda/2,\tP,\zeta_\gamma,d_\gamma)\cr
&= \overline{N}(\tP,\lambda/2) - \bigcup_{0\ne n\in\ZZ}K(\lambda/2,\tP,\zeta_n,d_n).
\end{aligned}$$
The set $(\inter \mathcal D) \cap \overline{N}(\tP,\lambda/2)$ has the same
volume as $(\inter \mathcal D) \cap N(\tP,\lambda/2)$, which in turn
is isometric to an open subset of  full measure in 
$N(P,\lambda/2)$.  Thus the result follows.
\EndProof

\Number\label{theta} Now let $z_0$ be a point of $\HH^3$, and let $\D$
and $R$ be real numbers with $0<\D<R$. We may regard
$\overline{K(R,z_0,\zeta,\D)}\cap S(R,z_0)$ as a metric ball about $\zeta$ in
$S(R,z_0)$, with respect to the spherical metric. The radius 
$\Theta$  of this metric ball is the angle formed at $z_0$ between
$\eta_\zeta$ and a segment joining $z_0$ to any point of the boundary of
the topological disk $\overline{K(R,z_0,\zeta,\D)}\cap S(R,z_0)$. We have
 $0<\Theta<\pi/2$.  Note that $\Theta$ is an angle in a
hyperbolic right triangle with hypotenuse $R$, in which the other side
adjoining the angle has length $\D$. Hence 
$$\cos\Theta = \frac{\tanh \D}{\tanh R}$$

It follows that $\Theta=\Theta(\D,R)$, 
where $\Theta(\D,R)$ is the real-valued function with domain
$$\{(\D,R):0<\D<R\}\subset{\RR}^2$$
defined by
$$\Theta(\D, R) = \arccos\left(\frac{\tanh \D}{\tanh R}\right),$$
which takes values in $(0,\pi/2)$.
\EndNumber

\Number\label{bloomsday} 
Let $z_0$ be a point in $\HH^3$, let $R>0$ be a real number, let $\zeta$
and $\zeta'$ be points of $S(R,z_0)$, let $\alpha$ denote the spherical distance from
$\zeta$ to $\zeta'$, and let $\D$ and $\D'$ be positive
numbers less than $R$. Set $K=K(R,z_0,\zeta,\D)$ and $K'=K(R,z_0,\zeta',\D')$. If 
 $\alpha\le
\Theta(\D,R)-\Theta(\D',R)
$  the spherical triangle inequality implies that
$\overline{K'}\cap S(R,z_0)\subset\overline{ K}\cap S(R,z_0)$. It follows easily that
$K'\subset K$.
\EndNumber

\Number\label{omigosh} We define a function $\kappa$ on
$(0,\infty)\times(0,\infty)\subset{\RR}^2$ by defining $\kappa(R,\D)$
to be the volume of $K(R,z_0,\zeta,\D)$, where $z_0$ and $\zeta$ are points in
$\HH^3$ separated by a distance $R$.  Note that if we fix any $R>0$,
the function $\kappa(R,\cdot)$ is monotone decreasing (in the weak
sense) on $(0,\infty)$, and takes the value $0$ when $\D\ge R$.

Next we define functions $\iota$ and $\sigma$  which give,
respectively, the volume of the intersection and the union of two
caps.  More precisely, we set
$$
\begin{aligned}
\iota(R,\D,\D',\alpha) &= \vol (K(R,z_0,\zeta,\D)\cap K(R,z_0,\zeta',\D'));\\
\sigma(R,\D,\D',\alpha) &= \vol (K(R,z_0,\zeta,\D)\cup K(R,z_0,\zeta',\D')),
\end{aligned}
$$
where $z_0$ is a point in $\HH^3$ and $\zeta$ and $\zeta'$ are points in
$S(R,z_0)$ separated by a spherical distance $\alpha$.
We regard $\iota$ and $\sigma$ as real-valued functions with
domain $(0,\infty)^3\times[0,\pi]\subset{\RR}^4$.
Note that
$\iota$ is symmetric in the second and third variables, and that
\Equation\label{bouquiniste}
\sigma(R,\D,\D',\alpha)=\kappa(R,\D)+\kappa(R,\D')-\iota(R,\D,\D',\alpha)
\EndEquation
for any $(R,\D,\D',\alpha)\in{\mathcal D}$.
\EndNumber

In the Appendix (Section \ref{cap comps}) we will give a
formula for the function $\kappa$ and a numerical procedure for
calculating $\iota$ and $\sigma$. These will be used in Sections
\ref{Numbers} and \ref{short section}.

The main result of this section, Proposition \ref{quai voltaire},
gives some monotonicity properties of the function $\sigma$ that will
be needed in Section \ref{New section}.  The proof of Proposition
\ref{quai voltaire} will involve the following technical lemma.

\Lemma\label{i was just sitting in a refrigerator, minding my own
  business}
 Let $z$ be a point of $\HH^3$, let $R$ be a positive number,
and let $\zeta_0$, $\zeta_1$ and $\zeta_2$ be points lying on a great
circle of $S(R,z)$.  Let $\D$ and $\D'$ be numbers such that $0<\D\le
\D'<R$.  For $i=1,2$, let $\alpha_i$ denote the spherical distance from
$\zeta_0$ to $\zeta_i$.  Suppose that $0 \le \alpha_1 < \alpha_2 \le
\pi$ and that the spherical distance from $\zeta_1$ to $\zeta_2$ is
$\alpha_2 - \alpha_1$.  Assume that $\alpha_2 - \alpha_1 <
\Theta(\D',R),$
where $\Theta$ is the function defined in \ref{theta}.
Set $K_0=K(R,z,\zeta_0,\D)$, and $K_i=K(R,z,\zeta_i,\D')$
for $i=1,2$. Set $Y=K_1\cap K_2$, and
$X_i=K_i-Y$ for $i=1,2$.  Then either $X_1\subset K_0$ or $X_2\cap
K_0=\emptyset$.
\EndLemma

\Proof We set $\theta=\Theta(\D,R)$ and
$\theta'=\Theta(\D',R)$, where $\Theta$ is the function defined
in \ref{theta}. Since $0<\D\le \D'<R$, we have
$0<\theta'\le\theta<\pi/2$. According to the hypothesis we have
$0<\alpha_2-\alpha_1<\theta'$.

We set $B=N(R,z)$ and $S=S(R,z)$.  Since the points $\zeta_0$,
$\zeta_1$ and $\zeta_2$ all lie on a great circle, the points $z$,
$\zeta_0$, $\zeta_1$ and $\zeta_2$ lie on a plane $W$. For $i=0,1,2$
the planes $\Pi_i = \Pi(z,\zeta_i,\D')$, are all perpendicular to $W$.
We will reduce the proof of the lemma to a $2$-dimensional
argument by considering the intersections of various sets with
$W$.

We set $\Delta = B\cap W$, $C = S\cap W$ and, for $i = 0,1,2$, we set
$k_i=K_i\cap W$ and $\pi_i = \Pi_i\cap W$.  Let $p:\HH^3\to W$ denote
the perpendicular projection.  Then we have $K_i=B\cap p^{-1}(k_i)$
for $i=0,1,2$.  We also set $\eta=k_1\cap k_2$, and $\xi_i=k_i-\eta$
for $i=1,2$; then we have $X_i=B\cap\pi^{-1}(\xi_i)$ for $i=1,2$.

Let us orient the circle $C$ in such a way that for $i=1,2$, the
clockwise angle from $\zeta_0$ to $\zeta_1$ is $\alpha_i$.  For
$i=0,1,2$ let $A_i$ denote the arc $k_i\cap C$ and let
$\chi_i$ denote the chord $\pi_i\cap \Delta$.  Note that the arc $A_0$
subtends an angle of $2\theta < \pi$ while $A_i$ subtends an angle
$2\theta' < \pi$ for $i = 1,2$.  For $i=0,1,2$ let $l_i$ and $r_i$
denote the endpoints of $A_i$, where $l_i$ is the initial endpoint
when $A_i$ is described in the clockwise direction. 

The clockwise angle from $l_1$ to $r_2$ is $\delta\doteq \alpha_2 -
\alpha_1 + 2\theta' < 4\theta'$, so the arcs $A_1$ and $A_2$ overlap
in a single non-degenerate sub-arc, and $A_1\cup A_2$ is an arc which
subtends an angle $\delta$.  When the arc $A_1\cup A_2$ is described
in the clockwise direction, the points of $T\doteq \partial
A_1\cup\partial A_2$ appear in the order $l_1, l_2, r_1, r_2$.

The clockwise angle from $l_0$ to $l_1$ is
$\gamma\doteq\theta+\alpha_1-\theta'$, and we have
$\gamma+\delta=\theta + \alpha_2 + \theta' < 2\pi$. Hence if we
describe $C$ in the clockwise direction starting at $l_0$, the points
of $T\cup\{l_0\}$ appear in the order $l_0, l_1, l_2, r_1, r_2$.  In
particular, the set $A_0\cap T$ consists of the terms of an initial
subsequence of $(l_1,l_2,r_1,r_2)$. It follows that the set $(C-A_0)\cap
T$ consists of the terms of a final subsequence of
$(l_1,l_2,r_1,r_2)$.

Since the arcs $A_1$ and $A_2$ overlap, the chords $\chi_1$ and
$\chi_2$ cross in a point $Q$ in the interior of $\Delta$.  We
distinguish two cases in the proof of the lemma, depending on whether
or not $Q$ lies in $k_0$.  The two cases are illustrated in the
figure below.

\begin{figure}[h]
\begin{picture}(0,0)%
\includegraphics{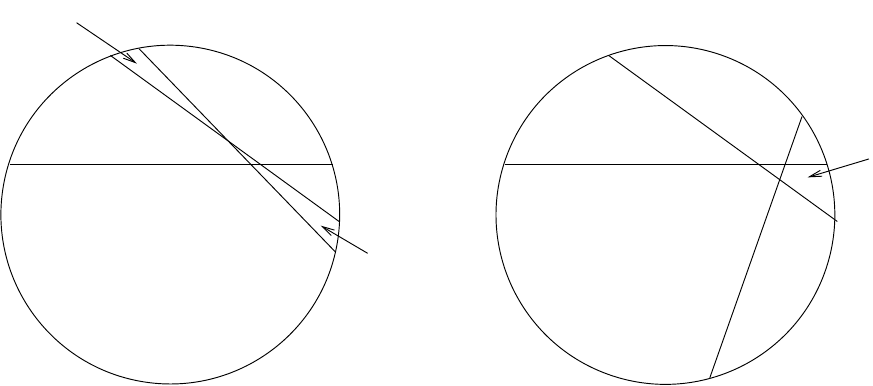}%
\end{picture}%
\setlength{\unitlength}{3947sp}%
\begingroup\makeatletter\ifx\SetFigFont\undefined%
\gdef\SetFigFont#1#2#3#4#5{%
  \reset@font\fontsize{#1}{#2pt}%
  \fontfamily{#3}\fontseries{#4}\fontshape{#5}%
  \selectfont}%
\fi\endgroup%
\begin{picture}(6975,3068)(112,-2397)
\put(5611,-2017){\makebox(0,0)[lb]{\smash{{\SetFigFont{12}{14.4}{\familydefault}{\mddefault}{\updefault}{\color[rgb]{0,0,0}$k_2$}%
}}}}
\put(4320,-849){\makebox(0,0)[lb]{\smash{{\SetFigFont{12}{14.4}{\familydefault}{\mddefault}{\updefault}{\color[rgb]{0,0,0}$k_0$}%
}}}}
\put(5058,-49){\makebox(0,0)[lb]{\smash{{\SetFigFont{12}{14.4}{\familydefault}{\mddefault}{\updefault}{\color[rgb]{0,0,0}$k_2$}%
}}}}
\put(5734,-111){\makebox(0,0)[lb]{\smash{{\SetFigFont{12}{14.4}{\familydefault}{\mddefault}{\updefault}{\color[rgb]{0,0,0}$\xi_1$}%
}}}}
\put(7087,-603){\makebox(0,0)[lb]{\smash{{\SetFigFont{12}{14.4}{\familydefault}{\mddefault}{\updefault}{\color[rgb]{0,0,0}$\eta$}%
}}}}
\put(6226,-1648){\makebox(0,0)[lb]{\smash{{\SetFigFont{12}{14.4}{\familydefault}{\mddefault}{\updefault}{\color[rgb]{0,0,0}$\xi_2$}%
}}}}
\put(3076,-1411){\makebox(0,0)[lb]{\smash{{\SetFigFont{12}{14.4}{\familydefault}{\mddefault}{\updefault}{\color[rgb]{0,0,0}$\xi_2$}%
}}}}
\put(1951,-211){\makebox(0,0)[lb]{\smash{{\SetFigFont{12}{14.4}{\familydefault}{\mddefault}{\updefault}{\color[rgb]{0,0,0}$\eta$}%
}}}}
\put(1126,-136){\makebox(0,0)[lb]{\smash{{\SetFigFont{12}{14.4}{\familydefault}{\mddefault}{\updefault}{\color[rgb]{0,0,0}$k_1$}%
}}}}
\put(2251,-1111){\makebox(0,0)[lb]{\smash{{\SetFigFont{12}{14.4}{\familydefault}{\mddefault}{\updefault}{\color[rgb]{0,0,0}$k_2$}%
}}}}
\put(526,539){\makebox(0,0)[lb]{\smash{{\SetFigFont{12}{14.4}{\familydefault}{\mddefault}{\updefault}{\color[rgb]{0,0,0}$\xi_1$}%
}}}}
\put(301,-811){\makebox(0,0)[lb]{\smash{{\SetFigFont{12}{14.4}{\familydefault}{\mddefault}{\updefault}{\color[rgb]{0,0,0}$k_0$}%
}}}}
\end{picture}%
\end{figure}

If $Q\in k_0$ then each of the chords $\chi_1,\chi_2$ has at least one
endpoint in $A_0$. In particular the set $A_0\cap T$, which consists
of the terms of an initial subsequence of $(l_1,l_2,r_1,r_2)$, has
cardinality at least $2$. Hence the arc $A_0$ contains $l_1$ and
$l_2$, which are the endpoints of the arc $\beta_1 = \xi_1\cap C$.
Moreover, since the arc $\beta_1$ subtends an angle $\alpha_2 -
\alpha_1 < \theta' < \pi/2$ while the arc $A_0$ subtends an angle
$2\theta<\pi$, the arc $A_0$ cannot contain $C-\beta_1$.  It therefore
follows that $\beta_1\subset A_0$.  Since $k_0$ is convex and since
$\xi_1$ is the convex hull of $\beta_1\cup\{Q\}$, we have $$X_1=B\cap
p^{-1}(\xi_1)\subset B\cap p^{-1}(k_0)=K_0,$$ and the lemma is proved
in this case.  If $Q\notin k_0$ then each of the chords
$\chi_1,\chi_2$ has at least one endpoint in $C-A_0$.  In particular
the set $(C-A_0)\cap T$, which consists of the terms of a final
subsequence of $(l_1,l_2,r_1,r_2)$, has cardinality at least
$2$. Hence the endpoints $r_1$ and $r_2$ of the arc $\beta_2 =
\xi_2\cap C$ are contained in $C-A_0$.  The arc $C-\beta_2$ subtends
an angle $2\pi - (\alpha_2 - \alpha_1) > 2\pi - \theta$, while the arc
$C - A_0$ subtends an angle $2\pi - 2\theta < 2\pi - \theta$.  It
follows that $C - \beta_2$ cannot be contained in $C - A_0$, and hence
that $\beta_2$ is contained in $C - A_0$.
 
Since $\Delta-k_0$ is convex, and since $\xi_2$ is the convex hull of
$\beta_2\cup\{Q\}$, we have $\xi_2\subset \Delta-k_0$.  We now have
$$X_2=B\cap\pi^{-1}(\xi_2)\subset B\cap\pi^{-1}(\Delta-k_0)=B-K_0,$$
and the lemma is proved in this case as well.
\EndProof

\Proposition\label{quai voltaire}
For any $R, \D>0$ and $\alpha\in[0,\pi]$, the function $\sigma$ is
monotone decreasing in its third variable and monotone increasing in
its fourth variable (in the sense of \ref{monotonicity}).
\EndProposition

\Proof The first assertion is easy, because if $0<\D_1'<\D_2'$, if $\zeta$
and $\zeta'$ are points of $S(R,z_0)$ whose spherical distance is
$\alpha$, and if we set $K=K(R,z_0,\zeta,\D)$, and $K_i'=K(R,z_0,\zeta,\D_i')$
for $i=1,2$, then we have $K_2'\subset K_1'$ and hence $K_2'\cup
K\subset K_1'\cup K$, so that
$\sigma(R,\D,\D_2',\alpha)\le\sigma(R,\D,\D_1',\alpha)$. 

To prove the second assertion we must show that
$\sigma(R,\D,\D',\cdot)$ is monotone increasing on $[0,\pi]$ for any
$R, \D, \D'\in (0,\infty)$.  We may assume that $\D \le \D'$, since
$\sigma$ is obviously symmetric in its second and third arguments. In
view of \ref{bouquiniste}, it suffices to prove that
$\iota(R,\D,\D',\cdot)$ is monotone decreasing. We may assume that
$\D'<R$, since otherwise the function $\iota(R,\D,\D',\cdot)$ is
identically zero.

We set $\theta=\Theta(\D,R)$ and $\theta'=\Theta(\D',R)$,
where $\Theta$ is the function defined in \ref{theta}.  It clearly
suffices to show that if $\alpha_1$ and $\alpha_2$ satisfy
$0\le\alpha_1<\alpha_2\le\pi$ and $\alpha_2-\alpha_1<\theta'$, then
$\iota(R,\D,\D',\alpha_2)\le\iota(R,\D,\D',\alpha_1)$.

Let $z_0$ be a point of $\HH^3$, and let $\zeta_0$, $\zeta_1$ and
$\zeta_2$ be points of $S(R,z_0)$ such that the spherical distance
from $\zeta_0$ to $\zeta_i$ is $\alpha_i$ for $i=1,2$, and such that
$\zeta_1$ lies on the great circular arc of length $\alpha_1$ with
endpoints $\zeta_0$ and $\zeta_2$. We set $K_0=K(R,z_0,\zeta_0,\D)$ and
$K_i=K(R,z_0,\zeta_i,\D')$ for $i=1,2$.  Then for $i=1,2$ we have
$$\vol K_i=\kappa(R,\D')$$
and
$$\vol(K_0\cap K_1)=\iota(R,\D,\D',\alpha_i).$$

Let us also set $Y=K_1\cap K_2$, and $X_i=K_i-Y$ for $i=1,2$.  Then
$K_i$ is set-theoretically the disjoint union of $X_i$ with $ Y$.  Hence
$$\begin{aligned}
\iota(R,\D,\D',\alpha_1)-\iota(R,\D,\D',\alpha_2)
= &\vol(K_0\cap K_1) - \vol(K_0\cap K_2) \\
= &\vol(K_0\cap X_1)+\vol(K_0\cap Y)\\
   &-\vol(K_0\cap X_2)-\vol(K_0\cap Y),
\end{aligned}$$
i.e.
\Equation\label{don't talk with you mouth full, charlie}
\iota(R,\D,\D',\alpha_1)-\iota(R,\D,\D',\alpha_2)=
\vol(K_0\cap X_1)-\vol(K_0\cap X_2).
\EndEquation
On the other hand, for $i=1,2$ we have
$$\vol X_i=\vol K_i-\vol Y=\kappa(R,\D')-\vol Y.$$
In particular,
\Equation\label{lugar, sam lugar}
\vol X_1=\vol X_2.
\EndEquation

Now according to Lemma \ref{i was just sitting in a refrigerator,
  minding my own business}, we have either $X_1\subset K_0$ or
$X_2\cap K_0=\emptyset$. If $X_1\subset K_0$, then from (\ref{don't
  talk with you mouth full, charlie}) and (\ref{lugar, sam lugar}) we
find that
$$\begin{aligned}\iota(R,\D,\D',\alpha_1)-\iota(R,\D,\D',\alpha_2)&=
\vol( X_1)-\vol(K_0\cap X_2)\cr
&=\vol( X_2)-\vol(K_0\cap X_2)\cr
&\ge0.\end{aligned}$$
On the other hand, if
$X_2\cap K_0=\emptyset$, then from (\ref{don't talk
  with you mouth full, charlie}) we find that
$$\iota(R,\D,\D',\alpha_1)-\iota(R,\D,\D',\alpha_2)= \vol( K_0\cap
X_1)\ge0.$$ 
Thus in both cases we have
$\iota(R,\D,\D',\alpha_2)\le\iota(R,\D,\D',\alpha_1)$, as required.
\EndProof

\section{Nearby volume}\label{New section}

Suppose that $M $ is a closed orientable hyperbolic $3$-manifold.
The goal of this section is to prove a technical result, Lemma
\ref{nested cubes}, which gives a lower bound for the volume of the
metric neighborhood $N(P,\lambda/2)\subset M$, where $\lambda$ is a
positive real number, $P$ is a point of $\XX_M$ such that
$\ff_M(P)\ge\lambda$, and certain inequalities are satisfied.
Proposition \ref{c.d. rivington} expresses $\vol(N(P,\lambda/2))$ in
terms of the volume of a union of caps associated to certain elements
of the cyclic group $C_P$.  The hypotheses of Lemma \ref{nested cubes}
ensure that the caps associated to powers of a generator of $C_P$ with
exponents greater than $3$ are empty, which leads to an estimate for
$\vol N(P, \lambda/2)$ involving the first three powers of the
generator.  An additional hypothesis implies that the caps
corresponding to the third powers are contained in the union of those
associated to the first and second powers.  This leads to an estimate
for $\vol N(P, \lambda/2)$ which involves only first and second powers
of generators of $C_P$.

We begin by introducing some conventions that will be used for the
statement and subsequent applications of Lemma \ref{nested cubes}.

\Number\label{monster dread our damages}
For each integer $n\ge1$, we define a function $\Phi_n$ on the domain
$$\{(\delta ,D):0<\delta \le D\}\subset{\RR}^2$$
by
$$\Phi_n(\delta ,D)=\arccosh\bigg(\cosh(n\delta )+\frac{(\cosh(n\delta
)-1)(\cosh D-\cosh \delta  )}{\cosh \delta +1 }\bigg).$$
\EndNumber

\Lemma\label{monotonicity of Phi} For each integer $n\ge 1$ the
function $\Phi_n$ is increasing in each of the variables $\delta$ and
$D$ (in the sense of \ref{monotonicity}).
\EndLemma

\Proof
It is clear that $\Phi$ is increasing in the variable $D$. To prove
that it is monotone increasing in $\delta$, we fix an $n\ge1$ and a
$D>0$, and set $A=\cosh D$. For any $\delta\in(0,D)$ we have
$\Phi_n(\delta,D)=\arccosh(\frac12f(e^\delta)$, where $f$ is
the function defined on $(1,e^D)$ by
$$f(u)=u^n+u^{-n}+\frac{(u^n+u^{-n}-2)(2A-u-u^{-1})}{u+u^{-1}+2}.$$
It therefore suffices to show that $f$ is monotone increasing on
$(1,e^D)$.

Simplifying, we find that
$$f(u)=\frac{(2A+2)(u^n+u^{-n}-2)}{u+u^{-1}+2}+2.$$
This gives
$$f'(u)=\frac{2A+2}{(u+u^{-1}+2)^2}g(u),$$ where
$$g(u)=(n-1)(u^n-u^{-n-2})+(n+1)(u^{n-2}-u^{-n})+
2n(u^{n-1}-u^{-n-1})+2(1-u^{-2}).$$
In this expression for $g(u)$, it is clear that each term is positive
when $u>1$ and $n\ge1$. Hence $f'$ is positive on $(1,\infty)$.
\EndProof

\Lemma\label{nth power} Let $\delta $ be a positive number, and let
$\gamma$ be a loxodromic isometry of $\HH^3$ whose translation length
is at least $\delta $. Then for any point $z\in\HH^3$ and each integer
$n\ge1$, we have
$$\dist(z,\gamma^n\cdot z)\le \Phi_n(\delta,\dist(z,\gamma\cdot z)).$$
\EndLemma

\Proof
Let $R$ denote the distance from $z$ to the axis of $\gamma$ and
set $D = \dist(z, \gamma\cdot z)$.   If $\gamma$ has translation
length $l$ and twist angle $\theta$ then
$R = \omega(l,\theta,D)$ (see \ref{omegastuff}).
In particular we have
$$\sinh^2R \ge \frac{\cosh D - \cosh l}{\cosh l + 1} .$$
Now since $\gamma^n$ has translation length
$nl$ and twist angle $n\theta$, we have
$$
\begin{aligned}
\cosh(\dist(z, \gamma^n\cdot z))
&=   \cosh(nl) + (\cosh(nl) - \cos(n\theta))\sinh^2R\\
&\ge \cosh(nl) + (\cosh(nl) - 1)\sinh^2R\\
&\ge   \cosh(nl) + \frac{(\cosh(nl) - 1)(\cosh D - \cosh l)}{\cosh l + 1}\\
&= \cosh \Phi_n(l, D)\\
&\ge \cosh \Phi_n(\delta, D),
\end{aligned}
$$
where the last inequality follows from Lemma \ref{monotonicity of Phi}.
\EndProof

\Lemma\label{Phish sandwich}
If $n$ is a positive integer and if $\delta$ and $D$ are real numbers
with $0 < \delta \le  D$, then $ n\delta \le  \Phi_n(\delta,D) \le  nD$.
\EndLemma

\Proof The inequality $ n\delta \le  \Phi_n(\delta,D)$ is immediate from
the definition. To prove the other inequality, note that since $0 <
\delta \le  D$, there exist a loxodromic isometry $\gamma$ of $\HH^3$
with translation length $\delta $, and a point $z\in\HH^3$ such that
$\dist(z,\gamma\cdot z)=D$. The triangle inequality implies that
\Equation\label{bis di yiden veln kummen} \dist(z,\gamma^n\cdot z) \le 
nD.  \EndEquation Since Lemma \ref{nth power} gives
$\Phi_n(\delta,\dist(z,\gamma\cdot z))\le \dist(z,\gamma^n\cdot z)$,
it follows that \Equation\label{i had tsuris with some flower}
\Phi_n(\delta,D) \le nD.  \EndEquation  \EndProof 

\Number\label{hyperbolic lore of cosines}
If $x$ and $y$ are real numbers  with $0< y\le2x$, we have
$$0<(\coth x)(\coth y-\cosech y)\le(\coth x)(\coth 2x-\cosech 2x)=1.$$
Hence on the domain
$$\{(x,y):0< y\le2x\}\subset{\RR}^2,$$
we may define a function $\Psi$, with values in $[0,\pi/2]$, by
$${\Psi}(x,y)=\arccos((\coth x)(\coth y-\cosech y)).$$ Since $\coth
x$ is monotone decreasing for $0<x<\infty$, and $\coth y-\cosech
y$ is monotone increasing for $0<y<\infty$, the function
${\Psi}$ is monotone increasing in its first argument and monotone
decreasing in its second.

If an isosceles hyperbolic triangle has base $y$ and has its other two
sides equal to $x$, the triangle inequality gives $y\le2x$, and the
hyperbolic law of cosines shows that the base angles are equal to
${\Psi}(x,y)$.
\EndNumber

The statement of the following lemma involves the function $\Theta$
which was defined in \ref{theta}, as well as the functions $\kappa$,
and $\sigma$ which were defined in \ref{omigosh}.

\Lemma\label{nested cubes} 
Let $M$ be a closed,
orientable hyperbolic $3$-manifold and let $P$  be a point of $\XX_M$.
Suppose that $\delta$ and $\lambda$ are constants with
$0 <  \lambda < 4\delta$.
Assume that
\Bullets
\item  $\ff_M(P) \ge \lambda$, and that
\item the conjugacy class of a generator of $C_P$ is represented
  by a closed geodesic of length at least $\delta$.
\EndBullets 
Set $D = D_M(P)$. Then $T_n=\Phi_n(\delta,D)$ is defined for every
$n\ge1$. Furthermore, 
 $(D,T_2)$ lies in the domain of $\Psi$, and we have
$$\vol N\big(P,\frac\lambda2\big)\ge
B\big(\frac\lambda2\big)-2\sigma\big(\frac\lambda2, \frac{ D }2,
\frac{T_2}2,{\Psi}( D ,T_2)\big)
-2\kappa\big(\frac\lambda2,\frac{T_3}2\big).$$ 
If in addition we have
$D < T_3 < \lambda$, 
so that, in particular,  the quantities
$\Theta(\frac{D}{2},\frac{\lambda}{2})$ and $\Theta(\frac{T_3}{2},\frac{\lambda}{2})$ are defined,
and if 
\Equation\label{grand duke}
\cos\big(\Theta\big(\frac{ D }2,\frac\lambda2\big)-
\Theta\big(\frac{T_3}2,\frac\lambda2\big)\big) < \frac{\cosh  D  \cosh
T_3-\cosh 2 D }{\sinh  D  \sinh T_3},
\EndEquation
then
$$\vol N\big(P,\frac\lambda2,\big)\ge
B\big(\frac\lambda2\big)-
2\sigma\big(\frac\lambda2, \frac{ D }2, \frac{T_2}2,{\Psi}( D ,T_2)\big).$$ 
\EndLemma

\Proof
We write $M=\HH^3/\Gamma$, where $\Gamma\le\isomplus(\HH^3)$ is
discrete, torsion-free and cocompact.
Set $q=q_\Gamma$ (\ref{same to me}), and choose $\tP\in
  q^{-1}(P)$. We use the
base point $\tP\in\HH^3$ to identify $\pi_1(M,P)$ with
$\Gamma$ (see \ref{same to me}). In particular, $C_P$ is identified
with a subgroup of $\Gamma$. We fix a generator $x$ of $C_P$.

The hypothesis implies that $x$ has translation length at least
$\delta$. In particular we have $D\ge\delta$, so that
$T_n=\Phi_n(\delta,D)$ is defined for every $n\ge1$.

According to Lemma \ref{Phish sandwich} we have $T_2 = \Phi_2(\delta,
 D ) \le 2 D $.  Thus $( D ,T_2)$ lies in the domain of ${\Psi}$.

For each integer $n\ne0$ we set $d_n=\dist(\tP,x^n\cdot\tP)$. We
observe that $d_{-n}=d_n$ for each $n\ne0$ and that $d_1=D$.
Moreover, since $x$ has translation length at least $\delta$, it
follows from Lemmas \ref{Phish sandwich} and \ref{nth power}, and the
triangle inequality, that
\Equation\label{bye bye love}
|n|\delta\le \Phi_{|n|}(\delta, D )\le d_n\le |n| D
\EndEquation
for every integer $n\ne0$, and in particular,
\Equation\label{i hear you knockin'}
 d_n\ge T_n\quad{\rm for}\  n>0.
\EndEquation

We let $\zeta_n$ denote the point of intersection of $S(\lambda/2,\tP)$
with the ray originating at $\tP$ and passing through $x^n\cdot\tP$.  We
set
$$\calk_n=K(\lambda/2,\tP,\zeta_n,d_n/2)$$
for each integer $n\ne0$,
$$\cals=\bigcup_{0\ne
  n\in\ZZ}K(\lambda/2,\tP,\zeta_n,d_n/2),$$ and
$$\cals_N=\bigcup_{0<|n|\le N}K(\lambda/2,\tP,\zeta_n,d_n/2)$$
for each
integer $N>0$.  
Since $\ff_M(P)
  \ge \lambda$, it follows from Proposition \ref{c.d. rivington} that
\Equation\label{cherry blossoms}
\vol(N(P,\lambda/2))=B(\lambda/2)-\vol\cals.
\EndEquation 

The estimates in the conclusion of the lemma will be deduced via
(\ref{cherry blossoms}) from suitable estimates for the volume of
$\cals$. As a first step we shall estimate the volume of $\cals_2$.
Note that
\Equation\label{onion soup}
\cals_2=\calt_{+1}\cup\calt_{-1},
\EndEquation
where
$$\calt_\epsilon=K(\lambda/2,\tP,\zeta_\epsilon,
D/2)\cup
K(\lambda/2,\tP,\zeta_{2\epsilon},{d_{2
}}/2)$$
for $\epsilon=\pm1$.

For $\epsilon=\pm1$ we consider the hyperbolic triangle with vertices
$\tP$, $x^{\epsilon}\cdot\tP$ and $x^{2\epsilon}\cdot\tP$. The sides
adjacent to the vertex $\tP$ have lengths $d_1={  D }$ and
$d_{2}$. The side opposite $\tP$ has length
$\dist(x^{\epsilon}\cdot\tP,x^{2\epsilon}\cdot\tP)=\dist(\tP,x^{\epsilon}\cdot\tP)={
 D }$.  Let $\alpha$ denote the angle at the vertex $\tP$.
According to the discussion in \ref{hyperbolic lore of cosines} we
have $d_2\le2{  D }$ and $\alpha={\Psi}({  D },d_2)$.
By (\ref{i hear you knockin'}) we have $d_2\ge T_2$. Hence the
monotonicity properties of ${\Psi}$ pointed out in
\ref{hyperbolic lore of cosines} give
$$\alpha\le{\Psi}({  D },T_2).$$ 

Our definitions of $\alpha$ and of the $\zeta_n$ imply that $\alpha$
is the spherical distance between $\zeta_\epsilon$ and
$\zeta_{2\epsilon}$. In view of the definition of $\calt_\epsilon$,
and the definition of the function $\sigma$ given in \ref{omigosh}, it
follows that $\vol(\calt_\epsilon)=\sigma(\lambda/2,{ D
}/2,d_2/2,\alpha)$ for $\epsilon=\pm1$. Since $\alpha\le{\Psi}({ D
},T_2)$, and since $d_2\ge T_2$, it follows from Proposition \ref{quai
  voltaire} that
$$\vol(\calt_\epsilon)\le\sigma(\lambda/2,{
 D }/2,T_2/2,{\Psi}({  D },T_2))$$ for $\epsilon=\pm1$. But
from (\ref{onion soup}) we have
$$\vol\cals_2=\vol(\calt_{+1}\cup\calt_{-1})\le\vol(\calt_{+1})+\vol(\calt_{-1}),$$
and hence 

\Equation\label{ralph branca}
\vol\cals_2\le2\sigma(\lambda/2,{  D }/2,T_2/2,{\Psi}({  D },T_2)).
\EndEquation

We now turn to the estimation of $\vol\cals$.  Since by hypothesis we
have $4\delta>\lambda$, it follows from (\ref{bye bye love}) that
$d_n>\lambda$ for any $n\ge4$.  In view of the remarks in \ref{basic
caps}, it follows that
$$\calk_n=K(\lambda/2,\tP,\zeta_n,d_n/2)=\emptyset\quad{\rm for\ any\ }n{\rm\ with\ }|n|\ge4.$$
Hence 
\Equation\label{princess ida}
\cals=\cals_3.
\EndEquation

It follows from (\ref{princess ida}) that
\Equation\label{bobby thompson}
\vol(\cals)=\vol(\cals_2\cup\calk_3\cup\calk_{-3}) \le
\vol(\cals_2)+\vol\calk_3+\vol\calk_{-3}.
\EndEquation
By the definition of the function $\kappa$ given in \ref{omigosh}, and
the monotonicity observed there, together with (\ref{i hear you
knockin'}), we find that
\Equation\label{roy campanella}
\vol\calk_3+\vol\calk_{-3}=
2\kappa(\lambda/2,d_3/2)\le2\kappa(\lambda/2,T_3/2).
\EndEquation
Combining (\ref{bobby thompson}) with (\ref{roy campanella}) and
(\ref{ralph branca}), we deduce that
\Equation\label{duck dodgers}
\vol(\cals)\le2\sigma(\lambda/2,{  D }/2,T_2/2,{\Psi}({
 D },T_2))+2\kappa(\lambda/2,T_3/2).
\EndEquation
 The first assertion of the lemma follows immediately from
(\ref{cherry blossoms}) and (\ref{duck dodgers}).

For the proof of the second assertion we
shall begin by showing that if (\ref{grand
duke}) holds then
\Equation\label{utopia limited}
K(\lambda/2,\tP,\zeta_{3\epsilon},d_{3}/2)\subset
K(\lambda/2,\tP,\zeta_{\epsilon},D/2)\quad{\rm for\ }\epsilon=\pm1. 
\EndEquation
If $d_3\ge\lambda$, then by the remarks in \ref{basic caps} we have
$K(\lambda/2,\tP,\zeta_{3\epsilon},d_{3}/2)=\emptyset$ for $\epsilon=\pm1$, so
that (\ref{utopia limited}) is true in this case. Now suppose that
$d_3<\lambda$.  Note that this implies that $\Theta(d_3/2,\lambda/2)$
is defined.
 
For $\epsilon=\pm1$ we consider the hyperbolic triangle
with vertices $\tP$, $x^{\epsilon}\cdot\tP$ and $x^{3\epsilon}\cdot\tP$. The
sides adjacent to the vertex $\tP$ have lengths $d_1={  D }$ and
$d_{3}$. The side opposite $\tP$ has length
$\dist(x^{\epsilon}\cdot\tP,x^{3\epsilon}\cdot\tP)=
\dist(\tP,x^{2\epsilon}\cdot\tP)=d_2$.
From the hyperbolic law of cosines it follows that the angle $\gamma$
at the vertex $\tP$ is determined by
$$ 
\cos\gamma=\frac{\cosh {  D }\cosh d_3-\cosh d_2}{\sinh {
 D }\sinh d_3}
$$
We have $d_2 \le 2D$ by (\ref{bye bye love}).  Thus 
$$
\cos\gamma\ge
\frac{\cosh {  D }\cosh d_3-\cosh 2D}{\sinh {
 D }\sinh d_3}.
$$

If we set $f(D, d_3)$ equal to the right-hand side of the inequality above, then $$
\begin{aligned}
\frac{\partial f}{\partial d_3}
 &= \frac{\cosh D \sinh D \sinh^2 d_3 -
              \cosh D\sinh D \cosh^2 d_3 +
              \sinh D\cosh 2D\cosh d_3}
             {\sinh^2D\sinh^2 d_3}\\
&= \frac{\sinh D(\cosh(2D)\cosh d_3 - \cosh D)}
             {\sinh^2D\sinh^2 d_3}
\end{aligned}
.$$
Since $\cosh(2D)\cosh d_3 - \cosh D > \cosh(2D) - \cosh D > 0$, the
function $f$ is increasing in $d_3$.
By (\ref{i hear you knockin'}) we
have $d_3\ge T_3$. Hence
$$\cos\gamma\ge\frac{\cosh {  D } \cosh T_3 -\cosh 2{
     D }}{\sinh {  D } \sinh T_3 }.
$$

In view of (\ref{grand duke}) it follows that \Equation\label{only in
  america} \cos(\Theta(D/2,\lambda/2)-\Theta(T_3/2,\lambda/2))
<\cos\gamma.
\EndEquation

It is clear from the definition given in \ref{theta} that
$\Theta$ is monotone decreasing in the first variable.
Since $\Theta$ takes values in $(0,\pi/2)$, and since
$D < T_3$ by hypothesis, we have
\Equation\label{for two cents plain}
 \pi/2>\Theta(D/2,\lambda/2)-\Theta(T_3/2,\lambda/2)>0.
\EndEquation 
Since $0<\gamma<\pi$, it follows from (\ref{only
in america}) and (\ref{for two cents plain})  that
$$\gamma<\Theta(D/2,\lambda/2)-\Theta({T_3}/2,\lambda/2).$$
Since $\Theta$ is monotone decreasing in the first variable  and
since $T_3\le d_3$ by (\ref{i hear you knockin'}),
we deduce that 
\Equation\label{cardinal alfred}
\gamma<\Theta(D/2,\lambda/2)-\Theta(d_3/2,\lambda/2).
\EndEquation

Our definitions of $\gamma$ and of the $\zeta_n$ imply that $\gamma$ is
the spherical distance between $\zeta_\epsilon$ and $\zeta_{3\epsilon}$. It
therefore follows from \ref{bloomsday} and from (\ref{cardinal
  alfred}) that $K(\lambda/2,\tP,\zeta_{3\epsilon},d_3/2)\subset
K(\lambda/2,\tP,\zeta_\epsilon,D/2)$.  This shows that (\ref{utopia
  limited}) holds in this case as well.

Now it follows from (\ref{princess ida}) and (\ref{utopia limited})
that, under the assumption (\ref{grand duke}), we have
\Equation\label{nearsighted mr magoo}
\cals=\cals_3=\cals_2.
\EndEquation
Combining (\ref{nearsighted mr magoo}) with (\ref{ralph branca}), we
deduce that \Equation\label{watch the drivers roll}
\vol\cals\le\sigma(\lambda/2,{  D }/2,T_2/2,{\Psi}({
   D },T_2)).
\EndEquation
The second assertion of the lemma follows immediately from (\ref{cherry
  blossoms}) and (\ref{watch the drivers roll}).
\EndProof

\section{Distant points}\label{distant point section}

The purpose of this section is to adapt some results proved in
\cite{cusp} to the context of the present paper.

\Number\label{chalk cliffs}

We define a set $\calx\subset{\RR^2}$ by
$$\calx=\{(D,\lambda):\frac1{1+e^D}+\frac1{1+e^\lambda}<\frac12\}.$$
For each integer $k>2$  we define a real-valued function $\rho_k$ on $\calx$ by
$$\rho_k(D,\lambda)=\frac12\log\bigg(\frac{k-2}{1/2-1/(1+e^D)-1/(1+e^\lambda)}-1\bigg).$$
Thus for any $(D,\lambda)\in\calx$ we have 
$$\frac{k-2}{1+e^{2\rho_k(D,\lambda)}}+ \frac{1}{1+e^{D}}+\frac{1}{1+e^{\lambda}}
=\frac12.$$ 

We define a real-valued function $\Sigma$ on $[0,\infty)^3$ by
$$\Sigma(h,R_1,R_2)=\arccosh(\sinh R_1\sinh R_2+\cosh R_1\cosh
R_2\cosh h).$$
According to \cite[p. 89]{Fe}, 
given a line $A$ in $\HH^2$ and points $z_1$ and $z_2$ of $\HH^2$
which lie in different components of $\HH^2-A$,
if $R_i$ is the distance from $z_i$ to $A$, and if $h$ is the
distance between the perpendicular projections of $z_1$ and $z_2$ to
$A$, then $\dist(z_1,z_2)=\Sigma(h,R_1,R_2)$.
\EndNumber

\Number\label{breyers calling}It follows that
if $z_1$ and $z_2$ are points in $\HH^3$, if $A$ is a line in $\HH^3$,
if $R_i$ is the distance from
$z_i$ to $A$, and $h$ is the distance between the perpendicular
projections of $z_1$ and $z_2$ to $A$, then
$\dist(z_1,z_2)\le\Sigma(h,R_1,R_2)$.
\EndNumber

\Lemma\label{fracklog too}Let $k>2$ be an integer,  let $M$ be a closed,
 orientable hyperbolic $3$-manifold such that
 $\pi_1(M)$ is $k$-free,  and let $\mu$ be a Margulis
 number for $M$. Let  $P$ be a point of $\XX_M$, and set
 $\lambda=\ff_M(P)$ and $D=D_M(P)$. 
Then we have $(D,\lambda)\in\calx$, and there is a point
$Q\in\Mthick(\mu)$ such that $\dist(P,Q)\ge\rho_k(D,\lambda)$.
\EndLemma

\Proof It follows from the definitions given in \ref{eggs and fries}
that a generator $x_0$ of $C_P$ is represented by a loop at
  $P$ of length
$D$, and that some element $x_1$ of $\Gamma- C_P$ is represented by a
loop  at $P$  of length $\lambda$.

Since $x_1\notin C_P$, the elements $x_0$ and $x_1$ do not commute.
Since $\Gamma\cong\pi_1(M)$ is in particular $2$-free, $x_0$ and $x_1$
are independent (in the sense of \ref{kurosh stuff}).

Since $\Gamma$ is $k$-free, and since there are two independent
elements of $\pi_1(M,P)$ represented by loops of length $D$ and
$\lambda$, it follows from the case $m=2$ of \cite[Corollary
\bigradcor]{cusp} that there is a point $Q\in \Mthick(\mu)$ such that
$\rho=\dist_M(P,Q)$ satisfies \Equation\label{pickle}
\frac{k-2}{1+e^{2\rho}}+ \frac{1}{1+e^{D}}+\frac{1}{1+e^{\lambda}}
\le\frac12.  \EndEquation This implies that $(D,\lambda)\in\calx$, and
that $\dist(P,Q)\ge\rho_k(D,\lambda)$.  \EndProof

\Proposition\label{far beak}Let $k>2$ be an integer,  let $M$ be a closed,
 orientable hyperbolic $3$-manifold such that
 $\pi_1(M)$ is $k$-free,  and let $\mu$ be a Margulis
 number for $M$. Let  $P$ be a point of $\XX_M$, and set
 $\lambda=\ff_M(P)$ and $D=D_M(P)$. Then we have
 $(D,\lambda)\in\calx$. Furthermore, if $s$ is a real number
 such that $\lambda/2\le s\le \rho_k(D,\lambda)$, then there is a point
$Y_s\in\Mthick(\mu)$ such that $\dist(P,Y_s)=s$.
\EndProposition

Note that the second assertion of Proposition \ref{far beak} is
vacuous if $\lambda/2> \rho_k(D,\lambda)$. 

\Proof[Proof of Proposition \ref{far beak}]
The first assertion, that $(D,\lambda)\in\calx$, is included in Lemma
\ref{fracklog too}. 

To prove the second assertion, we first define a continuous
function $\Delta :\Mthick(\mu)\to\RR$ by $\Delta
(Y)=\dist_M(Y,P)$. Since $\mu$ is a Margulis number for $M$, the set
$\Mthick(\mu)$ is connected by \ref{rev}. Since $\Delta$ is
continuous, the set $J=\Delta(\caly)\subset\RR$ is an interval.

According to Lemma \ref{fracklog too}, there is a point
$Q\in\Mthick(\mu)$ such that $\dist(P,Q)\ge\rho_k(D,\lambda)$.
If we set $d=\Delta(Q)$ it follows that $d\ge\rho_k(D,\lambda)$ and
that $d\in J$.

Now recall from \ref{eggs and fries} that there is a loop $\gamma_0$
of length $\ell_P<\ff_M(P)=\lambda$ based at $P$; that $[\gamma_0]$
lies in the maximal cyclic subgroup $C_P$ of $\pi_1(M,P)$; and that
there is a loop $\gamma_1$ of length $\ff_M(P)=\lambda$ based at $P$
such that $[\gamma_1]\notin C_P$. In particular, $[\gamma_0]$ and
$[\gamma_1]$ do not commute. Hence if we denote by $c_i$ the support
of $\gamma_i$ (i.e. $c_i=\gamma_i([0,1])$), and denote by $K$ the
path-connected set $c_0\cup c_1\subset M$, then the inclusion
homomorphism $\pi_1(K,P)\to\pi_1(M,P)$ has a non-abelian image. Since
$\mu$ is a Margulis number for $M$, each component of $\Mthin(\mu)$
has an abelian fundamental group, and hence
$K\not\subset\Mthin(\mu)$.

Let us fix a point $Q'\in K\cap\Mthick(\mu)$, and set
$d'=\Delta(Q')\in J$. The definition of $K$ implies that $Q'$ lies in
the support of a loop of length at most $\lambda$ based at $P$, and
hence that $d'=\dist(P,Q')\le\lambda/2$.

Now if $s$ is a real number such that $\lambda/2\le s\le
\rho_k(D,\lambda)$, then in particular we have $d'\le s\le d$. Since
the interval $J$ contains $d$ and $d'$, it also contains $s$. This
means that there is a point $Y_s\in\Mthick(\mu)$ such that
$s=\Delta(Y_s)=\dist(P,Y_s)$.
\EndProof

\section{The volume of a metric ball}
\label{Boroczky balls}

We have
\Equation\label{bees in my bonnet, pain in my heart}
B(r)=\pi(\sinh(2r)-2r).
\EndEquation

\Number\label{twain}
For $n\ge2$ and for $R>0$ we shall denote by $h_n(R)$ the distance from the barycenter to a vertex of the
regular hyperbolic tetrahedron $\Delta(R)$ with sides of length $2R$.
It is easy to verify, using hyperbolic trigonometry, that
$$\tanh h_2(R) = \frac{2\sinh^2R}{\sqrt{\cosh^2(2R) - \cosh^2R}},$$ 
and
\Equation\label{who would}
\tanh h_3(R) = \frac{2\sinh^2R}{\sqrt{\cosh^2(2R) - \cosh^2h_2(R)}}.
\EndEquation
\EndNumber

\Number\label{you don't remember me}
Let $R$ be any positive number. Consider 
an arbitrary sphere-packing in $\HH^3$ by spheres of
radius $R$. Let $\mathcal D$ denote the Dirichlet domain for
this packing, centered at a point $z\in\HH^3$.  The main result of
B\"or\"oczky's paper \cite{boroczky} states that \Equation\label{i lit
  another lunky} \vol \mathcal D \ge B(R)/d(R),  \EndEquation 
where $d$ is a function of which the
definition  will be reviewed below.

Now let $M$ be a hyperbolic $3$-manifold, 
and suppose that $Y$ is a point in $M$ which is the center of a
hyperbolic ball of radius $R$ (i.e.  $Y$ is a $2R$-thick point).  Let
us write $M=\HH^3/\Gamma$, where $\Gamma\le\Isom(\HH^3)$ is discrete
and torsion-free. Set $q=q_\Gamma$ (\ref{same to me}), and let $\tY$
denote a point in $q^{-1}(Y)$. Then $q^{-1}(Y)$ is the
set of center points for a sphere-packing in $\HH^3$ by spheres of
radius $R$. Applying (\ref{i lit another lunky}) with $z=\tY$, we find
that
$\vol M \ge B(R)/d(R)$.

B\"or\"oczky's proof of (\ref{i lit another lunky}) actually gives a
stronger conclusion.  It is shown in the proof that, in fact,

\Equation\label{cheroot}
\vol(\mathcal D\cap N(z, h_3(R))) \ge B(R)/d(R).
\EndEquation
While (\ref{cheroot}) may not give improved estimates for the
density of a sphere packing, it does have a very natural
interpretation for a sphere-packing defined as above in terms of a
hyperbolic $3$-manifold $M$ and a point $Y\in M$: it is equivalent to the statement
that $\vol N(Y,h_3(R)) \ge B(R)/d(R)$.

The main result of this section, Proposition \ref{boroczky
    bonus},  will give a stronger lower bound for the
volume of $N(Y,h_3(R))$ under the stronger hypothesis that there exists
a point $Q\in M$ which is sufficiently far away from $Y$.
\EndNumber

\Number\label{seven pillows of wisdom}
We now give the definition of B\"or\"oczky's density
function $d$.  Let 
$$\beta(r)=\arcsec(\sech(2r)+2)$$ denote the dihedral
angle of $\Delta(r)$. Let $$\tau(r)=3\int_{\beta(r)}^{\arcsec3}\arcsech((\sec t)-2)dt$$ denote
the volume of $\Delta_3(r)$. Then $$d(r)=(3\beta(r)-\pi)(\sinh(2r)-2r)/\tau(r).$$
\EndNumber

\Number
B\"or\"oczky's inequality (\ref{cheroot}) applies in particular to the case of a sphere-packing
in $\HH^3$ consisting of a single sphere of radius $R$ centered at a
point $z$. In this case we have $\cald=\HH^3$,  so that 
$\vol(\mathcal D\cap N(z, h_3(R)))=\vol(N(z, h_3(R)))=B(h_3(R))$. Thus
(\ref{cheroot}) becomes  
\Equation\label{old oaken frigate}
B(h_3(R)) \ge B(R)/d(R),
\EndEquation
which therefore holds for every $R>0$.
\EndNumber

The statement of Proposition \ref{boroczky bonus} involves a function
$\Vbor$ which we shall now define.

\Definition\label{bonus function def}Let $R$ and $\rho$ be positive real
numbers such that $\rho  > h_3(R)$. We define 
$$\phi_1(R,\rho ) = \arcsin\left(\frac{\sqrt{\cosh^2\rho  - \cosh^2R}}{\sinh
    \rho \cosh R}\right),$$
$$\phi_2(R) = \arcsin\left(\frac{\sqrt{\cosh^2h_3(R) - \cosh^2R}}{\sinh h_3(R)\cosh R}\right) ,$$
$$\phi(R,\rho ) = \phi_1(R,\rho ) - \phi_2(R),$$
 and
$$\Vbor(R,\rho )=
 \left(\frac{1-\cos\phi(R,\rho )}{2}\right)B(h_3(R)) +
 \left(\frac{1+\cos\phi(R,\rho )}{2}\right)B(R)/d(R) .$$
\EndDefinition

\Remark\label{Vborincreases}
For any fixed $R>0$ the function $\phi(R,\cdot)$
  is positive-valued and monotone
increasing on $ (h_3(R),\infty)$. In view of the inequality (\ref{old
  oaken frigate}), it follows that $\Vbor(R,\cdot)$ is also monotone
increasing on $ (h_3(R),\infty)$.
\EndRemark

\Proposition
\label{boroczky bonus}
Let $M$ be a hyperbolic $3$-manifold, let $R$ be a positive number,
and suppose that $Y$ is a $2R$-thick point in $M$.  Suppose that there
exists a point $P\in M$ such that $\rho  \dot= \dist_M(Y,P) > h_3(R)$.
Then
\Equation\label{which twin has the toni}
\vol N(Y,h_3(R)) \ge\Vbor(R,\rho ).
\EndEquation
\EndProposition

\Proof[Proof of Proposition \ref{boroczky bonus}] We set $q=q_\Gamma$
(\ref{same to me}), and we choose a point $\tY \in q^{-1}(Y)$.  Since
$\ell_Y\ge2R$, there is a sphere-packing in $\HH^3$ consisting of
spheres of radius $R$ centered at the points $q^{-1}(Y)$.  Let
$\mathcal D$ be a Dirichlet domain for the sphere centered at
$\tY$. Let $\mathcal B$ denote the ball of radius $R$ centered at
$\tY$ and let $\mathcal B'$ denote the ball of radius $h_3(R)$
centered at $\tY$.  If $X$ is a subset of $\HH^3$ we let $C(X)$ denote
the union of all rays from $\tY$ that contain a point of $X$.

The construction given in \cite[\S 5]{boroczky}, when restricted to
the $3$-dimensional case, begins by decomposing $\mathcal D$ as the
union of $D_0 = C(\mathcal D\cap\partial \mathcal B')$ and $D_1$,
where $D_1$ is the union of the sets of the form $C(F\cap \mathcal
B')$ as $F$ runs over the closed $2$-dimensional faces of $\mathcal
D$.  The set $D_1$ is then further subdivided to obtain a
decomposition of $\mathcal D$ into the union of $D_0$ and a certain
family of $3$-dimensional convex cells.  For each cell $E$ it is shown
that $\vol(\mathcal B\cap E )/\vol(\mathcal B' \cap E) \le d(R)$.
Moreover, since $D_0$ is the cone based at $\tY$ on the subset of
$\partial \mathcal B'$ that is contained in $\mathcal D$, we have
$\vol(\mathcal B\cap D_0 )/\vol(\mathcal B' \cap D_0) = \vol\mathcal
B/\vol \mathcal B' < d(R)$.  Thus 
$\vol(\mathcal B\cap \mathcal D )/\vol(\mathcal B' \cap \mathcal D)
\le d(R)$; since $\mathcal B \subset \mathcal D$ this implies
 B\"or\"oczky's stated result that 
$\vol\mathcal B/\vol\mathcal D \le d(R)$.

We may summarize the discussion above as follows.
If we set $t = \vol D_0/\vol\mathcal B'$, so that $\vol D_0 =
tB(h_3(R))$,
then B\"or\"oczky's argument implies that
\Equation\label{after an elephant}
\vol\mathcal D = \vol D_0 + \vol (\mathcal D - D_0)
\ge tB(h_3(R)) + (1 - t)B(R)/d(R) .
\EndEquation

(It is shown in \cite[Lemma 12]{boroczky} that the vertices of
$\mathcal D$ lie outside $\mathcal B'$, so the set $D_0$ is always
non-empty; that is, we have $t>0$ in (\ref{after an
    elephant}).  This observation does not strengthen B\"or\"oczky's
theorem about general sphere packings, since one has no {\it a priori}
information about the size of the set $D_0$. Our goal
here is to quantify the improvement given by (\ref{after an
    elephant}) in
terms of $\dist(Y,P)$.)

\begin{figure}[ht]
\begin{picture}(0,0)%
\includegraphics{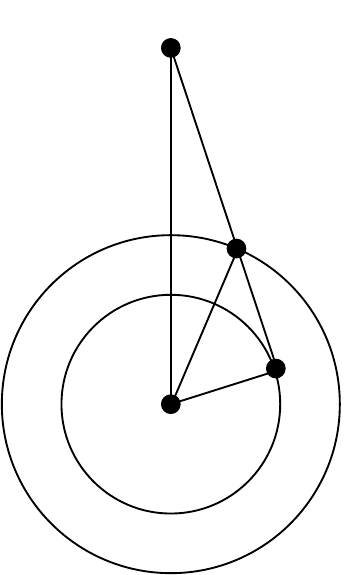}
\end{picture}%
\setlength{\unitlength}{3947sp}%
\begingroup\makeatletter\ifx\SetFigFont\undefined%
\gdef\SetFigFont#1#2#3#4#5{%
  \reset@font\fontsize{#1}{#2pt}%
  \fontfamily{#3}\fontseries{#4}\fontshape{#5}%
  \selectfont}%
\fi\endgroup%
\begin{picture}(2734,4600)(434,-3977)
\put(1696,479){\makebox(0,0)[lb]{\smash{{\SetFigFont{12}{14.4}{\familydefault}{\mddefault}{\updefault}{\color[rgb]{0,0,0}$\widetilde Q$}%
}}}}
\put(2776,-2416){\makebox(0,0)[lb]{\smash{{\SetFigFont{12}{14.4}{\familydefault}{\mddefault}{\updefault}{\color[rgb]{0,0,0}$U$}%
}}}}
\put(2521,-1216){\makebox(0,0)[lb]{\smash{{\SetFigFont{12}{14.4}{\familydefault}{\mddefault}{\updefault}{\color[rgb]{0,0,0}$V$}%
}}}}
\put(1696,-2986){\makebox(0,0)[lb]{\smash{{\SetFigFont{12}{14.4}{\familydefault}{\mddefault}{\updefault}{\color[rgb]{0,0,0}$\widetilde P$}%
}}}}
\end{picture}%
\end{figure}

In order to establish (\ref {which twin has the toni}), it
suffices to show that the quantity
$\vol D_0/\vol\mathcal B'$, which is denoted $t$ in (\ref{after an
  elephant}), is greater than $(1-\cos\phi)/2$, where $\phi$ is the angle
defined in the statement.  For this, it suffices to show that $D_1$
contains $C\cap \mathcal B'$ where $C$ is the convex region bounded by
a circular cone with apex at $\tY$, such that the angle between the
axis and a generating line of $\partial C$ is at least $\phi$.

Let $\tP$ be a point of $\mathcal D$ such that $q(\tP) = P$.  We know
that $\mathcal D$ contains the ball $\mathcal B$.  Since $\mathcal D$
is convex, the entire convex hull $H$ of $\mathcal B\cup \{\tP\}$ is
contained in $\mathcal D$.  We shall show that $H$ contains a conical
region $C$ of angle $\phi$.

Consider a line through $\tP$ which is tangent to $\partial \mathcal
B$ at a point $U$, and let $V$ be the point where the segment
$[\tP,U]$ meets $\partial \mathcal B'$.

We have $\dist(\tY, \tP) = \rho $, $\dist(\tY, V) = h_3(R)$
and $\dist(\tY, U) = R$.  The angle $\angle \tY U\tP$ is a
right angle.  Hence it follows easily from the hyperbolic Pythagorean
Theorem and the hyperbolic law of sines that the measure of $\angle
\tP\tY U$ is $\phi_1$ and the measure of $\angle V\tY U$
is $\phi_2$, where $\phi_1$ and $\phi_2$ are as defined in the
statement.  It is clear that $H$ contains $\mathcal B\cap C$ where $C$ is the
convex region bounded by the circular cone with apex at $\tY$,
whose axis contains $\tP$ and whose boundary contains $V$.  Since
$\phi=\phi_1-\phi_2$ is the angle between the axis and a generator of
the boundary cone of $C$, the result follows.
\EndProof

\section{Margulis numbers and diameter}\label{Margulis section}

\Lemma\label{violet} 
Let $m$ be a positive integer and let $H$ be
finitely generated group which is $m$-free and has rank $\ge m$.
Let $S$ be a finite generating set for $H$ and let $T_0\subset S$ be
an independent set (see \ref{kurosh
  stuff}).  Then there exists an  independent set $T$ with
$T_0\subset T\subset S$ and $|T|=m$.
\EndLemma

\Proof Among all independent sets $T$ such that $T_0\subset T\subset
S$,  we choose one, say $T_1$, which is maximal with
respect to inclusion.  Assume, for a contradiction, that $|T_1|<m$.

We are given that $H$ has rank $\ge m$.  Hence among all sets $T$ such that
$T_1\subset T\subset S$ and $\langle T\rangle$ has rank $\ge m$ we may
choose one, say $T_2$, which is minimal with respect to
inclusion. Since $\langle T_1\rangle$ has rank $<m$, we have $T_2\ne
T_1$. We choose an element $x_0\in T_2-T_1$. The minimality of $T_2$
implies that the group $J=\langle T_2-\{x_0\}\rangle$ has rank
$<m$. Since $\langle T_2\rangle=\langle J\cup\{x_0\}\rangle$ has rank
$\ge m$, the rank of $J$ must be $m-1$ and $T_2$ must have rank
$m$. We fix a generating set $\{x_1,\ldots,x_{m-1}\}$ for $J$. Then
$\{x_0,\ldots,x_{m-1}\}$ is a generating set for $\langle T_2\rangle$.
But the rank-$m$ group $\langle T_2\rangle$ is free since $H$ is
$m$-free, and so $\{x_0,\ldots,x_{m-1}\}$ is a basis for $\langle
T_2\rangle$ (cf. \ref{kurosh stuff}). Hence $\langle T_2\rangle$ is a
free product $J\star\langle x_0\rangle$, where $\langle x_0\rangle$ is
infinite cyclic. In particular, the group $\langle T_1\cup
\{x_0\}\rangle$ is a free product $\langle T_1\rangle\star\langle
x_0\rangle$, and hence the set $ T_1\cup\{ x_0\}$ is independent. This
contradicts the maximality of $T_1$.  \EndProof

\Proposition\label{better margulis}
Let $k$ and $m$ be integers with $2\le m<k$, and let $M$ be a closed,
orientable hyperbolic $3$-manifold such that $\pi_1(M)$ is $k$-free.
Let $\mu$ be a Margulis number for $M$, let $\lambda$ be a positive
real number such that
$$\frac{m-1}{1+e^\lambda}+\frac1{1+e^{\mu}}\ge\frac12,$$
and let
$\Delta$ denote the  extrinsic diameter (\ref{extrinsic}) of
$\Mthick(\mu)$ in $M$. Suppose that
$$\frac m{1+e^{\lambda}}+\frac{k-m}{1+e^{2\Delta}}\ge\frac12.$$
Let $P$ be any point in $M$, and let $H$ denote the subgroup of
$\pi_1(M,P)$ generated by all elements that are represented by
loops of length $<\lambda$. Then $H$ has rank $<m$.
\EndProposition

\Proof 
Let us write $M=\HH^3/\Gamma$, where $\Gamma\le\Isom_+(\HH^3)$ is
discrete, cocompact and torsion-free.  We set $q=q_\Gamma$ (\ref{same
  to me}) and choose a point $z \in q^{-1}(P)$. We use the base point
$z\in\HH^3$ to identify $\pi_1(M,P)$ with $\Gamma$ (see \ref{same to
  me}). We may then regard $H$ as the subgroup generated by the set
$S$ consisting of all elements $\xi\in\Gamma$ such that
$\dist(\xi\cdot z,z)<\lambda$. The discreteness of $\Gamma$ implies
that $S$ is finite.  As $\Gamma$ is $m$-free, $H\le\Gamma$ is in
particular $m$-free; this will allow applications of Lemma
\ref{violet} to $H$.

We distinguish two cases, depending on whether
$ P\in\Mthick(\mu)$ or $  P  \in\Mthin(\mu)$.

First suppose that $P \in\Mthick(\mu)$. In this case we apply Lemma
\ref{violet}, taking $T_0=\emptyset$.  In order to show that $\rank
H<m$, it suffices to show that there is no independent set $T\subset
S$ such that $|T|=m$. Suppose that $T_0=\{\xi_1,\ldots,\xi_m\}$ is
such a set.  We let $d_i<\lambda$ denote the minimal length of a loop
based at $P$ and representing $\xi_i$.
 
Since $\xi_1,\ldots,\xi_m$ are independent and $\pi_1(M)$ is $k$-free,
we may apply \cite[Corollary \bigradcor]{cusp} to obtain a point
$Q\in \Mthick(\mu)$ such that $\rho=\dist_M(P,Q)$ satisfies
\Equation\label{lenny} \frac{k-m}{1+e^{2\rho}}+\sum_{i=1}^m
\frac{1}{1+e^{d_i}} \le\frac12.
\EndEquation
Hence $\rho\le\Delta$.  Since $\xi_i\in S$ we have
$d_i<\lambda$ for $i=1,\ldots,m$. It therefore follows from
(\ref{lenny}) that
$$\frac m{1+e^{\lambda}}+\frac{k-m}{1+e^{2\Delta}}<\frac12.$$ This
contradicts the hypothesis.

Now suppose that $  P  \in\Mthin(\mu)$. We fix an element $\eta\ne1$ of
$\Gamma$ such that $\dist(z,\eta\cdot z)<\mu$.  In this case we apply
Lemma \ref{violet}, letting $S\cup\{\eta\}$ play the role of $S$ in
Lemma \ref{violet}, and taking $T_0=\{\eta\}$.  In order to show that
$\rank H<m$, it suffices to show that there is no independent set $T$
with $\eta\in T\subset S\cup\{\eta\}$ such that $|T|=m$. Suppose that
$T_0=\{\xi_1,\ldots,\xi_m\}$ is such a set, with $\xi_1=\eta$ and
$\xi_2,\ldots,\xi_m\in S$. We write $d_i=\dist(z,\xi_i\cdot z)$ for
$i=1,\ldots,m$.  Since $\xi_1,\ldots,\xi_m$ are independent, it
follows from \cite[Theorem 6.1]{accs}, together with the main result
of \cite{agol} or \cite{cg}, that \Equation\label{charlie brown}
\sum_{i=1}^{m} \frac{1}{1+e^{d_i}}\leq\frac12.
\EndEquation

Since $\xi_1,\ldots,\xi_m\in S$, we have
$d_i\le\lambda$ for $i=2,\ldots,m$. Our choice of $\xi_1=\eta$
gives $d_1<\mu$. Hence (\ref{charlie brown}) gives
$$ \frac{m-1}{1+e^{\lambda}}+\frac1{1+e^{\mu}}<\frac12.$$
This contradicts the hypothesis.
\EndProof

\Corollary\label{biddy biddy bum bum to san fernando}
Let $k>2$ be an integer and let $M$ be a closed,
orientable hyperbolic $3$-manifold such that $\pi_1(M)$ is $k$-free.
Let $\mu$ be a Margulis number for $M$, let $\lambda$ be a positive
real number such that
$$\frac{1}{1+e^\lambda}+\frac1{1+e^{\mu}}\ge\frac12,$$ and let
$\Delta$ denote the extrinsic diameter (\ref{extrinsic}) of
$\Mthick(\mu)$ in $M$. Suppose that $$\frac
2{1+e^{\lambda}}+\frac{k-2}{1+e^{2\Delta}}\ge\frac12.$$ Then $\lambda$
is itself a Margulis number for $M$. \EndCorollary
\Proof
This is the case $m=2$ of Proposition \ref{better margulis}.
\EndProof

\Lemma\label{morris margolin} Suppose that $M$ is a closed, orientable
hyperbolic $3$-manifold such that $\pi_1(M)$ is $4$-free and $\vol
M\le\whatslightlybiggervol$. Then $\whatMargulis$ is a Margulis number
for $M$.
\EndLemma

\Proof 
Set $\lambda=\whatMargulis$ and $\mu=1.078$. By direct
computation we find that \Equation\label{sir rupert murgatroyd}
\frac{1}{e^\lambda+1}+\frac1{e^{\mu}+1}>\frac12.
\EndEquation

Since $\pi_1(M)$ is $2$-free, it follows from  \cite[Corollary
  4.2]{rankfour} that $\log3$ is a Margulis number for $M$, so
that $\mu<\log3$ is also a Margulis number for $M$.

Since $\pi_1(M)$ is $3$-free, it follows from \cite[Corollary
 \threefreevol]{last} that some point $P\in M$ is the center of a
 hyperbolic ball of radius $(\log5)/2$.

Let $\Delta$ denote the extrinsic diameter of the compact subset
$\Mthick(\mu)$ of $M$, and set
$$\rho=\max_{x\in\Mthick(\mu)}\dist(P,x).$$ The triangle
inequality implies that $\Delta\le2\rho$.

Choose a point $Q\in\Mthick(\mu)$ such that $\dist(P,Q)=\rho$. Set
$V=N(Q,\whatfunnyradius)$ and $W=N(P,h_3((\log5)/2)$. Since
$\whatfunnyradius<\mu/2$, the set $V$ is (intrinsically) isometric to
a hyperbolic ball of radius $\whatfunnyradius$, and hence $\vol
V=B(\whatfunnyradius)=\whatfunnyvolume\ldots$. On the other hand,
since $P$ is the center of a ball of radius $(\log5)/2$,  it
follows from (\ref{i lit another lunky})  that $\vol W\ge
B((\log5)/2)/d((\log5)/2)=3.087\ldots$.  Hence $\vol V+\vol
W>\whatslightlybiggervol\ge\vol M$, and therefore $V\cap
W\ne\emptyset$. The triangle inequality therefore implies that
$$\rho=\dist(P,Q)\le h_3\big(\frac{\log5}2\big)+\whatfunnyradius<\whatelse.$$
Hence $\Delta<\twicewhatelse$. We therefore have
\Equation\label{duke of plaza-toro}
\frac 2{1+e^{\lambda}}+\frac{2}{1+e^{2\Delta}}>
\frac 2{1+e^{\whatMargulis}}+\frac{2}{1+e^{\fricewhatelse}}>
\frac12.
\EndEquation

It follows from (\ref{sir rupert murgatroyd}) and (\ref{duke of
  plaza-toro}) that the hypotheses of Corollary \ref{biddy biddy
  bum bum to san fernando} hold with $k=4$, and with $\lambda$ and
$\mu$ defined as above. Hence $\lambda$ is a Margulis number for
$M$.
\EndProof

\section{Distant volume}\label{distant volume section}
\Number\label{over the hills and far away}
We set
$$\mu_0=\whatMargulis$$
and
$$h=h_3\bigg(\frac{\mu_0}2\bigg)=0.67\ldots,$$
where $h_3$ is the function defined in \ref{twain} and
is calculated using (\ref{who would}).

We shall define a function $\Vfar(D,\lambda)$ on the set $\calx$ that
was defined in Subsection \ref{chalk cliffs}.  The definition will use
the function $\rho_4$ that was defined in \ref{chalk cliffs}, the
function $B$ defined by (\ref{bees in my bonnet, pain in my heart}),
and the function $\Vbor$ given by  Definition \ref{bonus
  function def}.  We set $Z(D,\lambda) = \rho_4(D,\lambda) -
\lambda/2$ and define

$$\Vfar(D,\lambda)=
\begin{cases}
\displaystyle{\Vbor\big(\frac{\mu_0}2,h+\frac\lambda2\big)+
B\big(\min\big(\frac{\mu_0}{2},
     \frac12\big(Z(D,\lambda)-h\big)\big)\big)} &
  \text{if $Z(D,\lambda)> h$;}\\
\displaystyle{ B\big(\min\big(\frac{\mu_0}{2},Z(D,\lambda)\big)\big)} &
  \text{if $0<Z(D,\lambda)\le h$;}\\
0 &
  \text{if $Z(D,\lambda)\le 0$.}
\end{cases}
$$ 
\EndNumber

(Note that $\Vbor(\mu_0/2,h+\lambda/2)$ is defined
  for every $\lambda>0$, since $h+\lambda/2>h=h_3(\mu_0/2)$.)

\Lemma\label{vus far monotone}
The function $\Vfar$ is monotone decreasing in its first argument.
\EndLemma

\Proof
Suppose $D_1 < D_2$.  Set $\rho^{(i)}=\rho_4(D_i,\lambda)$ for
$i=1,2$. It is clear from the definition given in \ref{chalk cliffs}
that $\rho_4$ is monotone decreasing in its first argument, and hence
\Equation\label{a bridge too far}
\rho^{(1)}\ge\rho^{(2)}.
\EndEquation

It is clear from the definition of $\Vfar$ that $\Vfar(D_1,\lambda)$ and
$\Vfar(D_2,\lambda)$ are both non-negative. In the case where
$\rho^{(2)}\le\lambda/2$, we have $\Vfar(D_2,\lambda)=0$, and
hence
\Equation\label{eggheads of the world unite}
\Vfar(D_1,\lambda)\ge\Vfar(D_2,\lambda). 
\EndEquation
If $\rho^{(2)}>\lambda/2$, then by (\ref{a bridge too far}) we have
$\rho^{(1)}>\lambda/2$. 

If $\rho^{(1)}$ and
$\rho^{(2)}$ both lie in the interval $(\lambda/2, h+\lambda/2]$
then (\ref{eggheads
of the world unite}) follows from (\ref{a bridge too far}).
If $\rho^{(1)}$ and
$\rho^{(2)}$ 
both lie in the interval $(h+\lambda/2,\infty)$, then (\ref{eggheads
of the world unite}) follows from (\ref{a bridge too far}) and
Remark \ref{Vborincreases}. 
Finally, suppose that 
$\lambda/2<\rho^{(2)}\le h+\lambda/2$ and that
$\rho^{(1)} >h+\lambda/2$. Since the definition of $\Vbor$ immediately
implies that 
$\Vbor(R,\rho )\ge B(R)$ for any $R>0$ and any $\rho  > h_3(R)$, we
have 
$\Vfar(D_1,\lambda)\ge\Vbor(\mu_0/2,h+\lambda/2)\ge
B(\mu_0/2)\ge\Vfar(D_2,\lambda)$, so that
(\ref{eggheads of the world unite}) holds in this case as well.
\EndProof 

\Lemma\label{spanky and his gang}Let $M$ be a closed, orientable
hyperbolic $3$-manifold such that $\pi_1(M)$ is $4$-free.  Suppose
that $\mu_0$ is a Margulis number for $M$. Suppose that $P$ is a point of
$\XX_M$, and set $\lambda=\ff_M(P)$.   Then 
$(D_M(P),\lambda)\in\calx$, so that $\Vfar(D_M(P),\lambda)$ is
defined, and 
$$\vol (M-N(P,\lambda/2))\ge  \Vfar(D_M(P),\lambda).  $$
\EndLemma

\Proof
It follows from Lemma \ref{fracklog too} that $(D_M(P),\lambda)\in\calx$.

Let us set $D=D_M(P)$ and $\rho=\rho_4(D,\lambda)$.

The lemma is trivial if $\rho\le \lambda /2$, as we have
$\Vfar(D,\lambda )= 0$ in that case. We shall therefore assume that
$\lambda /2<\rho $. It then follows from Proposition
\ref{far beak} that  there is a
point $Y_\rho\in\Mthick(\mu_0)$ such that
\Equation\label{petrapavlovsk} \dist(P,Y_\rho)=\rho.  \EndEquation

Consider the case in which $\lambda/2<\rho\le h+\lambda/2$.  In this
case we set $r=\min(\mu_0/2,\rho-\lambda/2)>0$, so that according to
the definition in \ref{over the hills and far away} we have
$\Vfar(D)=B(r)$.  Since $r+ \lambda /2\le\rho\le\dist(P,Y_\rho)$, we
have $N(P, \lambda /2)\cap N(Y_\rho,r)=\emptyset$. Hence
\Equation\label{mammy yokum} \vol (M-N(P, \lambda /2))\ge \vol
N(Y_\rho,r).
\EndEquation
On the other hand, since $Y_\rho\in\Mthick(\mu_0)$ and $0<r\le\mu_0/2$, the set
$N(Y_\rho,r)$ is intrinsically isometric to a hyperbolic ball of radius
$r$, and so
\Equation\label{stummicks out of joint}
\vol N(Y_\rho,r)=B(r)=\Vfar(D,\lambda).
\EndEquation
In this case the conclusion of the lemma follows from (\ref{mammy
  yokum}) and (\ref{stummicks out of joint}).

Now consider the case in which $\rho> h+\lambda/2$. In
this case we set $\nu=\frac{1}{2}(\rho -(h+ (\lambda/2))) > 0$.
We then have 
\Equation\label{how do you spell sciatica}
\rho= h+\frac\lambda2+2\nu.
\EndEquation
 If we set
$t= \lambda /2+\nu$, we therefore have $\lambda/2<t<\rho$, and it
follows from
 Proposition
\ref{far beak} that there is a point
$Y_t\in\Mthick(\mu_0)$ such that 
\Equation\label{we'll always have secaucus}
\dist(Y_t,P)=t=\frac{\lambda}2+\nu .
\EndEquation
From (\ref{petrapavlovsk}), (\ref{how do you spell sciatica}),
(\ref{we'll always have secaucus}) and the triangle inequality it
follows that \Equation\label{your job was not to fade}
\dist(Y_t,Y_\rho)\ge h+\nu.
\EndEquation

From (\ref{we'll always have secaucus}), (\ref{your job was not to
  fade}) and the triangle inequality we deduce that
\Equation\label{they told me this would keep the rain out}
N(P,\lambda/2)\cap N(Y_t,\nu)=\emptyset = N(Y_\rho,h)\cap N(Y_t,\nu).
\EndEquation

On the other hand, since
$\rho > h+\lambda/2$, it follows from (\ref{petrapavlovsk}) that
$$\dist(Y_\rho,P)>h+\lambda/2,$$
so that the triangle inequality gives
\Equation\label{met a mutt today}
N(P,\lambda/2)\cap N(Y_\rho,h)=\emptyset.
\EndEquation

From (\ref{met a mutt today}) and (\ref{they told me this would keep
  the rain out}) we deduce that
\Equation\label{an out-of-the-way place}
\vol (M-N(P,\lambda/2))\ge\vol N(Y_t,\nu)+\vol N(Y_\rho,h).
\EndEquation

We now apply Proposition \ref{boroczky bonus}, taking $R=\mu_0/2$,
taking for $Y$ the $\mu_0$-thick point $Y_\rho$, and defining $P$ and
$\rho$ as
above. This gives $\vol N(Y,h_3(R)) \ge\Vbor(\mu_0/2,\rho)$. Since $\rho >
h+\lambda/2$, it follows from Remark \ref{Vborincreases} that
$\Vbor(\mu_0/2,\rho)>\Vbor(\mu_0/2,h+(\lambda/2))$. Hence
\Equation\label{cows were out of season}
\vol N(Y_\rho,h))>\Vbor(\frac{\mu_0}2,h+\frac\lambda2).
\EndEquation

Set $m=\min(\mu_0/2,\nu)$.
Since $m\le\mu_0/2$ and $Y_t\in\Mthick(\mu_0)$, the set $N(Y_t,m)$ is
intrinisically isometric to a hyperbolic ball of radius
$m$. Hence
\Equation\label{adieu, mr. pugh}
\vol N(Y_t,\nu) \ge B(m).
\EndEquation
From (\ref{an out-of-the-way place}), (\ref{cows were out of season})
and (\ref{adieu, mr. pugh}) it follows immediately that
$$\vol (M-N(P,\lambda/2)) >\Vbor(\frac{\mu_0}2,h+\frac\lambda2)+ B(m)
=\Vfar(D),$$
which gives the conclusion of the lemma in this case.
\EndProof

\section{The case where there is no short geodesic}\label{Numbers}

We set $\delta_0=0.58$ and $\lambda_0=\log7$.

\Lemma\label{your evil star's in the ascendant} Let $D$ be a number
with $\delta_0\le D \le.7$, and set $T_3=\Phi_3(\delta_0,D)$.  Then we
have $D <T_3<\lambda_0$, so that, in particular, $\Theta(D/2,\lambda_0/2)$ and
$\Theta(T_3/2,\lambda_0/2)$ are defined,
and
\Equation\label{i'm the defendant}
\cos\bigg(\Theta\bigg( \frac{D}{2},\frac{\lambda_0}{2} \bigg)-
\Theta\bigg(\frac{T_3}{2}, \frac{\lambda_0}{2} \bigg)\bigg)<
\frac{\cosh D \cosh T_3 - \cosh 2D} {\sinh D \sinh T_3}.
\EndEquation
\EndLemma

\Proof
Let $\phi$ denote the function defined on $[\delta_0,\infty)$
by $\phi(x)=\Phi_3(\delta_0,x)$.  Let $\theta$ denote the function
defined on $(0,\lambda_0/2)$ by $\theta(x)=\Theta(x,\lambda_0/2)$.
From the definition of $\Phi_n$ given in \ref{monster dread our
  damages} and the definition of $\Theta$ given in \ref{theta}, it is
clear that $\phi$ is monotone increasing on its domain and that
$\theta$ is monotone decreasing on its domain

If $D$ is a point of $[\delta_0,.7]$ then the monotonicity of
$\phi$ implies that 
$$
D \le .7 < 3\delta_0 = \phi(\delta_0) \le \phi(D) \le
 \phi(.7) = 1.766\ldots < \lambda_0 .
$$ 
Since $T_3 = \phi(D)$, this proves the first assertion.

To prove the second assertion we consider an arbitrary sub-interval
$[a,b]$ of $[\delta_0,.7]$.  From the definition given in
\ref{theta} it is clear that the function $\Theta$ is monotone
decreasing in its first argument, and hence for any $D\in[a,b]$ we
have
$$\theta( D/2)\ge\theta( b/2)\ge\theta(.7/2)=1.1\ldots,$$
and
$$\theta(\phi(D)/2)\le\theta(\phi( a)/2)\le\theta(\phi(\delta_0)/2)=0.36\ldots.$$
It follows that
$$\theta( D/2)-\theta(\phi(D)/2)\ge\theta( b/2)-\theta(\phi( a)/2)>0.$$
Since $\Theta$ takes values in $(0,\pi/2)$ we have
$\theta(D/2)-\theta(\phi(D)/2)<\pi/2$. Hence
\Equation\label{evil bert}
\cos(\theta( D/2)-\theta(\phi(D)/2))\le\cos(\theta( b/2)-\theta(\phi(a)/2))
\EndEquation
for any 
$D\in[a,b]$.

On the other hand, 
for any 
$D\in[a,b]$, using the monotonicity of the hyperbolic cotangent we
find that
\Equation\label{ahura mazda}
\begin{aligned}
\frac{\cosh D \cosh T_3 - \cosh 2D}
{\sinh D \sinh T_3}
&=&\coth D\coth\phi(D)-\frac{\cosh2D}{\sinh D\sinh\phi(D)}\cr
&\ge&\coth b\coth\phi(b)-\frac{\cosh2b}{\sinh a\sinh\phi(a)}.
\end{aligned}
\EndEquation

If for all $a$ and $b$ with $\delta_0\le a<b\le .7$ we set
$$\Delta(a,b)=\cos(\theta( b/2)-\theta(\phi(a)/2))-\bigg(\coth
b\coth\phi(b)-\frac{\cosh2b}{\sinh a\sinh\phi(a)}\bigg),$$
then it follows from (\ref{evil bert}) and (\ref{ahura mazda}) that,
for every interval $[a,b]\subset[\delta_0,.7]$ and every point
$D\in[a,b]$, we have
$$
\cos(\theta({D}/{2})-
\theta(\phi(D)/{2}))
-\bigg( \frac{\cosh D
  \cosh \phi(D) - \cosh 2D} {\sinh D \sinh \phi(D)}\bigg)
 \le\Delta(a,b).$$ 
In
particular, (\ref{i'm the defendant}) will hold for $D\in[a,b]$
provided that $\Delta(a,b)<0$.  But by direct computation we find that
$\Delta(.58,.63)$, $\Delta(.63,.67)$, $\Delta(.67,.68)$,
$\Delta(.68,.69)$, and $\Delta(.69,.7)$ are all  negative.
 This
completes the proof.
\EndProof

\Number\label{is this the court of the exchecquer} 
We define functions $\Vnearnought(D)$ and $\Vnear(D)$ on
$[\delta_0,\infty)$ as follows.  For $n = 2, 3$ we set
$T_n = T_n(D) = \Phi_n(\delta_0, D)$.
According to Lemma \ref{Phish sandwich} we have 
$T_2 = \Phi_2(\delta_0,D)\le2D$, so that
$(D, T_2)$ is contained in the domain of $\Psi$.
We may therefore define
$$\Vnearnought(D)=
B\bigg(\frac{\lambda_0}{2}\bigg)-2\sigma\bigg(\frac{\lambda_0}{2}, \frac{D}{2},
\frac{T_2}{2},\Psi(D,T_2)\bigg).$$ 
Finally, we define
$$
\Vnear(D)=\begin{cases}
 \Vnearnought(D) & \text{for $\delta_0\le D <  .7  $}\\
 \Vnearnought(D)-2\kappa(\lambda_0/2,T_3/2) &
 \text{for $D\ge  .7  $.}
\end{cases}
$$
\EndNumber

\Lemma\label{no bandanas} Let $M$ be a closed, orientable hyperbolic
 $3$-manifold which contains no closed geodesic of length $<\delta_0$.
 Suppose that $P$ is a point of $\XX_M$ such that $\ff_M(P)\ge\lambda_0$.
 (see \ref{eggs and fries}). Then $D_M(P)\ge\delta_0$, so that
  $\Vnear(D_M(P))$  is defined, and
$$\vol N(P,\lambda_0/2)\ge \Vnear(D_M(P)).$$
\EndLemma

\Proof We set $D=D_M(P)$.

We shall apply Lemma \ref{nested cubes}, taking $\lambda=\lambda_0$
and $\delta=\delta_0$. By direct computation we find that
$3\delta_0<\lambda_0<4\delta_0$.  The hypothesis that $M$ contains no
closed geodesic of length $<\delta_0$ implies that a generator $x$ of
$C_P$ (see \ref{eggs and fries}) has translation length $\ge\delta_0$.
Since we have assumed that $\ff_M(P)\ge\lambda_0$, the hypotheses needed
for the first assertion of Lemma \ref{nested cubes} are satisfied.

It follows from the first assertion of Lemma \ref{nested cubes} and
the definition of $\Vnearnought$ that
 $$\vol N\bigg(P,\frac{\lambda_0}{2}\bigg)\ge \Vnearnought
- 2\kappa\bigg(\frac{\lambda_0}{2},\frac{T_3}{2}\bigg).$$

In view of the definition of $\Vnear$, it follows that
$\vol N(P,\lambda_0/2)\ge\Vnear(D)$ if $D\ge  .7   $.

Now suppose that $D< 0.7 $.  In this case it follows from Lemma
\ref{your evil star's in the ascendant} that $\max(D,T_3)<\lambda_0$,
so that $\Theta(D/2,\lambda_0/2)$ and $\Theta(T_3/2,\lambda_0/2)$
are defined, and that
$$\cos\bigg(\Theta\bigg(
\frac{D}{2},\frac{\lambda_0}{2}\bigg)
-\Theta\bigg(\frac{T_3}{2},\frac{\lambda_0}{2}\bigg)\bigg) <
\frac{\cosh D \cosh T_3 - \cosh 2D}{\sinh D \sinh T_3}.
$$
It therefore follows from the second assertion of Lemma \ref{nested
cubes} and the definition of $\Vnearnought$ that 
$\vol N(P,\lambda_0/2) \ge \Vnearnought(D)$. 
Since $\Vnear(D)=\Vnearnought(D)$ in this case, the present lemma is
now proved in all cases.
\EndProof

\Lemma\label{if the parents ate the spinach}For any
$D\in[\delta_0,\infty)$ we have  $ (D,\lambda_0)\in\calx$, so that
$\Vfar(D,\lambda_0)$ is defined, and
\Equation\label{he wouldn't have caught pneumonia}
\Vnear(D)+\Vfar(D,\lambda_0)>\whatvol.
\EndEquation
\EndLemma

\Proof 
 If $D\ge\delta_0$ we have 
$$\frac1{1+e^D}+\frac1{1+e^{\lambda_0}}\le\frac1{1+e^{\delta_0}}+\frac18=.483\ldots<\frac12,$$
so that $ (D,\lambda_0)\in\calx$.

Before turning to the proof of (\ref{he wouldn't have caught
  pneumonia}), we shall summarize the monotonicity properties of
various functions that will be used in the proof.  It is clear from
the definition given in \ref{monster dread our damages} that the
function $\Phi_n$ is monotone increasing in its second argument for
each $n$.  We pointed out in \ref{omigosh} that $\kappa$ is monotone
decreasing in its second argument. We pointed out in
\ref{hyperbolic lore of cosines} that ${\Psi}$ is monotone increasing
in its first argument and monotone decreasing in its second.
According to Lemma \ref{vus far monotone}, the function $\Vfar$ is
monotone decreasing in its first argument.  Finally, since $\kappa$ is
monotone decreasing in its second argument, the function $x\mapsto
\kappa(\lambda_0/2,\Phi_3(\delta_0,x)/2)$ is monotone decreasing on
its domain.

We now turn to the proof of (\ref{he wouldn't have caught pneumonia}).
We first consider the case $D\ge\lambda_0$. In this case we have 
$\kappa(\lambda_0/2,D/2)=0$, and hence
$$\begin{aligned}\Vnearnought(D)&=B(\frac{\log
    7}{2})-2\sigma(\lambda_0/2, D/2,
  \Phi_2(\delta_0,D)/2,\Psi(D,\Phi_2(\delta_0,D)))\cr &=B(\frac{\log
    7}{2})-2\kappa( \lambda_0/2, \Phi_2(\delta_0,D)/2)\cr&\ge
  B(\frac{\lambda_0}{2})-2\kappa( \lambda_0/2, \Phi_2(\delta_0,\lambda_0)/2)\cr
  &=4.015\ldots\end{aligned}$$
Furthermore, we have
$$\Phi_3(\delta_0,D/2)\ge\Phi_3(\delta_0,\lambda_0/2)=2.307\ldots>\lambda_0.$$
Hence $\kappa(\lambda_0/2,\Phi_3(\delta_0,D)/2)=0$, and so
$\Vnear(D)=\Vnearnought(D)=4.015\ldots$.  Since $\Vfar(D,\lambda_0)=0$ in this
case, (\ref{he wouldn't have caught pneumonia}) holds.

For the rest of the proof we shall restrict
attention to the case $\delta_0\le D<\lambda_0$. 

Let us define a {\it useful interval} to be a half-open interval
$I=[a,b)\subset[\delta_0,\lambda_0)$ whose interior $(a,b)$ does not
contain $0.7$.

For any useful interval $[a,b)$, we have
$\Phi_2(\delta_0,a)\le2a\le2b$ according to \ref{is this the court of the
  exchecquer}. Hence $(b,\Phi_2(\delta_0,a))\in{\RR}^2$ lies in the domain of
${\Psi}$. In view of the monotonicity properties pointed out above,
for any $D\in[a,b)$ we have
$$\begin{aligned}\Vnearnought(D)&=B(\frac{\lambda_0}{2})-2\sigma(\lambda_0/2, D/2,
  \Phi_2(\delta_0,D)/2,\Psi(D,\Phi_2(\delta_0,D)))\cr &\ge
  B(\frac{\lambda_0}{2})-2\sigma(\lambda_0/2,  a/2,
  \Phi_2(\delta_0,D)/2,{\Psi}( b,\Phi_2(\delta_0,D)))\cr &\ge B(\frac{\lambda_0}{2})-2\sigma(\lambda_0/2,
  a/2, \Phi_2(\delta_0,a)/2,{\Psi}(b,\Phi_2(\delta_0,a))).\end{aligned}$$ Thus if for
  every useful interval $I=[a,b)$ we set
$$\mnearnought(I)=B(\frac{\lambda_0}{2})-2\sigma(\lambda_0/2, a/2,
  \Phi_2(\delta_0,a)/2,{\Psi}(b,\Phi_2(\delta_0,a))),$$
we have
\Equation\label{if i hadn't shot him}
\Vnearnought(D)\ge\mnearnought(I)\quad{\rm whenever}\quad D\in I.
\EndEquation

Let us associate a number $\mnear(I)$ to any useful interval $I=[a,b)$
by setting $\mnear(I)=\mnearnought(I)$ if $I\subset[\delta_0,  .7  )$ and
$\mnear(I)=\mnearnought(I)-\kappa(\lambda_0/2,\Phi_3(\delta_0,a)/2)$ if
$I\subset[  .7  ,\lambda_0)$. 

 It follows from (\ref{if i hadn't shot him}), the definition of
$\Vnear$ (see \ref{is this the court of the exchecquer}), and the
monotonicity of $x\mapsto \kappa(\lambda_0/2,\Phi_3(\delta_0,x)/2)$ 
that for any useful interval $I$ we have \Equation\label{we got 13 and
  they got 13} \Vnear(D)\ge\mnear(I)\quad{\rm for\ every}\ D\in I.
\EndEquation

Since $\Vfar$ is monotone decreasing in its first argument,
for any useful interval $I=[a,b)$ we have \Equation\label{i'm not a
cab driver} \Vfar(D,\lambda_0)\ge\Vfar(b,\lambda_0)\quad{\rm for\ every}\ D\in
I=[a,b).  \EndEquation

It follows from (\ref{we got 13 and they got 13}) and (\ref{i'm not a
  cab driver}) that for any useful interval $I=[a,b)$ we have
\Equation\label{i'm a coffee pot}
\Vnear(D)+\Vfar(D,\lambda_0)\ge\mnear(I)+\Vfar(b,\lambda_0)\quad{\rm for\ every}\ D\in I.
\EndEquation
Hence in order to complete the proof the lemma, it suffices to write
$[\delta_0,\lambda_0)$ as a union of a family $\cali$ of useful
intervals such that $\mnear(I)+\Vfar(b,\lambda_0)>\whatvol$ for every
$I=[a,b)\in\cali$.

The proof is completed by direct calculation.  A separate calculation
was done on each of the intervals $[\delta_0=0.58,0.598), [0.598,0.608),
[0.608, 0.618), [0.618,0.7), [0.7, \log 7=\lambda_0)$.  Each of these intervals
was further subdivided into 20 equal-sized subintervals $I_i =
[x_n,x_{n+1})$ for $n = 0, \ldots, 19$.  The subintervals $I_n$ are
useful, and it was verified that
$\mnear(I_n)+\Vfar(x_n,\lambda_0)>\whatvol$ for each $n$.

The minimum value computed in this way was $3.4409\ldots$, which arises for
the subinterval $[0.5971,0.598)$.
Note that the calculation of
$\mnear$ requires calculating a value of the function $\iota$.
See Subsection \ref{how we did iota} for an explanation of the methods
that were used to make this calculation.
\EndProof

\Lemma\label{dix kilometres a pied}
Let $M$ be a closed, orientable hyperbolic $3$-manifold such that
$\pi_1(M)$ is $4$-free.  Suppose  that $M$ contains no closed
geodesic of length $< \delta_0$. Then $\vol M>\whatvol$.
\EndLemma

\Proof 
We may assume that $\vol M\le\whatslightlybiggervol$, as otherwise
there is nothing to prove. Then by Lemma 
\ref{morris margolin}, $\mu_0=\whatMargulis$ is a Margulis number
for $M$.

Since $\pi_1(M)$ is $4$-free we may apply Corollary \ref{Key theorem}.
If alternative (i) of Corollary \ref{Key theorem} holds, i.e. if $M$
contains an embedded hyperbolic ball of radius $\lambda_0/2$, then we
have $\vol M\ge B(\lambda_0/2)=4.65\ldots$. (Recall that we
have set $\lambda_0=\log7$ in this section.)

Now suppose that alternative (ii) of Corollary \ref{Key theorem}
holds, i.e. that there is a point $P\in\XX_M$ with
$\ff_M(P)=\lambda_0$.  We set $D=D_M(P)$.  Since in particular
$\ff_M(P)\ge\lambda_0$, and since $M$ contains no closed geodesic of length
$<\delta_0$, it follows from the definition of $D_M$ (see
\ref{eggs and fries} )
that $D\ge\delta_0$, and it follows from Lemma \ref{no bandanas} that
\Equation\label{ca use les souliers} \vol N(P,\lambda_0/2)\ge
\Vnear(D).
\EndEquation
On the other hand, since $\ff_M(P)=\lambda_0$, and since $\mu_0$ is a
Margulis number for $M$, Lemma \ref{spanky and his gang} gives \Equation\label{the
  frog and the mouse} \vol (M-N(P,\lambda_0/2))\ge \Vfar(D,\lambda_0).
\EndEquation
From (\ref{ca use les souliers}) and (\ref{the frog and the mouse}) it
follows that
$$\vol M\ge \Vnear(D)
+ \Vfar(D,\lambda_0).$$
The conclusion of the present lemma now follows from Lemma \ref{if the
  parents ate the spinach}.
\EndProof

\section{The case where there is a short geodesic}\label{short section}

As in Section \ref{Numbers}, we set $\delta_0=.58$.
As in Section \ref{distant volume section}, we set $\mu_0=1.119$.

\Proposition\label{waiting for lefty} Let $M$ be a closed, orientable
hyperbolic $3$-manifold, and let $\mu$ be a Margulis number for $M$.
Suppose that $c$ is a closed geodesic in $M$ of length $l<\mu$, and
let $P$ be any point of $c$. Then $P\in\XX_M$ and $D_M(P)=l$.
Furthermore, we have
$$\vol N(P,\ff_M(P)/2)=B(\ff_M(P)/2)-2\kappa(\ff_M(P)/2,l/2).$$
\EndProposition

\Proof Let $C\le\pi_1(M,P)$ denote the image of $\pi_1(c,P)$ under the
inclusion homomorphism, and let $x$ denote a generator of
$C$. Since $c$ is a closed geodesic of length $l$, the cyclic subgroup
$C$ of $\pi_1(M,C)$ is maximal, and for any integer $n\ne0$, the
minimal length of a loop based at $P$ and representing $x^n$ is
$nl$. Now if $\alpha$ is a loop based at $P$ that represents an
element of $\pi_1(M,P)-C$, then $[\alpha]$ does not commute with
$x$; and since $l<\mu$, and $\mu$ is a Margulis number for $M$,
the length of $\alpha$ is $\ge\mu$, and in particular $>l$. It follows
that $C=C_M(P)$, that $P\in\XX_M$, and that $D_M(P)=l$.

The second assertion is an application of Proposition
\ref{c.d. rivington}. Let us write $M=\HH^3/\Gamma$ where
$\Gamma\le\isomplus(\HH^3)$ is discrete and torsion-free, and set
$q=q_\Gamma$. Set $\lambda=\ff_M(P)$. Let $\tP$ be a point of
$q^{-1}(P)$, and let us identify $\pi_1(M,P)$ with $\Gamma$ via the
isomorphism determined by the base point $\tP\in\HH^3$ (see \ref{same
to me}).  Since $c$ is a closed geodesic, the component $\tc$ of
$q^{-1}(c)$ containing $P$ is the axis of $C$. Let $\rho_{+1}$ and
$\rho_{-1}$ be the closed rays emanating from $\tP$ and contained in
$\tc$. We may suppose them to be labeled in such a way that $x^{\epsilon
n}\cdot\tP\in\rho_\epsilon$ for every $n>0$ and for $\epsilon=\pm1$.

For each integer $n\ne0$, we have
$\dist(\tP,x^n\cdot\tP)=|n|l$. For $\epsilon=\pm1$, let
$\zeta_\epsilon$ denote the point of intersection of $S(\lambda/2,\tP)$
with $\rho_\epsilon$.
According to Proposition \ref{c.d. rivington} we have
\Equation\label{For those who like that sort of thing}
\vol(N(P,\lambda/2)) =B(\lambda/2)-\vol(\bigcup_{n>0,\epsilon=\pm1}
K(\lambda/2,\tP,\zeta_\epsilon,nl/2)).
\EndEquation

In the notation of \ref{basic caps}, the plane
$\Pi(\tP,\zeta_\epsilon,nl)$ is orthogonal to $\tc$ for each $n>0$ and
for $\epsilon=\pm1$. Hence the half-spaces $H(z_0,\zeta_{+1},l/2)$ and
$H(z_0,\zeta_{-1},l/2)$ are disjoint, and
$$ H(z_0,\zeta_\epsilon,nl/2) \subset H(z_0,\zeta_\epsilon,l/2)$$
for each
$n>1$ and for $\epsilon=\pm1$. It follows that 
\Equation\label{that is the sort of thing that they like}
K(\lambda/2,\tP,\zeta_{+1},l/2)\cap
K(\lambda/2,\tP,\zeta_{-1},l/2)=\emptyset,
\EndEquation
and that
\Equation\label{ain't it}
K(\lambda/2,\tP,\zeta_\epsilon,nl/2))\subset
K(\lambda/2,\tP,\zeta_\epsilon,l/2))
\EndEquation
for each $n>0$ and for $\epsilon=\pm1$.
From (\ref{For those who like that sort of thing}), 
(\ref{that is the sort of thing that they like}) and (\ref{ain't
  it}), we find that 
$$\begin{aligned}\vol(N(P,\lambda/2)) &=B(\lambda/2)-(\vol(
K(\lambda/2,\tP,\zeta_{+1}),l/2)+\vol(
K(\lambda/2,\tP,\zeta_{-1},l/2))\\
&=B(\lambda/2)-2\kappa(\lambda/2,l/2).
\end{aligned}
$$
\EndProof

We define a function $W$ with domain $\calx\subset\RR^2$ (see
\ref{chalk cliffs}) by
$$W(l,\lambda)=\Vfar(l,\lambda)+B(\lambda/2)-2\kappa(\lambda/2,l/2).$$

\Lemma\label{sea beach} Let $k$ be an integer $>2$, let $M$ be a
closed, orientable hyperbolic $3$-manifold such that $\pi_1(M)$ is
$k$-free, and suppose that $\mu_0$ is a Margulis number for $M$.
Suppose that $c$ is a closed geodesic in $M$ of length $l<\mu_0$, and
let $P$ be any point of $c$. Then $P\in\XX_M$ and
$(l,\ff_M(P))\in\calx$. Furthermore, we have
$$\vol M\ge W(l,\ff_M(P)).$$
\EndLemma

\Proof
Set $\lambda=\ff_M(P)$.  According to Lemma \ref{waiting for lefty},
we have $P\in\XX_M$, $D_M(P)=l$, and
$$\vol N(P,\lambda/2 )=B(\lambda/2 )-2\kappa(\lambda /2,l/2).$$

According to Lemma \ref{spanky and his gang}, we have
$(D_M(P),\lambda)\in\calx$, so that $\Vfar(D_M(P),\lambda)$ is
defined, and
$$\vol (M-N(P,\lambda/2))\ge  \Vfar(D_M(P),\lambda)=\Vfar(l,\lambda).  $$
Hence 
$$\begin{aligned}\vol M&\ge \vol N(P,\lambda/2 )+\vol(M-N(P,\lambda/2))\\
&\ge B(\lambda/2 )-2\kappa(\lambda /2,l/2)+\Vfar(l,\lambda)\\
&=W(l,\lambda).\end{aligned}$$
\EndProof

\Lemma\label{brighton beach} For any  $(l,\lambda)\in\calx$ with
$.003\le l\le\delta_0$, we have $W(l,\lambda)>\whatvol$.  
\EndLemma

\Proof
Let us define $\rho(l) = \frac12\log((e^l+3)/(e^l-1))$, so that
$$\frac{1}{1+e^{l}}+\frac{1}{1+e^{2\rho(l)}} = \frac12.$$ By
definition the set $\calx$ consists of all points in $\RR^2$ such that
$\lambda \ge 2\rho(l)$; or equivalently, of all points
of the
form $(l, 2\rho(l) + y)$ with $l>0$ and $y\ge0$.

We define $W^*(l,y) = W(l, 2\rho(l) + y)$ for $l>0$ and $y\ge0$.  It suffices to show
that for any $l\in[.003,\delta_0]$ and any $y>0$ we have 
$W^*(l,y) > \whatvol$.  We observe that the function $\rho$ is
monotone decreasing for $l > 0$.

For $(l, \lambda) \in \calx$ we define
$$V_N (l, \lambda) = B(\lambda/2) - 2\kappa(\lambda/2, l/2).
\leqno{\hskip 10pt 13.3.1}$$
Thus the function $W$ defined before Lemma 13.2 can be expressed as
$W = V_{\rm far} + V_N $.

The function $V_N$ is increasing in its first argument because $\kappa$ is
(weakly) decreasing in its second argument while $B$ is increasing.  We
claim that $V_N$ is also increasing in its second argument. To prove this,
suppose that numbers $\lambda_1$, $\lambda_2$ and $l$ are given, with $l>0$
and $0<\lambda_1\le\lambda_2$.  Fix a line $L\subset\HH^3$ and a point
$P\in L$.  Consider the the two planes that are perpendicular to $L$ and
have distance $l/2$ from $P$, and let $Y$ denote the closed neighborhood of
$P$ bounded by these two planes.  For $i=1,2$, let $N_i$ denote the
ball of radius $\lambda_i/2$ centered at $P$.  Then, by the definition of
$\kappa$, we have $\vol(Y\cap N_i)=V_N(l,\lambda_i)$. Since
$\lambda_1\le\lambda_2$, we have $N_1\subset N_2$ and therefore
$Y\cap N_1\subset Y\cap N_2$; hence
$V_N (l,\lambda_1)\le V_N (l,\lambda_2)$, and our claim is proved.

Next, for $l > 0$ and $y \ge 0$ we set $V_N^*(l,y) = V_N(l,2\rho(l)+y)$ and
$V_{\rm far}^* (l,y) = V_{\rm far}(l,2\rho(l)+y)$, so that
$W^* = V_{\rm far} ^*+ V_N^ *$.

We claim that for every $l\in[0.003,\delta_0]$ and every $y>0.5$, we have
$V_N^ *(l, y) > 3.44$. Before
proving the claim we observe that if we are given an arbitrary $y>0.5$, an
arbitrary closed interval $[a,b]\subset(0,\infty)$, and an arbitrary
$l\in[a,b]$ then, since $\rho$ is decreasing, and since we have shown that
$V_N$ is increasing in each of its arguments, we have
$V_N^ *(l, y)= V_N(l,2\rho(l)+y) \ge V_N(a,2\rho(b)+0.5)$.

Hence, since $\delta_0=0.58$, to prove our claim it suffices to
exhibit $[0.003,0.58]$ as the union of a finite collection of closed
subintervals in such a way that, for each interval $[a,b]$ in the
collection, we have $V_N(a,2\rho(b)+0.5) > 3.44$.  We take our
collection to consist of the subintervals $[x_{i-1},x_i]$ for
$i=1,\ldots,24$, where $(x_0,\ldots,x_{24})$ is the $25$-tuple whose
terms are $0.003$, $0.004$, $0.006$, $0.0085$, $0.012$, $0.017$,
$0.023$, $0.03$, $0.04$, $0.05$, $0.06$, $0.08$, $0.11$, $0.15$,
$0.2$, $0.24$, $0.3$, $0.36$, $0.42$, $0.44$, $0.47$, $0.5$, $0.535$,
$0.56$, and $0.58$. For each interval $[a,b]$ in this collection, the
formula for the function $\kappa$ given by Proposition A.3 can be used
to directly compute $V_N(a,2\rho(b)+0.5)$. The smallest value produced
by these computations is $3.441\ldots$, which arises for $a = 0.004$
and $b = 0.006$.  Thus our claim is established.

In particular, for every $l\in[0.003,\delta_0]$ and every $y>.5$, we have
shown that $W^*(l,y) > 3.44$.  To complete the proof of the lemma, it
therefore suffices to show that for any $(l, y)$ in the rectangle
$R: = [0.003, \delta_0] \times [0, 0.5]$ we have $W^*(l,y) > 3.44$.

We denote by $\cals$ the set of all subrectangles of $R$ of the form
$[l_0,l_1]\times[y_0,y_1]$ with $.003<l_0<l_1\le\delta_0$ and $0\le
y_0<y_1\le.05$.

We define a continuous function $\chi$ on $R$ by
  $$\chi(l,y)=\rho_4(l,2\rho(l)+y)-\bigg(
  h+\frac{2\rho(l)+y}{2}\bigg).$$
It follows from 
the definition of $\Vfar$ given in Subsection \ref{over the hills
  and far away} that we have 
\Equation\label{your homework}
\Vfar^*(l,y)\ge B\big(\max\big(0,\big(\min\big(\frac{\mu_0}{2},\rho_4(l,{2\rho(l)+y})-
  \dfrac{{2\rho(l)+y}}{2}\big)\big)\big)\big)
\EndEquation
for every $(l,y)\in R$, with equality when 
$\chi(l,y)\le0$. Furthermore, when $\chi(l,y)\ge0$ we have
\Equation\label{my homework}
\Vfar^*(l,y)=\Vbor\big(\frac{\mu_0}2,h+\frac{2\rho(l)+y}2\big)+B\big(\min\big(\frac{\mu_0}{2},
     \frac12\big(\rho_4(l,{2\rho(l)+y})-\big(h+
     \dfrac{{2\rho(l)+y}}{2}\big)\big)\big) \big).
\EndEquation

The function $\rho$ is monotone decreasing, and the function
$\rho_4$ is  decreasing in each of its arguments. Hence for each
$S=[l_0,l_1]\times[y_0,y_1]\in\cals$, the function
$\chi|S$ is bounded below by
$$\chi_S\dot=\rho_4(l_1,2\rho(l_0)+y_1)-\bigg(
  h+\frac{2\rho(l_0)+y_1}{2}\bigg).$$

Suppose that $\cals=[l_0,l_1]\times[y_0,y_1]$ is a rectangle such that
$\chi_S>0$.  Then in particular $\chi$
takes only positive values on $S$, and hence
(\ref{my
  homework}) holds for every $(l,y)\in S$.
Since $\rho$ is a monotone decreasing function, 
$\rho_4$ is  decreasing in each of its arguments, 
and $B$ is increasing,
and since $\Vbor$ is increasing in its second argument according to
Remark \ref{Vborincreases}, the quantity
\Equation\label{a magic pencil}
V^S_+\dot=\Vbor\big(\frac{\mu_0}2,h+\frac{2\rho(l_1)+y_0}2\big)+B\big(\min\big(\frac{\mu_0}{2},
     \frac12\big(\rho_4(l_1,{2\rho(l_0)+y_1})-\big(h+
     \dfrac{{2\rho(l_0)+y_1}}{2}\big)\big)\big)\big )
\EndEquation
is a lower bound for $\Vfar^*|S$ for any $S=[l_0,l_1]\times[y_0,y_1]\in\cals$ with
$\chi_S>0$.

On the other hand, if
$S=[l_0,l_1]\times[y_0,y_1]$ is an {\it arbitrary} rectangle in
$\cals$ then (\ref{your
  homework}) holds for every $(l,y)\in S$. Hence
\Equation\label{that will grade}
V^S_-\dot=
B\big(\max\big(0,\big(\min\big(\frac{\mu_0}{2},\rho_4(l_1,{2\rho(l_0)+y_1})-
  \dfrac{{2\rho(l_0)+y_1}}{2}\big)\big)\big) \big) 
\EndEquation
is a lower bound for $\Vfar^*|S$, for any
  $S=[l_0,l_1]\times[y_0,y_1]\in\cals$.

We have observed that $W^* = \Vfar^* + V_N^*$, and that
$V^*_N(l,y)$ is decreasing in $l$
and increasing in $y$. It follows that 
$V_N^S\doteq V_N^*(l_1,y_0)$ is a lower bound for $V_N^*|S$ for any
$S=[l_0,l_1]\times[y_0,y_1]\in\cals$.
Hence 
$V_N^S+V_+^S$ 
is a lower bound for $W^*|S$ 
for every 
$S=[l_0,l_1]\times[y_0,y_1]\in\cals$
with $\chi_s>0$, and 
$V_N^S+V_-^S$ 
is a lower bound for $W^*|S$ 
for every 
$S=[l_0,l_1]\times[y_0,y_1]\in\cals$.

Hence in order to complete the proof of the
lemma, it suffices to 
specify subsets
$\cals_+$ and $\cals_0$ of $\cals$ with
 $\cals_+\subset\cals_0$,
such that
\begin{itemize}
\item $R$ is the union of the rectangles in $\cals_0$;
\item $\chi_S > 0$ for every $S\in \cals_+$; 
\item $V_N^S+V_+^S>3.44$ for every $S\in\cals_+$; and
\item $V_N^S+V_-^S>3.44$ for every $S\in\cals_0-\cals_+$.
\end{itemize}

To define $\cals_0$ we begin with the rectangles $R_1 =
[0.003,0.103]\times[0,0.5]$, $R_2 = [0.1,0.5]\times[0,0.5]$, and $R_3
= [0.5,0.58]\times[0,0.5]$.  We subdivide each of the $R_i$ into
equal-sized subrectangles, where the subrectangles are chosen to form a
$40\times 100$ grid on $R_1$; a $50\times 100$ grid on $R_2$ and an
$80\times 100$ grid on $R_3$.  We define $\cals_0$ to be the union of
the rectangles in these three grids.

To specify $\cals_+\subset\cals_0$ we compute $\chi_S$ numerically for
each $S\in\cals_0$ and take $\cals_+$ to consist of those rectangles
for which the computed value of $\chi_S$ exceeds $0.1$.   Numerical
  computation of $\chi_s$ involves evaluation of elementary functions
  and arithmetic.  Since the round-off errors in the arithmetic and
  the errors inherent in the standard approximations of elementary
  functions by rational functions combine to give a margin of error
much less than $0.1$, we indeed have $\chi_S > 0$ for every $S\in
\cals_+$.

We computed $W_S$ numerically for each $S\in\cals_0$, using one of the
formulas \ref{a magic pencil} or \ref{that will grade} and using
Proposition \ref{formula for kappa} to compute the value of the
function $\kappa$ which appears in formula \ref{V sub N}.  The
minimum value of $W_S$ obtained in this manner is
$3.4511\ldots$, which arises for $S=[0.579,0.58]\times[0.145,0.15]$.
This shows that $W^*$ is bounded below by $\whatvol$ on $R$, as
required.
\EndProof

\Lemma\label{v'la le vitrier}Let $M$ be a closed, orientable
hyperbolic $3$-manifold such that $\pi_1(M)$ is $4$-free.  Suppose
that $M$ contains a closed geodesic of length $< \delta_0$. Then $\vol
M>\whatvol$.
\EndLemma

\Proof
 We may assume that $\vol M\le\whatslightlybiggervol$, as otherwise
there is nothing to prove. Then by Lemma 
\ref{morris margolin}, $\mu_0=\whatMargulis$ is a Margulis number
for $M$. 

Let $c$ be a closed geodesic in $M$ of length $l<\delta_0$. Let $P$ be
any point of $c$.
Since $\delta_0<\mu_0$, we may apply Lemma \ref{sea beach} to deduce that $P\in\XX_M$, that
$(l,\ff_M(P))\in\calx$, and that
$$\vol M\ge W(l,\ff_M(P)).$$
If $l\ge.003$, it follows from Lemma \ref{brighton beach} that
$W(l,\ff_M(P))>\whatvol$, and hence that $\vol M>\whatvol$.

There remains the case in which $l<.003$. In this case, let $T$ denote
the maximal embedded tube about $c$. It follows from \cite[Corollary
10.5]{accs} that \Equation\label{margaret and the giraffe} \vol T\ge
V(.003), \EndEquation where $V$ is the function defined in Section 10
of \cite{accs}. Computing $V(.003)$ from the definition given in
\cite{accs} we find that \Equation\label{who blew up da owl}
V(.003)=3.1345\ldots \EndEquation On the other hand, it follows from a
result of Przeworski's \cite[Corollary 4.4]{prez} on the density of
cylinder packings that \Equation\label{wintergreen for prezident} \vol
T < 0.91\vol M.  \EndEquation From (\ref{margaret and the giraffe}),
(\ref{who blew up da owl}) and (\ref{wintergreen for prezident}) it
follows that
$$\vol M>\frac{3.134}{.91}>\whatvol.$$
\EndProof

\Proof[Proof of Theorem \ref{the big one}] The theorem is an
  immediate consequence of Lemmas \ref{dix kilometres a pied} and
  \ref{v'la le vitrier}.
\EndProof

The proof of Theorem \ref{sweet land of flub-a-dub} will involve
combining Theorem \ref{the big one} with the results of
\cite{lastplusone}. We refer the reader to \cite[Section
\nonfibroidsection]{lastplusone} for the definition of a fibroid. As
in \cite{lastplusone}, we will use a result due to Agol, Storm, and Thurston
from \cite{ast}. The information from \cite{ast} that we need is summarized
in Theorem \firstAST\ of \cite{last}, which states that if $M$ is a
closed orientable hyperbolic $3$-manifold containing a connected
incompressible closed surface which is not a fibroid, then
$\vol(M)>3.66$.

 \Proof[Proof of Theorem \ref{sweet land of flub-a-dub}] Assume
that $\dim_{\Z_2}H_1(M;\Z_2)\ge8$. Then according to \cite[Proposition
\ifnotwhynot]{lastplusone}, either $\pi_1(M)$ is $4$-free, or $M$
contains a closed incompressible surface of genus at most $3$ which is
not a fibroid. If $\pi_1(M)$ is $4$-free, it follows from Theorem
\ref{the big one} that $\vol(M)>\whatvol$. If $M$ contains a closed
incompressible surface which is not a fibroid, it follows from
\cite[Theorem \firstAST]{last} that $\vol(M)>3.66$. In either case the
hypothesis is contradicted.  \EndProof

\section{Appendix: Computations with caps}\label{cap comps}

In this section we will describe the methods used for numerical
computation of particular values of the functions
$\iota(R,\D,\D',\alpha)$ and $\kappa(R,D)$ which were needed for the
proofs of Lemmas \ref{if the parents ate the spinach} and
\ref{brighton beach}. The main results are Propositions \ref {formula
  for kappa}, \ref {special intersection}, \ref {complements} and \ref
{formula for iota}. In Subsection \ref{how we did iota} we will show
how to combine these results to calculate $\iota$ and
$\kappa$ for any values of the arguments.

The following well known special case of the distance formula in
the upper half-space model of $\HH^2$ will be needed.

\Lemma\label{ez distance} Let $\gamma$ be a geodesic in the hyperbolic
plane modeled in the upper half-plane by a semicircle $s$
with center $X = (x,0)$ and radius $\rho$.  Let $P = (x,\rho)$.
Suppose that $Q$ is a point of $s$ such that the arc of $s$ from $P$
to $Q$ subtends an angle $\theta$ in the Euclidean plane.  Then $\cosh
d_h(P,Q) = \sec\theta$.
\EndLemma

\Proof After applying a hyperbolic isometry which fixes
$\infty$, we may assume that $\rho = 1$, $x = 0$, $P = (0,1)$ and $Q
= (\sin\theta, \cos\theta)$.  The arc of the unit circle from $P$ to
$Q$ is a hyperbolic geodesic arc of length
$$\int_0^\theta \frac{d\theta}{\cos\theta} =
\log|\sec\theta+\tan\theta| .$$
Thus, since $0 <\theta \le \pi/2$, we have
$$
\begin{aligned}
\cosh d_h(P,Q)
&= \frac12\bigg(\sec\theta+\tan\theta+\frac1{\sec\theta+\tan\theta}\bigg)\\
&=\sec\theta.
\end{aligned}
$$
\EndProof

We will also need the following fact from Euclidean geometry, of which
the proof is an easy exercise:

\Lemma\label{harvest moon} Let $S_1$ and $S_2$ be two spheres in
$\EE^3$ such that $S_1\cap S_2$ is a circle $\cali $. Let $r_i$ denote
the radius of $S_i$, and let $D$ denote the distance between the
centers of $S_1$ and $S_2$. Then the radius $r$ of $\cali $ satisfies
$$r^2=r_1^2-\frac{(r_1^2+D^2-r_2^2)^2}{4D^2}.$$
Furthermore, the distance between the centers of $S_1$ and $\cali$ is
$(r_1^2+D^2-r_2^2)/(2D)$.
\NoProof
\EndLemma

\Proposition\label{formula for kappa}
Let $R$ and $\D$ be positive numbers with $0 < \D < R$.
Set 
 $$\epsilon = e^{-\D}\sqrt{1 - \frac{\cosh^2\D}{\cosh^2 R}}$$
 and
$$L(r) =  \cosh R - \sqrt{\sinh^2R  - r^2}.$$
Then
$$\kappa(R,\D) = \pi(e^R\cosh R - \frac{\cosh R}{L(\epsilon)} - \log L(\epsilon) - R
+ \frac12\log(e^{-2\D} - \epsilon^2) + \D).$$
\EndProposition

\Proof
We use the notation of Subsection \ref{basic caps}. By definition $\kappa(R,\D)$
is the volume of $K(R, Z ,\zeta,\D)$, where $ Z $ and $\zeta$ are points in
$\HH^3$ separated by a distance $R$. We may identify $\HH^3$
conformally with the upper half-space
$\UU^3=\RR^2\times(0,\infty)\subset\RR^3$ in such a way that
$ Z =(0,0,1)$  and $\zeta=(0,0,e^{-R})$. Then $\bar N( Z ,R)$ is
identified with the Euclidean ball $B$ of radius $\sinh R$ centered at
$(0,0,\cosh R)$.
The half-space $H( Z ,\zeta,\D)$ is identified with the intersection
of $\UU^3$ with the Euclidean ball $B'$ of radius $e^{-\D}$ centered at
$(0,0,0)$. 

The boundaries of $B$ and $B'$ intersect in a circle $\cali $. According to
Lemma \ref{harvest moon}, the square of the radius of $\cali $ is equal to
$$e^{-2\D}-\frac{
(e^{-2\D}+\cosh^2R-\sinh^2R)^2}
{4\cosh^2R},$$
which implies that the radius of $\cali $ is the quantity $\epsilon$
defined in the statement of the present proposition.

Since $ Z $ and $\zeta$ lie on the vertical axis $\{(0,0\}\times\RR$,
the circle $\cali $ lies in a horizontal plane and has center on the
vertical axis.  It follows that the vertical projection
$p:(x,y,t)\mapsto (x,y)$ maps $K(R, Z ,\zeta,\D)
  =\overline{N}( Z ,R)\cap H( Z ,\zeta,\D)$ onto a disk $\Delta$ of
  radius $\epsilon$ about $(0,0)$, and that for every $P\in\Delta$ the
  set $p^{-1}(P)$ is a line segment whose lower and upper endpoints
  lie, respectively, in the lower hemisphere of $\partial B$ and the
  upper hemisphere of $\partial B'$. If $r$ denotes the distance from
  $P$ to $(0,0)$, the definitions of the balls $B$ and $B'$ imply that
  the vertical coordinates of these endpoints are respectively equal
  to $L(r)$ and $U(r)$, where $L(r)$ is defined as in the statement of
  the proposition, and $U(r) = \sqrt{e^{-2\D} - r^2}$.

Using cylindrical coordinates in $\RR^3$, we therefore find that
$$\vol K(R,Z,\zeta,\D)=
\int_0^{2\pi}\int_0^\epsilon\int_{L(r)}^{U(r)}\frac{r}{t^3}dt\,dr\,d\theta .$$

We have
$$
\begin{aligned}
\int_0^{2\pi}\int_0^\epsilon\int_{L(r)}^{U(r)}\frac{r}{t^3}dt\,dr\,d\theta
&= \pi\int_0^\epsilon\left(\frac{r}{L(r)^2} - \frac{r}{U(r)^2}\right)dr\\
&= \pi\int_0^\epsilon\frac{r\,dr}{L(r)^2} -
\pi\int_0^\epsilon\frac{r\,dr}{e^{-2x} - r^2}\\
&= \pi\int_0^\epsilon\frac{r\,dr}{L(r)^2} + \pi(\frac12\log(e^{-2x} -
\epsilon^2) + x) .
\end{aligned}
$$

To complete the proof we observe that
$r\,dr = (\cosh R  - L(r))d\,L(r)$, and $L(0) = e^{-R}$.  Thus
we have
$$
\begin{aligned}
\int_0^\epsilon\frac{r\,dr}{L(r)^2} &=
\int_{L(0)}^{L(\epsilon)}\frac{\cosh R - u}{u^2}du\\
&= -\cosh R\left(\frac{1}{L(\epsilon)} - \frac{1}{L(0)}\right) -
\log L(\epsilon) + \log L(0) \\
&= e^R\cosh R - \frac{\cosh R}{L(\epsilon)} - \log L(\epsilon) - R,
\end{aligned}
$$
which completes the proof. 
\EndProof

The following
result, Proposition \ref{special intersection}, gives an 
integral formula that can be used to compute $\iota(R,\D,0,\alpha)$,
when $\alpha\in(\pi/2, \pi)$. The subsequent results of this section,
Propositions \ref{complements} and \ref{formula for iota}, will be
used to reduce the
general calculation of $\iota(R,\D,\D',\alpha)$ to the special case
covered by Proposition \ref{special intersection}.

\Proposition\label{special intersection}
Let $R$ and $\D$ be positive numbers with $0 < \D < R$, and let
$\alpha\in(\pi/2, \pi)$ be given.  Set $\phi = \alpha - \pi/2\in (0,\pi/2)$. 
\begin{enumerate}
\item \label{fat frat}If we define quantities 
$\radtwo$ and $c$ by
$\radtwo=(\sinh w+\cosh w\cos\phi)^{-1}$ and
$c = v\sin\phi\cosh w$,
then $${v^2 - \frac{(v^2+c^2+1)^2}{4(c^2+\cosh^2R)}}\ge0.$$
\item\label{with a dismal headache} Let us define quantities $\mu$, $\rho$ and $ m $ by 
$$\mu = \sqrt{v^2 - \frac{(v^2+c^2+1)^2}{4(c^2+\cosh^2R)}},\quad
\rho = \frac{\cosh R}{\sqrt{c^2+\cosh^2R}},\quad\text{and\ }
m  = c - \frac{c(v^2+c^2+1)}{2(c^2+\cosh^2R)}.$$
If $ m  > \rho\mu$ then $\iota(R,\D,0,\alpha) = 0$.
\item\label{and repose is tabooed}
Suppose that $ m  \le \rho\mu$, and define a quantity $\theta_0$ and a
function $\epsilon(\theta)$ by
$$\theta_0 = \arctan
(\frac{\sqrt{\rho^2\mu^2- m ^2}}{ m \rho}
)\text{\qquad and \qquad}
\epsilon(\theta) =\frac{ \rho\mu}{\sqrt{\cos^2\theta+ \rho^2\sin^2\theta} }.$$
Then $m\sec\theta\le\epsilon(\theta)$ for every
$\theta\in[-\theta_0,\theta_0]$. Furthermore, if we set
$$V=\{\theta,r\in
[-\pi,\pi]\times\RR\ \colon -\theta_0\le\theta\le\theta_0
\text{\  and \ } m \sec\theta\le r\le\epsilon(\theta)\},$$
then for every $(\theta,r)\in V$ we have 
$$
\sinh^2R  - (r\cos\theta -  m )^2 - r^2\sin^2\theta\ge0$$
and
$$\radtwo^2 - (r\cos\theta + c -  m )^2 - r^2\sin^2\theta\ge0.$$
Finally, if we define 
functions
$U(r,\theta)$ and $L(r,\theta)$ on $V$ by
$$U(r,\theta) =
 \sqrt{\radtwo^2 - (r\cos\theta + c - m)^2 - r^2\sin^2\theta}$$ and
$$L(r,\theta) =
 \cosh R - \sqrt{\sinh^2R  - (r\cos\theta - m)^2 - r^2\sin^2\theta},$$
then
$L(r,\theta)\le
U(r,\theta)$ for each $(r,\theta)\in V$, and
we have
\Equation\label{we bombed in new haven}
 \iota(R,\D,0,{\alpha}) =
\int_0^{\theta_0}\int_{m\sec\theta}^{\epsilon(\theta)}
\bigg(\frac{r}{L(r,\theta)^2}-\frac{r}{U(r,\theta)^2}\bigg)\,dr d\theta.
\EndEquation
\end{enumerate}
\EndProposition

\Proof
In this proof it will be understood that $v$, $c$, $\mu$, $\rho$, $m$,
$\theta_0$, $\epsilon(\theta)$, $V$, $U(r,\theta)$ and  $L(r,\theta)$
are defined as in the statement of the proposition. Many of these objects are defined
only subject to certain assertions in the statement, and will not be
mentioned until after these assertions have been proved.
 
We use the notation of Subsection \ref{basic caps}. By definition we
have
$$ \iota(R,\D,0,{\alpha})=\vol(K(R, Z ,\zeta_1,0)\cap K(R, Z ,\zeta_2,\D),$$
where $ Z $ is a point in $\HH^3$ and $\zeta_1$ and $\zeta_2$ are
points in $S(R, Z )$ separated by a spherical distance $\alpha$.  For
$i=1,2$ we set $\eta_i=\eta_{\zeta_i}$, and we let $\ell_i$ denote the
line containing $\eta_i$. We set $\Pi_1=\Pi( Z ,\zeta_1,0)$,
$H_1=H( Z ,\zeta_1,0)$, $\Pi_2=\Pi( Z ,\zeta_2,\D)$,
$H_2=H( Z ,\zeta_2,\D)$.  $K_1=K(R, Z ,\zeta_1,0)$, and
$K_2=K(R, Z ,\zeta_2,\D)$.

Since $\alpha>\pi/2$, we have $\zeta_1\notin H_2$. In particular,
\Equation\label{not!}
K_1\not\subset K_2.
\EndEquation

Since the rays $\eta_1$ and $\eta_2$ form an angle $\alpha\in(\pi/2,\pi)$ at
$ Z $, there
is a unique plane $U\subset\HH^3$ containing $\ell_1$ and
$\ell_2$. In particular $U$ is perpendicular to $\Pi_1$, and the line
$U\cap\Pi_1$ is perpendicular to $\ell_1$. Hence there is a
unique ray $\tau\subset U\cap\Pi_1$ with origin $ Z $ which forms an
angle $\phi=\alpha-\pi/2$ with $\eta_2$.

Since $0<x<R$, the sphere $S$ and the plane $\Pi_2$ meet in a circle
$\cali$.

Now let $\UU^3$ denote the upper half-space in $\RR^3$. We denote the
coordinates in $\RR^3$ by $x$, $y$ and $t$, so that $\UU^3$ is defined
by $t>0$.
We identify $\HH^3$ conformally with 
$\UU^3$ in such a way that
$ Z =(0,0,1)$; 
$\Pi_1$, $H_1$ and $U$ are the subsets of $\UU^3$ defined,
respectively, by $x=0$, $x\ge0$ and $y=0$; and $\tau$ is defined by
$x=y=0$ and $0<t\le1$.

If $X$ and $Y$ are points of $\bar\UU^3\subset\RR^3$ we shall denote
the Euclidean distance between $X$ and $Y$ by $|XY|$. When
$X,Y\in\UU^3$ we shall write $d_h(X,Y)$ for the hyperbolic distance
between $X$ and $Y$.

We will treat the $xz$-plane, which contains $U$, as a
Cartesian plane with coordinates $x$ and $z$; in particular each
Euclidean line in this plane  has a well-defined {\it
  slope}\/.

Since $ Z =(0,0,1)$, the ball $\bar N( Z ,R)$ is
identified with the Euclidean ball $\calb$ of radius $\sinh R$ centered at
$A\dot=(0,0,\cosh R)$. We set $S=\partial\calb$.

Since $\phi>0$ we have $\ell_2\not=U\cap\Pi_1$, so that the tangent vector to $\ell_2$ at
$ Z =(0,0,1)$ is a non-vertical vector in $\RR^3$. Hence the
hyperbolic line
$\ell_2$ is identified with
a semicircle in the half-plane $U$, whose center $B$ lies on the
$x$-axis. We write $B=(b,0,0)$.  Since
$\eta_2$ forms an angle $\alpha=\phi+(\pi/2)$ with $\eta_1$ (which
has tangent vector $(1,0,0)$) and 
forms an angle $\phi$ with $\tau$ (which
has tangent vector $(0,0,-1)$),
the semicircle $\ell_2$ has positive
slope at $Z$, and hence $b>0$.

If $E$ denotes the point of intersection of $\ell_2$ with $\Pi_2$,
then since $E\in\zeta_2$, the $x$-coordinate of $E$ is negative. It
follows that the tangent line to the semicircle $\ell_2$ has
positive slope at $E$ and hence that the
tangent line at $E$ to the semicircle $\Pi_2\cap U$,
which meets $\ell_2$ orthogonally  at the point $E$, has
negative slope.   This implies that 
$\Pi_2 $ is identified with a hemisphere whose center $C$ lies on the
$x$-axis and has a negative $x$-coordinate.

\begin{figure}[h]
\begin{picture}(0,0)%
\includegraphics{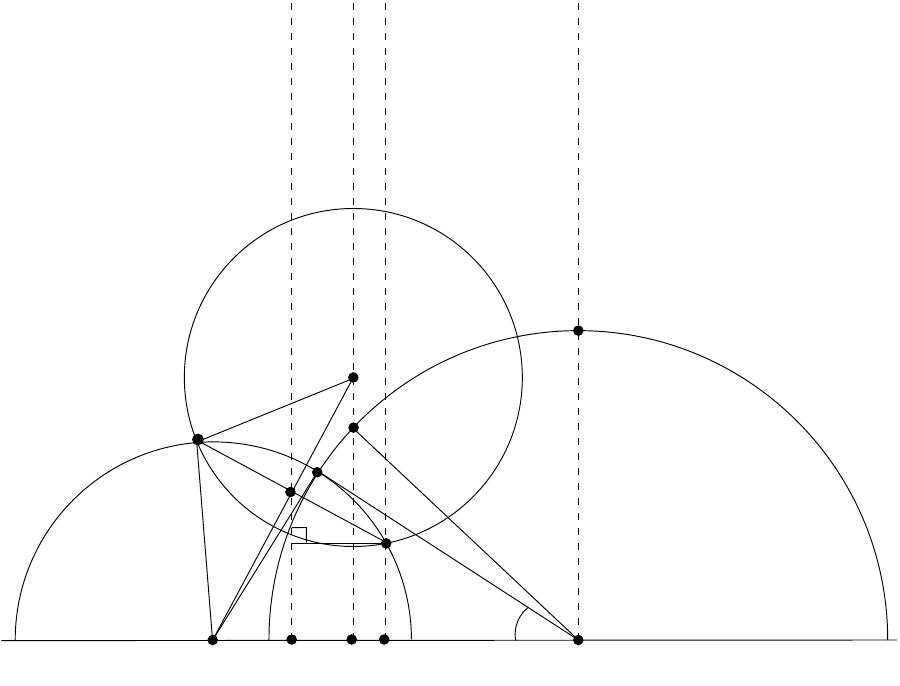}%
\end{picture}%
\setlength{\unitlength}{3947sp}%
\begingroup\makeatletter\ifx\SetFigFont\undefined%
\gdef\SetFigFont#1#2#3#4#5{%
  \reset@font\fontsize{#1}{#2pt}%
  \fontfamily{#3}\fontseries{#4}\fontshape{#5}%
  \selectfont}%
\fi\endgroup%
\begin{picture}(7189,5476)(24,-5525)
\put(4576,-5461){\makebox(0,0)[lb]{\smash{{\SetFigFont{12}{14.4}{\familydefault}{\mddefault}{\updefault}{\color[rgb]{0,0,0}$B$}%
}}}}
\put(1351,-3511){\makebox(0,0)[lb]{\smash{{\SetFigFont{12}{14.4}{\familydefault}{\mddefault}{\updefault}{\color[rgb]{0,0,0}$F$}%
}}}}
\put(3911,-4381){\makebox(0,0)[lb]{\smash{{\SetFigFont{12}{14.4}{\familydefault}{\mddefault}{\updefault}{\color[rgb]{0,0,0}$u$}%
}}}}
\put(4801,-2611){\makebox(0,0)[lb]{\smash{{\SetFigFont{12}{14.4}{\familydefault}{\mddefault}{\updefault}{\color[rgb]{0,0,0}$P$}%
}}}}
\put(3001,-3586){\makebox(0,0)[lb]{\smash{{\SetFigFont{12}{14.4}{\familydefault}{\mddefault}{\updefault}{\color[rgb]{0,0,0}$Z$}%
}}}}
\put(2626,-2986){\makebox(0,0)[lb]{\smash{{\SetFigFont{12}{14.4}{\familydefault}{\mddefault}{\updefault}{\color[rgb]{0,0,0}$A$}%
}}}}
\put(6601,-3511){\makebox(0,0)[lb]{\smash{{\SetFigFont{12}{14.4}{\familydefault}{\mddefault}{\updefault}{\color[rgb]{0,0,0}$\ell_2$}%
}}}}
\put(2701,-3886){\makebox(0,0)[lb]{\smash{{\SetFigFont{12}{14.4}{\familydefault}{\mddefault}{\updefault}{\color[rgb]{0,0,0}$E$}%
}}}}
\put(3994,-5050){\makebox(0,0)[lb]{\smash{{\SetFigFont{12}{14.4}{\familydefault}{\mddefault}{\updefault}{\color[rgb]{0,0,0}$\beta$}%
}}}}
\put(3301,-4561){\makebox(0,0)[lb]{\smash{{\SetFigFont{12}{14.4}{\familydefault}{\mddefault}{\updefault}{\color[rgb]{0,0,0}$G$}%
}}}}
\put(1426,-4486){\makebox(0,0)[lb]{\smash{{\SetFigFont{12}{14.4}{\familydefault}{\mddefault}{\updefault}{\color[rgb]{0,0,0}$v$}%
}}}}
\put(2101,-4111){\makebox(0,0)[lb]{\smash{{\SetFigFont{12}{14.4}{\familydefault}{\mddefault}{\updefault}{\color[rgb]{0,0,0}$M$}%
}}}}
\put(1621,-5446){\makebox(0,0)[lb]{\smash{{\SetFigFont{12}{14.4}{\familydefault}{\mddefault}{\updefault}{\color[rgb]{0,0,0}$C$}%
}}}}
\put(2151,-5461){\makebox(0,0)[lb]{\smash{{\SetFigFont{12}{14.4}{\familydefault}{\mddefault}{\updefault}{\color[rgb]{0,0,0}$M_0$}%
}}}}
\put(2776,-5461){\makebox(0,0)[lb]{\smash{{\SetFigFont{12}{14.4}{\familydefault}{\mddefault}{\updefault}{\color[rgb]{0,0,0}$Q$}%
}}}}
\put(151,-3886){\makebox(0,0)[lb]{\smash{{\SetFigFont{12}{14.4}{\familydefault}{\mddefault}{\updefault}{\color[rgb]{0,0,0}$\Pi_2\cap U$}%
}}}}
\put(1426,-2011){\makebox(0,0)[lb]{\smash{{\SetFigFont{12}{14.4}{\familydefault}{\mddefault}{\updefault}{\color[rgb]{0,0,0}$S\cap U$}%
}}}}
\end{picture}%
\end{figure}

We let
$\radone$ denote the radius of the semicircle $\ell_2$, and we set $\beta=\angle CBE$.

We have $\angle CBZ=\phi$.  Since the Euclidean distance from $(0,0)$
to $Z$ is $1$, we have $\radone= \csc\phi$ and $b= \cot\phi$.

Let $P$ denote the point of intersection of $\ell_2$ with the vertical
ray in $U$ originating at $B$.  By Lemma \ref{ez distance}, we have
$\cosh(d_h(P,Z)) = \csc\phi$.  Since $d_h(P,E) = \D + d_h(P,Z)$,
Lemma \ref{ez distance} implies that $ \csc\beta= \cosh(\D +
\arccosh\radone) $. Hence
$$\cot\beta=\sinh(\D + \arccosh\radone)=\cosh w\sqrt{u^2-1}+u\sinh w
=\cosh w\cot\phi+\sinh w\csc\phi$$
and
$$\cos\beta= \tanh(\D + \arccosh\radone)=\frac{\tanh w
+\sqrt{1-u^{-2}}}{1+\sqrt{1-u^{-2}}\tanh w}=\frac{\tanh w+\cos\phi}{1+\tanh w\cos\phi}.$$
Since $BEC$ is a Euclidean right triangle, we have 
$$|CE|=\radone\tan\beta =\frac{\csc\phi}{\cosh w\cot\phi+\sinh w\csc\phi}=\radtwo.$$
Thus
\Equation\label{whistlin' dixie}
  {\rm radius}(\Pi_2) =\radtwo  .
\EndEquation

The difference of the $x$-coordinates of $B$ and $C$ is $|BC| =
\radone\sec\beta$. Hence  the
$x$-coordinate of $C$ is
$$b-\radone\sec\beta=\cot\phi-(\csc\phi)\,\frac
{1+\tanh w\cos\phi}
{\tanh
  w+\cos\phi}=-\frac{\sin\phi}{\tanh w+\cos\phi}=-c,$$
i.e. 
\Equation\label{nijinsky who's doing a rhumba}
C=(-c,0,0).
\EndEquation

We have seen that the center of the Euclidean sphere $S$ is
$(0,0,\cosh R)$. It therefore follows from (\ref{nijinsky who's doing
  a rhumba}) that the Euclidean distance between the centers of the
sphere $S$ and the hemisphere $\Pi_2$ is $\sqrt{c^2+\cosh^2 R}$. We
have also seen that $S$ has radius $\sinh R$, while $\Pi_2$ has radius
$\radtwo $ by (\ref{whistlin' dixie}).  In the notation of Lemma
\ref{harvest moon}, taking $S_2=S$ and taking $S_1$ to be the sphere
containing $\Pi_1$, we have $r_1=v$, $r_2=\sinh R$ and $D=
\sqrt{c^2+\cosh^2R}$. The first assertion of the lemma gives
$${\rm
    radius}(\cali)^2=
v^2-\frac{(v^2+c^2+1)^2}{4(c^2+\cosh^2R)}.
$$
This shows that $v^2-(v^2+c^2+1)^2/(4(c^2+\cosh^2R))>0$, which is
assertion (\ref{fat frat}) of the proposition. It also shows that if
we define $\mu$ as in the  statement of the proposition then
\Equation\label{grover whelan unveilin'}
{\rm radius}(\cali)=\mu.
\EndEquation
The second assertion of Lemma \ref{harvest moon} gives the distance
between the center $C$ of $\Pi_2$ and the center of $\cali$, which we
denote by $ M $:  

\Equation\label{kankakee or paree}
|C M | = 
 \frac{\radtwo^2+c^2+ 1}{2\sqrt{c^2 + \cosh^2R}}.
\EndEquation

Let $p:\RR^3\to\RR^2$ denote the projection $(x,y,t)\mapsto
(x,y)$. We have $p( M )=(- m _0,0)$ for some $ m _0\in\RR$.
By similar triangles we have
$(c- m _0)/|C M |=c/|AC|$, so
$$ m _0 = c(1 - \frac{|C M |}{|AC|}) =c - \frac{c(v^2+c^2+1)}{2(c^2+\cosh^2(R))}= m .$$
Hence:
\Equation\label{or washington crossing the delaware}
p( M )=(- m ,0).
\EndEquation

Since the sphere $S$ and the hemisphere $\Pi_2$ are centered on the
closed Euclidean half-plane $ \bar U$, their circle of intersection
$\cali$ meets $U$ in two points $F$ and $G$, which are diametrically
opposite points of $\cali$. Hence $ M $  is the midpoint of
the segment $\overline{FG}$.  We write $p(F)=-f$ and $p(G)=g$, and we label
$F$ and $G$ in such a way that $-f<g$. In view of
(\ref{or washington crossing the delaware}) we have
$ m =(f-g)/2$. 

Note that $\cali$ is contained in the subset of $\UU^3$ defined by
$-f\le x\le g$. Note also that $AC$, the line joining the centers of
$S$ and $\Pi_1$, meets the diameter $\overline{FG}$ of $\cali$
perpendicularly at the center $M$.  Furthermore, since $A\notin\Pi_2$
and $M\in\Pi_2$, the point $M$ lies on the segment $\overline{AC}$.
Let $\delta$ denote the diameter of $\cali$ which is perpendicular to
the line $FG$. Then $\delta$ is invariant under reflection about $U$
and is therefore contained in a horizontal plane.

We now derive an expression for $g+ m $. We define points $ M _0,Q,Z\in U$
by setting $ M _0=(- m ,0,0)$ and
$Q=(0,0,0)$, and defining $Z$ to be the intersection of the vertical
line through $ M $ with the horizontal line in $U$ through $G$. As the
triangles $ACQ$ and $ M C M _0$ are similar, we have
$\angle ACQ= \angle  M C M _0=(\pi/2)-\angle C M  M _0=\angle C M G-\angle
C M  M _0=\angle G M  M _0$.
Hence the right triangles $ACQ$ and $G M  M _0$ are similar, and so
$$|AC|(g+ m ) = |G M |\cosh R.$$
It follows from (\ref{grover whelan unveilin'}) that $|G M |=\mu$, and
we have $|AC|^2=\cosh^2\mu+c^2$. Hence if $\rho$ is defined as in the
statement of the proposition, we have
\Equation\label{captain spaulding exploring the amazon}
g+ m =\rho\mu.
\EndEquation

We can now prove assertion (\ref{with a dismal headache}) of the
proposition. If $ m  > \rho\mu$, then by (\ref{captain spaulding
  exploring the amazon}) we have $g<0$. As we have observed that
$\cali$ is contained in the subset of $\UU^3$ defined by $-f\le x\le
g$, it follows that $\cali$ is contained in the open half-space
$\UU^3-H_1$ of $\UU^3$. Hence the set $\Pi_2\cap B$, which is the
hyperbolic convex hull of $\cali$, is also contained in
$\UU^3-H_1$, and is in particular disjoint from $K_1$. But $H_2\cap B$ is the frontier of $K_2$ relative to $B$,
and so $K_2$ either contains $K_1$ or is disjoint from it. The former
alternative is ruled out by (\ref{not!}); hence $K_1\cap
K_2=\emptyset$, and so
$\iota(R,\D,0,\alpha) = 0$. This is assertion (\ref{with a dismal headache}) of the
proposition. 

We now turn to the proof of Assertion (\ref{and repose is tabooed}) of
the proposition.
Assume that $ m  \le \rho\mu$. Then by (\ref{captain spaulding
  exploring the amazon}) we have $g\ge0$. 
This means that $G\in H_1$,
so that 
\Equation\label{lassie come-home}
\cali\cap H_1\ne\emptyset.
\EndEquation

It is clear from the definition of $m$ that $m>0$. Since $g\ge0$ and
$m=(f-g)/2$, it follows that $f>0$.

Let $T\subset\RR^3$ denote the
half-space $z\le\cosh R$, whose boundary plane contains $A$. Since $M$
lies on the segment $\overline{AC}$ we have $M\in T$.

Since $0<\alpha<\pi/2$, the definition of $c$ implies that $c>0$. In
view of (\ref{nijinsky who's doing a rhumba}), it follows that the
line $AC\subset U$ has positive slope, and that the segment
$\overline{AC}$ is disjoint from $H_1$. In particular we have $M\notin
H_1$. Since the line $FG$ is perpendicular to $AC$, it has negative
slope in $U$, and hence the ray originating at $M$ and passing through
$G$ is contained in $T$. In particular $FG\cap H_1\subset T$.

Let $W$ denote the plane containing $\cali$. Since $W$ contains $FG$
and is perpendicular to $U$, it now follows that $W\cap H_1\subset T$.

In particular, if $S_-$ and $S_+$ denote the
lower and upper hemisphere of $S$, we have 
\Equation\label{albel}
\cali\cap H_1\subset S_-.
\EndEquation

The projection $p$ maps $\cali$ onto an ellipse $\mathcal E$ in
$\RR^2$. Since $\delta$ and $FG$ are mutually perpendicular diameters
of $\cali$, and $\delta$ is contained in a horizontal plane, $p$ maps
$\delta$ and $FG$ onto the major and minor axes of $\cale$,
respectively. In particular the minor axis of $E$ is contained in the
$x$-axis; the length of the semi-major axis of $\cale$ is the radius of $\cali$,
which by (\ref{grover whelan unveilin'}) is equal to $\mu$; and the
semi-minor axis of $\cale$ has length $g+ m $, which by (\ref{captain
  spaulding exploring the amazon}) is equal to $\rho\mu$. The center
of $\cale$ is $p( M )=(- m ,0)$.

 We denote by $\hat\cale$ the compact set bounded
by $\cale$. Since $f>0$  and $g\ge0$, the line
  $x=0$ meets $\hat\cale$ in a possibly degenerate line segment
  $\nu$. We let $V_0$ denote the intersection of $\hat\cale$ with the
half-plane $x\ge0$.

The hemisphere $S_\pm$ is the graph of the function
$L_0^\pm(x,y)=\cosh R\pm\sqrt{\sinh^2R-x^2-y^2}$ on the disk
$\cald\subset\RR^2$ which has radius $\sinh R$ and is centered at
$(0,0)$. Likewise, it follows from (\ref{whistlin' dixie}) and
(\ref{nijinsky who's doing a rhumba}), $p(\Pi_2)$ is the graph of the
function $U_0(x,y)=\sqrt{v^2-(x+c)^2-y^2)}$ on the disk
$\cald'\subset\RR^2$ which has radius $v$ and is centered at $(-c,0)$.

We have $\cale=p(\cali)= p(S\cap\Pi_2)\subset p(S)\cap
p(\Pi_2)=\cald\cap\cald'$.  Since $\cald\cap\cald'$ is convex, we
have $\hat\cale\subset \cald\cap\cald'$.

The functions $\psi^\pm= U_0-L_0^\pm$ have domain $\cald\cap\cald'$.
Let $\calr$ denote the right half-plane in $\RR^2$, defined by
$x\ge0$.  Set $\calw=\cald\cap\cald'\cap\calr$.  It follows from
(\ref{albel}) that the function $\psi^-$ is identically zero on $\calr
\cap \cale=\calw\cap \cale$, and is non-zero on $\calw\setminus\cale$,
while the function $\psi^+$ is non-zero on $\calw$. Since $-\psi^-$ is
clearly a convex function, and since $\calw \cap \cale$ is the
frontier relative to $\calw$ of the convex set $\calr \cap\hat
\cale=\calw \cap\hat \cale$, the function $\psi^-$ must be
non-negative on $\calw \cap \hat\cale$ and negative on
$\calw\setminus\hat\cale$.  The function $\psi^+$ is non-zero on the
connected domain $\calw$, and is bounded above by $\psi^-$. But
$\psi^-$ vanishes on the subset $\calw\cap\cale$ of $\calw$, and this
subset is non-empty by (\ref{lassie come-home}). Hence $\psi^+$ is
negative-valued on $\calw$.

It follows that we have
\Equation\label{queen of tattoo}
L_0^-(x,y)\le U_0(x,y)<L_0^+(x,y)\text { when  }(x,y)\in\calw \cap\hat \cale,
\EndEquation
and
\Equation\label{mazurka in jazz}
L_0^-(x,y)> U_0(x,y)\text{ when }(x,y\in \calw\setminus\hat\cale.
\EndEquation
We have $K_1\cap K_2=B\cap H_1\cap H_2$, and $p(B\cap
H_2)\subset\cald\cap\cald'$. Set $q=p|(K_1\cap K_2):K_1\cap K_2\to\cald\cap\cald'$. It follows from 
(\ref{mazurka in jazz}) that $q^{-1}(x,y)=\emptyset$ when $(x,y)\in 
\calw\setminus\hat\cale$; and it follows from
(\ref{queen of tattoo}) that when $(x,y)\in\calr \cap\hat \cale=\calw
\cap\hat \cale$
we have $L_0^-(x,y)\le U_0(x,y)$, and $q^{-1}(x,y)$ is a
vertical line segment whose endpoints have $t$-coordinates 
$L_0^-(x,y)$ and $ U_0(x,y)$. Since the element of hyperbolic volume
on $\UU^3$ is ${dt\,dA}/{t^3}$, where $dA$ is the Euclidean area
element on $\RR^2$, we now find
\Equation\label{a view of niagara}
\begin{aligned} \iota(R,\D,0,\alpha)&=\vol(K_1\cap K_2)\\
&=\int_{\calr \cap\hat \cale}\int_
{L_0^-(x,y)}^{U_0(x,y)}\frac{dt\,dA}{t^3}\\
&=\int_{\calr \cap\hat \cale}\bigg(\frac1{2 L_0^-(x,y)^2}-
\frac1{2 
U_0(x,y)^2}\bigg)\,dA.
\end{aligned}
\EndEquation

We shall complete the proof by re-interpreting
the facts proved above in terms of polar coordinates $(r,\theta)$ in
$\RR^2$, taking the origin for the polar coordinates to be the point
$(- m ,0)$. Since $\cale$ is centered at $(- m ,0)$, has its minor
axis contained in the $x$-axis, and has semi-major axis of length $\mu$ and
semi-minor axis of length $g+ m $, it is defined in these coordinates by the
equation $r = \epsilon(\theta)$.  The set $\hat\cale$ is defined by
$r \le \epsilon(\theta)$. The half-plane $\calr$ is
defined by $ r \ge m \sec\theta$. The endpoints of the segment
  $\nu$ are the intersections of $\cale$ with the $y$-axis, whose
  polar equation is $r = m \sec\theta$. In view of the definitions of
  $\epsilon$ and $\theta_0$ it is clear that these intersection points
  are $(\pm\theta_0,m\sec\theta_0)$. Since $\nu\subset\hat\cale$, it
  follows that we have $ m
    \sec\theta\le\epsilon(\theta)$ whenever
    $-\theta_0\le\theta\le\theta_0$.

It now follows that  $V$ is the set of polar coordinate-pairs of points in $\calr
\cap\hat \cale$.  Since $\calr \cap\hat
\cale\subset\hat\cale\subset\cald\cap\cald'$, we have $ {\sinh^2R -
  (r\cos\theta -  m )^2 - r^2\sin^2\theta}\ge0$ and 
$${\radtwo^2
  - (r\cos\theta + c -  m )^2 - r^2\sin^2\theta}\ge0$$
 for every
$(\theta,r)\in V$.

The transition to polar coordinates
transforms $U_0(x,y)$ and $L_0^-(x,y)$ to the  functions $U(r,\theta)$
and $L(r,\theta)$. It follows that
$L(r,\theta)\le U(r,\theta)$ for all $(r,\theta)\in V$. The area element is given by $r\,dr\,d\theta$, and 
(\ref{a view of niagara}) becomes
$$
\begin{aligned}
\iota(R,\D,0,\alpha)&=
\int_{-\theta_0}^{\theta_0}
\int_{ m \sec\theta}^{\epsilon(\theta)}
\bigg(\frac{r}{2 L(r,\theta)^2}-
\frac{{r}}{2 
U(r,\theta)^2}\bigg)\,dr\,d\theta\\
&=
\int_0^{\theta_0}
\int_{ m \sec\theta}^{\epsilon(\theta)}
\bigg(\frac{r}{ L(r,\theta)^2}-
\frac{r}{
U(r,\theta)^2}\bigg)\,dr\,d\theta.\\
\end{aligned}
$$
\EndProof

\Proposition\label{complements}
If $\alpha \in [0,\pi/2]$  and $\D\le0$ then
$$\iota(R,\D,0,\alpha) + \iota(R,\D,0,\pi-\alpha)
= \kappa(R,\D) .$$
\EndProposition

\Proof 
Let $ Z \in\HH^3$ be given, set $B=\overline{ N( Z ,R)}$, and let $\zeta_1$ and $\zeta_2$ be points in $S\dot=S(R, Z )$
separated by a spherical distance $\alpha$. Let $\zeta_1^*$ denote the
antipode of $\zeta_1$ on $S$, so that the spherical distance between
$\zeta_1^*$ and $\zeta_2$ is $\pi-\alpha$. 
Then $K_1\dot=K(R, Z ,\zeta_1,0)$ and $K_1^*\dot=K(R, Z ,\zeta_1^*,0)$ are the two
half-balls bounded by the plane $H_1\dot=H( Z ,\zeta_1,0)$, so that
$K_1\cup K_1^*=B$ and $K_1\cap K_1^*=H_1$. In particular, setting
$K_2=K(R, Z ,\zeta_2,\D)$, we have
$(K_1\cap K_2)\cup (K_1^*\cap K_2)=K_2$ and
$(K_1\cap K_2)\cap (K_1^*\cap K_2)=K_2\cup H_1$. Hence
$$\iota(R,\D,0,\alpha) + \iota(R,\D,0,\pi-\alpha)=
\vol(K_1\cap K_2)+\vol (K_1^*\cap K_2)=\vol K_2
= \kappa(R,\D) .$$
\EndProof

\Proposition\label{formula for iota} Let $Z$ be a point in $\HH^3$,
let $R>0$ be a real number,
and let $\alpha$, $\D_1$
and $\D_2$ be real numbers with  $0\le\alpha\le\pi$ and $0 < \D_1 \le \D_2 < R$. 
Define quantities $\Psi_1$ and $\Psi_2$ 
by $\Psi_i = \arccos(\tanh \D_i/\tanh R)$ (so  that $0<\Psi_2
  \le \Psi_1<\pi/2$).
\begin{enumerate}
\item If $\alpha \le \Psi_1 - \Psi_2$ then 
  $\iota(R,\D_1,\D_2,\alpha) =\kappa(R,\D_2)$.
\item If $\alpha > \Psi_1 + \Psi_2$ then 
$\iota(R,\D_1,\D_2,\alpha) =0$.
\item If $\Psi_1 - \Psi_2 < \alpha \le \Psi_1 + \Psi_2$ then
there exists a unique pair
$(\alpha_1,\alpha_2)$ of numbers with $-\pi/2 \le \alpha_i\le
  \pi/2$ such that 
\Equation\label{truck talk}
\alpha_1+\alpha_2 = \alpha
\EndEquation
 and 
\Equation\label{i hold that on the seas}
 {\tanh \D_1}{\cos \alpha_2} = {\tanh \D_2}{\cos \alpha_1} .
\EndEquation
Moreover, 
\Equation\label{that great street}
\iota(R,\D_1,\D_2,\alpha) =
 \iota(R,\D_1,0,\alpha_1+\pi/2) + \iota(R,\D_2,0,\alpha_2+\pi/2).
\EndEquation
\end{enumerate}
\EndProposition

\Proof 
Let $Z$ be a point in $\HH^3$, and let $\zeta_1$ and $\zeta_2$ be
points of $S\doteq S(R,Z)$ separated by a spherical distance $\alpha$.
We set $\calb=N(R,Z)$.  For $i=1,2$, set 
$\eta_i=\eta_{\zeta_i}$,
$\Pi_i=\Pi( Z ,\zeta_i,\D_i)$, $H_i=H( Z ,\zeta_i,0)$, and $K_i=K(R,
Z ,\zeta_i,\D_i)$.

Let $\Pi$ be a hyperbolic plane containing $Z$, $\zeta_1$ and
$\zeta_2$.  We set $s = S(R,Z)\cap \Pi$, $D=\calb\cap\Pi$,
$\lambda_i=\Pi_i\cap\Pi$, and $k_i = K_i \cap \Pi$.
We have $K_1\subset K_2$ if and only if $k_1\subset k_2$, and we have
$K_1\cap K_2=\emptyset$ if and only if $k_1\cap k_2=\emptyset$.

Let $A_i$ and $B_i$ denote the two points where $\lambda_i$ meets $s$.
For $i = 1,2$ let $P_i$ denote the point $\eta_i\cap\lambda_i$.  Using
the right triangle $A_iP_iZ$ we find that $\cos\angle A_iZ\zeta_i =
(\tanh ZP_i)/(\tanh ZA_i)=(\tanh w_i)/(\tanh R)$, so that
\Equation\label{not so duh} \Psi_i = \angle{ A_iZ\zeta_i} = \angle{
  \zeta_iZB_i}=\frac12 \angle{ A_iZB_i}.
\EndEquation

Since $0<\Psi_2 \le \Psi_1<\pi/2$, there is a unique
arc $\tau_i\subset s$ which has endpoints $A_i$ and $B_i$ 
and has length less than half the length of $s$.

Let $E:t\mapsto e^{2\pi t}$ denote the standard covering map $:\RR\to
S^1$, and let $p:\RR\to s$ denote the composition of $E$ with a
Euclidean similarity transformation of $S^1$ onto $s$. Since the
spherical distance between $\zeta_1$ and $\zeta_2$ is $\alpha$,
we may choose the similarity transformation defining $p$ in
such a way that $p(0)=\zeta_1$ and $p(\alpha)=\zeta_2$.

Note that if $t_1,t_2\in\RR$ are given, and we set $T_i=p(t_i)$ for
$i=1,2$, then $\angle T_1ZT_2=|t_1-t_2|$.

For $i=1,2$ it follows from (\ref{not so duh}) that $p$ maps the
unordered pair $\{\pm\Psi_1\}$ to $\{A_1,B_1\}$, and maps the unordered
pair $\{\alpha\pm\Psi_2\}$ to $\{A_2,B_2\}$. After possibly re-labeling
the $A_i$ and $B_i$ we may assume that $p(-\Psi_1)=A_1$, $p
(\Psi_1)=B_1$, $p(\alpha-\Psi_2)=A_2$ and $p(\alpha+\Psi_2)=B_2$.

Let $J\subset\RR$ denote the smallest closed interval containing the
four numbers $\pm\Psi_1$ and $\alpha\pm\Psi_2$. Since
$0\le\Psi_i\le\pi/2$ and $0\le\alpha\le\pi$, we have
$J\subset[\pi/2,\pi]$. Hence $p$ maps $J$ homeomorphically onto an arc
in $s$. Since the intervals $[-\Psi_1,\Psi_1]$ and
$[\alpha-\Psi_2,\alpha+\Psi_2]$ have length $<\pi$, the map $p$ sends
$[-\Psi_1,\Psi_1]$ and $[\alpha-\Psi_2,\Psi_2]$) onto $\tau_1$ and
$\tau_2$ respectively.

Consider the case in which $\alpha \le \Psi_1 - \Psi_2$. In this case
we have $-\Psi_1\le\alpha-\Psi_2<\alpha+\Psi_2\le\Psi_1$, so that
$[\alpha-\Psi_2,\alpha+\Psi_2]\subset[-\Psi_1,\Psi_1]$ and hence
$\tau_2\subset \tau_1$.  We therefore have $k_1\subset k_2$ and hence
$K_1\subset K_2$, so that $\iota(R,\D_1,\D_2,\alpha) =\kappa(R,\D_2)$.
This proves (1).

Next suppose that $\alpha > \Psi_1 + \Psi_2$.  In this case we have
$-\Psi_1<\Psi_1<\alpha-\Psi_2<\alpha+\Psi_2$, so that
$[\alpha-\Psi_2,\alpha+\Psi_2]\cap[-\Psi_1,\Psi_1]=\emptyset$ and
hence $\tau_2\cap \tau_1=\emptyset$.  We therefore have $k_1\cap
k_2=\emptyset$ and hence $K_1\cap K_2=\emptyset$, so that
$\iota(R,\D_1,\D_2,\alpha) =0$.  This proves (2).

We now turn to the case $\Psi_1 - \Psi_2 < \alpha \le \Psi_1 +
\Psi_2$.  (In particular $\alpha$ is then non-zero.)  Since
$\Psi_1\ge\Psi_2$, in this case we have
$-\Psi_1<\alpha-\Psi_2<\Psi_1<\alpha+\Psi_2$, so that
$[\alpha-\Psi_2,\alpha+\Psi_2]$ and $[-\Psi_1,\Psi_1]$ overlap in the
common sub-interval $[\alpha-\Psi_2,\Psi_1]$ which is proper in both
of them. Hence the arcs $\tau_1$ and $\tau_2$ overlap in a common
sub-arc $\tau$ which has endpoints $A_2$ and $B_1$ and is proper in
both the $\tau_i$. It follows that the lines $\lambda_1$ and
$\lambda_2$ meet at a point $Y$ lying in the disk $D$, and that the
ray originating at $Z$ and passing through $Y$ meets the arc $\tau$ at
some point $X$.  We may write $X=p(\alpha_1)$ for some
$\alpha_1\in(\alpha-\Psi_2,\Psi_1)$.  We set
$\alpha_2=\alpha-\alpha_1\in[\alpha-\Psi_1,\Psi_2]$. In particular we
have $\alpha_1,\alpha_2\in(-\pi/2,\pi/2)$, and (\ref{truck talk})
obviously holds with our choice of the $\alpha_i$.
The cases where $\alpha_2>0$ and where $\alpha_2<0$ are illustrated in
the diagram.

\begin{figure}[h]
\begin{picture}(0,0)%
\includegraphics{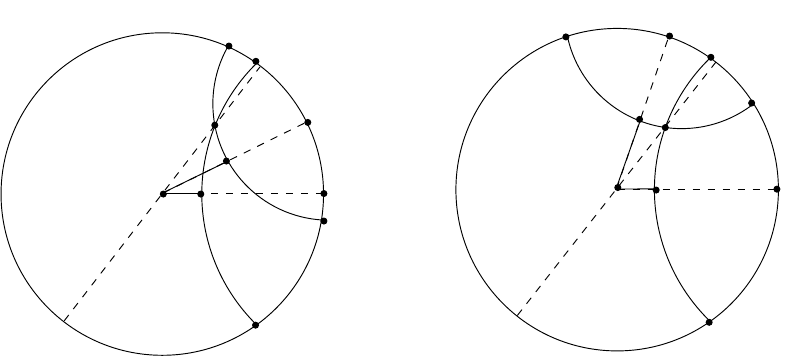}%
\end{picture}%
\setlength{\unitlength}{3947sp}%
\begingroup\makeatletter\ifx\SetFigFont\undefined%
\gdef\SetFigFont#1#2#3#4#5{%
  \reset@font\fontsize{#1}{#2pt}%
  \fontfamily{#3}\fontseries{#4}\fontshape{#5}%
  \selectfont}%
\fi\endgroup%
\begin{picture}(6296,2840)(593,-2567)
\put(5874,102){\makebox(0,0)[lb]{\smash{{\SetFigFont{12}{13.2}{\familydefault}{\mddefault}{\updefault}{\color[rgb]{0,0,0}$\zeta_2$}%
}}}}
\put(1647,-1315){\makebox(0,0)[lb]{\smash{{\SetFigFont{12}{13.2}{\familydefault}{\mddefault}{\updefault}{\color[rgb]{0,0,0}$Z$}%
}}}}
\put(5254,-1264){\makebox(0,0)[lb]{\smash{{\SetFigFont{12}{13.2}{\familydefault}{\mddefault}{\updefault}{\color[rgb]{0,0,0}$Z$}%
}}}}
\put(2330,  4){\makebox(0,0)[lb]{\smash{{\SetFigFont{12}{13.2}{\familydefault}{\mddefault}{\updefault}{\color[rgb]{0,0,0}$B_2$}%
}}}}
\put(2680,-233){\makebox(0,0)[lb]{\smash{{\SetFigFont{12}{13.2}{\familydefault}{\mddefault}{\updefault}{\color[rgb]{0,0,0}$B_1$}%
}}}}
\put(2485,-1113){\makebox(0,0)[lb]{\smash{{\SetFigFont{12}{13.2}{\familydefault}{\mddefault}{\updefault}{\color[rgb]{0,0,0}$P_2$}%
}}}}
\put(2251,-1439){\makebox(0,0)[lb]{\smash{{\SetFigFont{12}{13.2}{\familydefault}{\mddefault}{\updefault}{\color[rgb]{0,0,0}$P_1$}%
}}}}
\put(3216,-1552){\makebox(0,0)[lb]{\smash{{\SetFigFont{12}{13.2}{\familydefault}{\mddefault}{\updefault}{\color[rgb]{0,0,0}$A_2$}%
}}}}
\put(5887,-1410){\makebox(0,0)[lb]{\smash{{\SetFigFont{12}{13.2}{\familydefault}{\mddefault}{\updefault}{\color[rgb]{0,0,0}$P_1$}%
}}}}
\put(5917,-934){\makebox(0,0)[lb]{\smash{{\SetFigFont{12}{13.2}{\familydefault}{\mddefault}{\updefault}{\color[rgb]{0,0,0}$Y$}%
}}}}
\put(2056,-779){\makebox(0,0)[lb]{\smash{{\SetFigFont{12}{13.2}{\familydefault}{\mddefault}{\updefault}{\color[rgb]{0,0,0}$Y$}%
}}}}
\put(6300,-2443){\makebox(0,0)[lb]{\smash{{\SetFigFont{12}{14.4}{\familydefault}{\mddefault}{\updefault}{\color[rgb]{0,0,0}$A_1$}%
}}}}
\put(2653,-2487){\makebox(0,0)[lb]{\smash{{\SetFigFont{12}{14.4}{\familydefault}{\mddefault}{\updefault}{\color[rgb]{0,0,0}$A_1$}%
}}}}
\put(4944, 98){\makebox(0,0)[lb]{\smash{{\SetFigFont{12}{13.2}{\familydefault}{\mddefault}{\updefault}{\color[rgb]{0,0,0}$B_2$}%
}}}}
\put(6874,-1287){\makebox(0,0)[lb]{\smash{{\SetFigFont{12}{13.2}{\familydefault}{\mddefault}{\updefault}{\color[rgb]{0,0,0}$\zeta_1$}%
}}}}
\put(5432,-818){\makebox(0,0)[lb]{\smash{{\SetFigFont{12}{13.2}{\familydefault}{\mddefault}{\updefault}{\color[rgb]{0,0,0}$P_2$}%
}}}}
\put(6310,-155){\makebox(0,0)[lb]{\smash{{\SetFigFont{12}{13.2}{\familydefault}{\mddefault}{\updefault}{\color[rgb]{0,0,0}$B_1$}%
}}}}
\put(6677,-554){\makebox(0,0)[lb]{\smash{{\SetFigFont{12}{13.2}{\familydefault}{\mddefault}{\updefault}{\color[rgb]{0,0,0}$A_2$}%
}}}}
\put(3100,-681){\makebox(0,0)[lb]{\smash{{\SetFigFont{12}{13.2}{\familydefault}{\mddefault}{\updefault}{\color[rgb]{0,0,0}$\zeta_2$}%
}}}}
\put(3244,-1322){\makebox(0,0)[lb]{\smash{{\SetFigFont{12}{13.2}{\familydefault}{\mddefault}{\updefault}{\color[rgb]{0,0,0}$\zeta_1$}%
}}}}
\end{picture}%
\end{figure}

Since $p(\alpha_1)=X$ and $p(0)=\zeta_1$, we have
$\angle P_1ZY=\angle\zeta_1ZX=|\alpha_1|$. Since $p(\alpha_1)=X$
and $p(\alpha)=\zeta_2$, we have
$\angle P_2ZY=\angle\zeta_2ZX=|\alpha-\alpha_1|=|\alpha_2|$. 
Thus
for $i=1,2$ we have 
\Equation\label{enfin bref}
\angle  P_iZY=|\alpha_i|.
\EndEquation

For $i = 1,2$, the hyperbolic line
segment $\overline{ZY}$ is the common hypotenuse of the two
right triangles with vertices at $Z$, $Y$ and $P_i$. Hence for $i=1,2$,
the hyperbolic tangent of the length of this segment is 
$(\tanh \D_i)/\cos (\angle  P_iZY)$, which by (\ref{enfin bref})
is equal to
$(\tanh \D_i)/{\cos \alpha_i} $. In particular we have
$ {\tanh \D_1}\,{\cos \alpha_2} = {\tanh
  \D_2}\,{\cos \alpha_1} $.
This is (\ref{i hold that on the seas}).

For $i=1,2$, since $\Pi_i$ meets $\Pi$ perpendicularly in the line
$\ell_i$, we have $\Pi_1\cap\Pi_2=L$, where $L$ denotes the line
meeting $\Pi$ perpendicularly at $Y$.  Let $\ell_0$ denote the
hyperbolic line containing $Z$, $Y$ and $X$, and let $\Pi_0s$ denote
the plane which meets $\Pi$ perpendicularly in the line $\ell_0$.
Since $X\in\tau$, 
and since $\tau$ subtends an angle $\Psi_1+\Psi_2-\alpha<2\Psi_2<\pi$,
the endpoints $B_1$ and $A_2$ of $\tau$ lie in different components of
$\HH^3-\Pi_0$. We index the two half-spaces  bounded by
$\Pi_0$ as $H_0^1$ and $H_0^2$ in such a way that $B_1\in H_0^1$
 and $A_2\in H_0^2$.
We set $K_0^i=\calb\cap H_0^i$ for $i=1,2$.

Let $\ell$ denote the line which is perpendicular to $\Pi$ at $Y$.
Since the lines $\ell_0$, $\ell_1$ and $\ell_2$ meet at the point $Y$,
the planes $\Pi_0$, $\Pi_1$ and $\Pi_2$ meet in the line $\ell$.  Let
$h_0$, $h_1$ and $h_2$ denote the half-planes of $\Pi_0$, $\Pi_1$ and
$\Pi_2$ which are bounded by $\ell$ and contain $X$, $B_1$ and $A_2$
respectively. The definition of the half-spaces $H_0^1$ and $H_0^2$
implies that $H_0^1\cap H_0^2$ has frontier $h_1\cup h_2$ and contains
$h_0$. It follows that $H_1\cap H_2 = (H_1\cap H_0^1) \cup (H_2\cap
H_0^2).$ In particular we have
$$K_1\cap K_2 = (K_1\cap K_0^1) \cup (K_2\cap
K_0^2).$$
Since $(K_1\cap K_0^1) \cap (K_2\cap K_0^2)\subset\Pi_0$, it follows that
\Equation\label{why oh why oh}
\iota(R,\D_1,\D_2,\alpha)=
\vol (K_1\cap K_2)=
\vol(K_i\cap K_0^1)+
\vol(K_i\cap K_0^2).
\EndEquation

Now set $\zeta_0^1=p(\alpha_1+\pi/2)$ and 
$\zeta_0^2=p(\alpha_1-\pi/2)$. For $i=1,2$, 
set $\eta_0^i=\eta_{\zeta_0^i}$. 
Since $p(0)=\zeta_1$ and $p(\alpha_1+\pi/2)=\zeta_0^1$, we have
$\angle\zeta_1Z\zeta_0^1=|(\alpha_1+\pi/2)-0|=\alpha_1+\pi/2$. Since
$p(\alpha)=\zeta_2$ and $p(\alpha_1-\pi/2)=\zeta_0^2$, we have
$\angle\zeta_2Z\zeta_0^2=|(\alpha_1-\pi/2)-\alpha|=\alpha_2+\pi/2$. Thus
for $i=1,2$ we have
\Equation\label{all right, you name them}
\angle\zeta_iZ\zeta_0^i=\alpha_i+\pi/2.
\EndEquation

Since $p(\alpha_1)=X$, $p(\alpha_1+\pi/2)=\zeta_0^1$, and
$p(\alpha_1-\pi/2)=\zeta_0^2$, the rays 
$\eta_0^1$ and $\eta_0^2$ are perpendicular to $\ell_0$ and hence to
$\Pi_0$, and they point in opposite directions in the line through $Z$
which is perpendicular to $\Pi_0$. Hence each of the half-spaces 
$H(Z,\zeta_0^i,0)$ is bounded by $\Pi_0$. It follows that
$H(Z,\zeta_0^1,0)$ and $H(Z,\zeta_0^2,0)$ are equal to $H_0^1$ and $H_0^2$
in some order. But since $p(\alpha_1+\pi/2)=\zeta_0^1$, $p(\Psi_1)=B_1$ and
$\alpha_1<\Psi_1$, we have 
$\angle\zeta_0^1ZB_1=|\alpha_1+\pi/2-\Psi| =
\alpha_1-\Psi+\pi/2>\pi/2$. It follows that 
$H(Z,\zeta_0^1,0)=H_0^1$, and therefore that $H(Z,\zeta_0^2,0)=H_0^2$. Hence
\Equation\label{mobile mystics}K_0^i=K(R,Z,\zeta_0^i,0) \text{ for }i=1,2.
\EndEquation

From (\ref{all right, you name them}),
(\ref{mobile mystics}), and the definition of 
$\iota$, it follows that
\Equation\label{let alone}
\vol(K_i\cap K_0^i)=\iota(R,\D_i,0,\alpha_i +
\pi/2).
\EndEquation
The equality (\ref{that great street}) follows from (\ref{why oh why
  oh}) and (\ref{let alone}).

It remains to show that $\alpha_1,\alpha_2\in[-\pi/2,\pi/2]$ are uniquely
determined by the conditions \ref{truck talk} and \ref{i hold that on
  the seas}.
For this purpose it suffices to show that there is at most one number
$x\in[-\pi/2,\pi/2]$ such that
 ${\tanh \D_1}{\cos x} = {\tanh \D_2}\cos (\alpha-x)$, i.e. such that
\Equation\label{the meaning of romance}
{\tanh \D_1}\,
{\cos x}=
(\tanh \D_2)(\cos \alpha\cos x+\sin \alpha\sin x.)
\EndEquation

Since $\alpha$ and $\D_2$  are non-zero, the equality
(\ref{the meaning of romance}) cannot hold with $x=\pm\pi/2$. If
$\pi/2 <x<\pi/2$ then (\ref{the meaning of romance}) is equivalent to
$(\tanh\D_2)(\sin\alpha\,\tan x+\cos\alpha)=\tanh\D_1$, which can have
at most one solution for $x\in(-\pi/2,\pi/2)$ since $\tan x$ increases
monotonically on that interval.
\EndProof

\Number\label {how we did iota} The results of this appendix can be
used to compute arbitrary values of the functions $\kappa$ and
$\iota$. For any $R>0$ it is clear that $\kappa(R,\D)=0$ for $\D\ge R$
and that $\kappa(R,0)=B(R)/2=(\pi/2)(\sinh(2R)-2R$ (see \ref{bees in
  my bonnet, pain in my heart}).  When $0<\D<R$, we can calculate
$\kappa(R,\D)$ directly from Proposition \ref{formula for kappa}.

Similarly,  it is
clear that for any $R>0$ and any $\alpha\in[0,\pi]$, we have
$\iota(R,\D_1,\D_2,\alpha)=0$ whenever either $\D_1$ or $\D_2$ is at least
$R$.  If each of the $\D_i$ is less than  $R$ then Proposition
\ref{formula for iota} reduces the calculation of
$\iota(R,\D_1,\D_2,\alpha)$ to the special case in which $\D_2=0$. It
is clear that $\iota(R,\D,0,0)=\kappa(R,\D)$ and that
$\iota(R,\D,0,\pi)=0$, and it follows from Proposition
\ref{complements} that $\iota(R,\D,0,\pi/2)=\kappa(R,\D)/2$.
If $\pi/2<\alpha<\pi$, we can compute
$\iota(R,\D,0,\alpha)$ directly from Proposition \ref{special
  intersection}. If 
$0<\alpha<\pi/2$, Proposition \ref{complements} reduces the
calculation of $\iota(R,\D_1,\D_2,\alpha)$ to the calculation of
$\iota(R,\D,\pi-\alpha)$ and $\kappa(R,\D)$, which can be carried out
using Propositions \ref{formula for iota} and \ref{formula for kappa}.

There are two points in the body of this paper where the formulas
given in this appendix were used in rigorous sampling arguments to
establish numerical bounds for certain functions.  The proof of Lemma
\ref{if the parents ate the spinach} required calculating $100$
different numerical values of $\iota(R,\D,\D',\alpha)$, and the proof
of Lemma \ref{brighton beach} required calculating several thousand
values of $\kappa$.

The evaluation of $\kappa$ from Proposition \ref{formula for kappa} is
straightforward as it is given by a closed-form expression.  In
evaluating $\iota$, there are two steps which require somewhat more
elaborate numerical methods.  First, the application of Proposition
\ref{formula for iota} in the case where $\Psi_1 - \Psi_2 < \alpha \le
\Psi_1 + \Psi_2$ requires the numerical solution of the equations
(\ref{truck talk}) and (\ref{i hold that on the seas}).  We used the
routine {\tt hybrd}, which is included in the MINPACK library, to find
approximate solutions to this system.  The step involving Proposition
\ref{special intersection} seems to require numerical integration, as
we do not know of a closed-form expression for the integral in
(\ref{we bombed in new haven}).  For this purpose we relied on the
adaptive Gaussian quadrature method which is implemented in the python
scipy module and uses the FORTRAN quadrature routines from the
QUADPACK library.
\EndNumber


\begin{thebibliography}{10}

\bibitem{agol}
Ian Agol.
\newblock Tameness of hyperbolic 3-manifolds.
\newblock arXiv:math.GT/0405568.

\bibitem{last}
Ian Agol, Marc Culler, and Peter~B. Shalen.
\newblock Singular surfaces, mod 2 homology, and hyperbolic volume, {I}.
\newblock arXiv:math.GT/0506396, to appear in {\it Trans. Amer. Math. Soc.}

\bibitem{rankfour}
Ian Agol, Marc Culler, and Peter~B. Shalen.
\newblock Dehn surgery, homology and hyperbolic volume.
\newblock {\em Algebr. Geom. Topol.}, 6:2297--2312, 2006.

\bibitem{ast}
Ian Agol, Peter~A. Storm, and William~P. Thurston.
\newblock Lower bounds on volumes of hyperbolic {H}aken 3-manifolds.
\newblock {\em J. Amer. Math. Soc.}, 20(4):1053--1077 (electronic), 2007.
\newblock With an appendix by Nathan Dunfield.

\bibitem{accs}
James~W. Anderson, Richard~D. Canary, Marc Culler, and Peter~B. Shalen.
\newblock Free {K}leinian groups and volumes of hyperbolic {$3$}-manifolds.
\newblock {\em J. Differential Geom.}, 43(4):738--782, 1996.

\bibitem{bjorner}
Anders Bj{\"o}rner.
\newblock Nerves, fibers and homotopy groups.
\newblock {\em J. Combin. Theory Ser. A}, 102(1):88--93, 2003.

\bibitem{boroczky}
K.~B{\"o}r{\"o}czky.
\newblock Packing of spheres in spaces of constant curvature.
\newblock {\em Acta Math. Acad. Sci. Hungar.}, 32(3-4):243--261, 1978.

\bibitem{burns}
Robert~G. Burns.
\newblock On the intersection of finitely generated subgroups of a free group.
\newblock {\em Math. Z.}, 119:121--130, 1971.

\bibitem{cg}
Danny Calegari and David Gabai.
\newblock Shrinkwrapping and the taming of hyperbolic 3-manifolds.
\newblock {\em J. Amer. Math. Soc.}, 19(2):385--446 (electronic), 2006.

\bibitem{lastplusone}
Marc Culler and Peter~B. Shalen.
\newblock Singular surfaces, mod 2 homology, and hyperbolic volume, {II}.
\newblock math.GT/0701666.

\bibitem{cusp}
Marc Culler and Peter~B. Shalen.
\newblock Volume and homology of one-cusped hyperbolic $3$-manifolds.
\newblock http://www.math.uic.edu/~shalen/cusp1.pdf.

\bibitem{splittings}
Marc Culler and Peter~B. Shalen.
\newblock Varieties of group representations and splittings of {$3$}-manifolds.
\newblock {\em Ann. of Math. (2)}, 117(1):109--146, 1983.

\bibitem{paradoxical}
Marc Culler and Peter~B. Shalen.
\newblock Paradoxical decompositions, {$2$}-generator {K}leinian groups, and
  volumes of hyperbolic {$3$}-manifolds.
\newblock {\em J. Amer. Math. Soc.}, 5(2):231--288, 1992.

\bibitem{log5}
Marc Culler and Peter~B. Shalen.
\newblock Hyperbolic volume and mod {$p$} homology.
\newblock {\em Comment. Math. Helv.}, 68(3):494--509, 1993.

\bibitem{betti2}
Marc Culler and Peter~B. Shalen.
\newblock The volume of a hyperbolic {$3$}-manifold with {B}etti number {$2$}.
\newblock {\em Proc. Amer. Math. Soc.}, 120(4):1281--1288, 1994.

\bibitem{Fe}
Werner Fenchel.
\newblock {\em Elementary geometry in hyperbolic space}, volume~11 of {\em de
  Gruyter Studies in Mathematics}.
\newblock Walter de Gruyter \& Co., Berlin, 1989.
\newblock With an editorial by Heinz Bauer.

\bibitem{JS}
William~H. Jaco and Peter~B. Shalen.
\newblock Seifert fibered spaces in {$3$}-manifolds.
\newblock {\em Mem. Amer. Math. Soc.}, 21(220):viii+192, 1979.

\bibitem{kent}
Richard~P. Kent, IV.
\newblock Intersections and joins of free groups.
\newblock arXiv:0802.0033.

\bibitem{kurosh}
A.~G. Kurosh.
\newblock {\em The theory of groups}.
\newblock Chelsea Publishing Co., New York, 1960.
\newblock Translated from the Russian and edited by K. A. Hirsch. 2nd English
  ed. 2 volumes.

\bibitem{l-mcr}
Larsen Louder and D.~B. McReynolds.
\newblock Graphs of subgroups of free groups.
\newblock {\em Algebr. Geom. Topol.}, 9(1):327--335, 2009.

\bibitem{prez}
Andrew Przeworski.
\newblock A universal upper bound on density of tube packings in hyperbolic
  space.
\newblock {\em J. Differential Geom.}, 72(1):113--127, 2006.

\bibitem{spanier}
Edwin~H. Spanier.
\newblock {\em Algebraic topology}.
\newblock Springer-Verlag, New York, 1981.
\newblock Corrected reprint.

\bibitem{thurstonnotes}
William~P. Thurston.
\newblock The geometry and topology of $3$-manifolds.
\newblock Lecture notes, 1978.

\bibitem{Tret}
Marvin~D. Tretkoff.
\newblock A topological approach to the theory of groups acting on trees.
\newblock {\em J. Pure Appl. Algebra}, 16(3):323--333, 1980.

\end{thebibliography}
\end{document}